\documentclass[preprint,12pt,authoryear]{elsarticle}

\usepackage{epsfig, graphicx, subfigure, color}
  \usepackage[margin=1in]{geometry}
\usepackage{amssymb,amsmath}
 \usepackage{url}
 \usepackage{setspace} 
 \usepackage{subcaption}
 \usepackage{graphicx}
 \usepackage{makecell}
 \usepackage{float}
\usepackage{bm}
\usepackage{txfonts}
 \usepackage{multicol}
 \usepackage{multirow}
 \usepackage{array}
  \usepackage{subfigure}
  \usepackage{caption}
\usepackage{lineno}
\usepackage{booktabs}
\usepackage{placeins}
  
\graphicspath{{./art/}}

\def\E{\mathbb{E}}

\def\l{\langle}

\journal{Environmetrics}
\newtheorem{prop}{Proposition}
\newtheorem{example}{Example}
\newtheorem{cor}{Corollary}
 \newtheorem{thm}{Theorem}
 \newtheorem{definition}{Definition}
 \newtheorem{lem}[thm]{Lemma}
 \newdefinition{rem}{Remark}
 \newproof{proof}{Proof}
 \newproof{pot}{Proof of Theorem \ref{thm2}}

\newtheorem{suppthm}{Theorem}[section]
\newtheorem{supplem}[suppthm]{Lemma}
\newtheorem{suppprop}[suppthm]{Proposition}
\newtheorem{suppcor}[suppthm]{Corollary}
\newtheorem{suppdefinition}[suppthm]{Definition}
\newtheorem{suppexample}[suppthm]{Example}
\newtheorem{supprem}[suppthm]{Remark}

\begin{document}
\begin{frontmatter}



\title{A new family of Gaussian processes for modeling animal movement: application to bat telemetry data }
\author[label1]{J.H. Ram\'{\i}rez-Gonz\'{a}lez}
\author[label2]{A. Murillo-Salas}
\author[label1]{Ying Sun}

\ead{josehermenegildo.ramirezgonzalez@kaust.edu.sa, amurillos@ugto.mx, ying.sun@kaust.edu.sa}

\affiliation[label1]{organization={CEMSE Division, Statistics Program, King Abdullah University of Science and Technology},
            city={Thuwal},
            postcode={23955-6900},
            state={Makkah},
            country={Saudi Arabia}}

\affiliation[label2]{organization={Universidad de Guanajuato},
            city={Guanajuato},
            postcode={36240},
            state={Guanajuato},
            country={Mexico}}





\begin{abstract}
Modeling animal movement is essential for addressing various ecological and biological questions. However, developing an effective predictive model for animal movement is a challenging task. In this paper, we introduce a new family of Gaussian processes, derived from the limiting fluctuations of the rescaled occupation-time process of certain branching particle systems, and study its applicability to real animal movement data. We examine two subfamilies and show that these processes exhibit long-range dependence and covariance functions with logarithmic asymptotic growth. For the exponential subfamily used in the applied analysis, the process is also non-stationary and not intrinsically stationary on compact time intervals. These properties are relevant when dealing with animal trajectories that exhibit strong memory. Finally, we illustrate the practical applicability of the proposed model by analyzing bat movement data.
\end{abstract}


\end{frontmatter}





\doublespacing
\section{Introduction}\label{introduction}

Numerous studies have documented that migratory species are increasingly threatened by habitat loss, climate change, and anthropogenic pressures (\cite{Harris2009}, \cite{Williams2021}, \cite{LezamaOchoa2024}). Consequently, understanding animal migration is crucial for designing protected areas and effective conservation strategies (\cite{Runge2014}). In this context, Gaussian processes have become useful tools in probabilistic modeling and have also been employed to describe animal movement trajectories. For example, \cite{Rieber} used treed Gaussian processes to model telemetry paths. In the present work, we are particularly interested in trajectories exhibiting strong memory-dependent movement. Bats provide a representative example of such behavior (\cite{Sarel2017}), and, more generally, memory-based formulations have been recognized as an important component of stochastic models for animal movement; see, for instance, \cite{Smouseetal}.

Long-range dependence (LRD) and long memory have been extensively studied because of their relevance in a wide variety of disciplines (\cite{samo}). Their importance is reflected in numerous contributions, with applications ranging from fractional random fields (\cite{Anh1999}) and Internet traffic modeling (\cite{Karagiannis2004}) to DNA sequencing (\cite{Karmeshu2004}), econometrics (\cite{Guegan2005}), climate studies (\cite{Varotsos2006}), and linguistics (\cite{Alvarez2006}). More recently, in animal telemetry, LRD models have been used to represent movement patterns that exhibit dependence beyond the short-memory framework, thereby allowing greater flexibility in trajectory modeling (\cite{JY}, \cite{RCMS24}). These works address, among other topics, the detection of long memory in data, parameter estimation, limit theorems under LRD, and the simulation of long-memory processes.

A further aspect relevant to the present article is the large-time behavior of the covariance
itself. In the covariance literature, the asymptotic
regime of covariance functions has long been recognized as a central feature
in the modeling of dependence over large temporal or spatial scales
(\cite{GayHeyde1990}, \cite{Gneiting2000}, \cite{HaslettRaftery1989}).
In particular, Ma and Bhadra (\cite{MaBhadra2023}) show that covariance
models with exponentially decaying tails may be inadequate when distant
observations remain substantially informative, thereby motivating covariance
classes with polynomial tail behavior. In the non-stationary setting considered below, the covariance between a fixed
past position and a remote future position grows rather than decays;
accordingly, we refer to this property as asymptotic covariance growth rather
than covariance tail decay. Thus, in applications where remote past positions
continue to influence future movement, the asymptotic growth regime of the
covariance is not merely a technical aspect of the model, but may substantially
affect how dependence is represented across large temporal scales.

Motivated by these considerations, and also by covariance structures arising in the study of occupation-time fluctuations of branching particle systems (\cite{BGT-SPL2004}, \cite{BGT2007}, \cite{BGT2008}, \cite{LMR2024}), we consider a class of centered non-stationary Gaussian processes \(\zeta=(\zeta_t)_{t\geq 0}\) with covariance function
\small
\begin{eqnarray}\label{critical-wegithed-covariance}
  &&\lefteqn{K_f(s,t):=}\\
  &&2\int_{0}^{s\wedge t}f(r)\left[(s+t-2r)\log(s+t-2r)-(s-r)\log(s-r)\right.\nonumber\\
  &&\hspace{3.6cm}\left.-(t-r)\log(t-r)\right]dr,\nonumber
\end{eqnarray}
\normalsize
where \(f:\mathbb{R}_{+}\to \mathbb{R}_{+}\) is a measurable and locally integrable function. This formulation defines a broad family of Gaussian processes in which the weighting function \(f\) determines the covariance structure. The objective of this article is to study the applicability of this family to animal telemetry data, particularly in situations where the observed trajectories exhibit evidence of long-term dependence, as has been reported for bats (\cite{Conenna}, \cite{Voigt2017}).

As shown in Section~\ref{theory}, the covariance of the proposed Gaussian process with remote
future positions grows at logarithmic order. Since the logarithm is a slowly varying function (\cite{BinghamGoldieTeugels1987}), this yields an alternative LRD structure in which the contribution of the remote past is retained, but in a milder way than in models whose covariance growth is polynomial.

This distinction is particularly relevant in applications where past trajectory values remain influential over long time horizons, but polynomial covariance growth may overstate the contribution of distant past positions.


The application considered in this paper is also motivated by previous
models for telemetry data based on the integral of an
Ornstein--Uhlenbeck process. The integrated Ornstein--Uhlenbeck model driven by Brownian motion was considered in \cite{John}, whereas
\cite{JY} and \cite{RCMS24} consider the fractional extension
$\mu_H=(\mu_H(t))_{t\geq 0}$, where the velocity process
$v_H=(v_H(t))_{t\geq 0}$ is defined through a stochastic differential
equation driven by a fractional Brownian motion
$W_H=(W_H(t))_{t\geq 0}$ with Hurst parameter
$H\in(0,1)\setminus\{\frac{1}{2}\}$; see
(\ref{velocity-sde}) and (\ref{position-process}) below.
By Proposition~2.7 of \cite{RCMS24}, the process $\mu_H$ exhibits LRD.

A particularly relevant subclass of the family
\eqref{critical-wegithed-covariance} is obtained by taking
\[
f_{\sigma^2,\beta}(u)=\sigma^2e^{-\beta u},
\qquad \sigma^2>0,\quad \beta\ge0.
\]
This choice is motivated by the covariance representation of the process $\mu_H$.
As discussed in Section~\ref{Preliminaries}, the parameters $\sigma^2$
and $\beta$ retain, respectively, the covariance-scale and
exponential-attenuation roles that arise in the covariance representation
of $\mu_H$.

The theoretical comparison is complemented by the application to bat
telemetry data presented in Section~\ref{Application}. As documented in
Appendix~C, the telemetry series exhibit DFA-based evidence of
dependence over large temporal scales, together with evidence against weak
stationarity and mixed evidence concerning intrinsic stationarity. These
empirical features motivate the consideration of non-stationary position
models while retaining both stationary-increment and non-stationary-increment
specifications in the model comparison.

Accordingly, we fit the process $\zeta$ to the bat telemetry trajectories
as a non-stationary LRD model with logarithmic-order asymptotic covariance
growth. The code implementing the inferential procedure is available at
\url{https://github.com/joseramirezgonzalez/f-wsfBm}. To facilitate the
exploration of the proposed stochastic models, we also developed an
interactive R Shiny application, available at
\url{https://jhramgon.shinyapps.io/aplication_r_log_w/}.

Finally, Section~\ref{Discu} presents the main conclusions concerning the
proposed model and methodology, together with their potential
generalizations, advantages, and limitations.

\section{The model: some theoretical results}\label{theory}

In this section, we present theoretical results for the Gaussian process $\zeta=(\zeta_t)_{t\geq 0}$ with covariance function $K_f$ given in \eqref{critical-wegithed-covariance}. These results include sufficient conditions ensuring that the kernel $K_f$ is positive definite, as well as some memory properties of the process. In particular, we show that the covariance $\mathbb{E}(\zeta_r\zeta_{s+T})$ exhibits asymptotic logarithmic growth; see Proposition~\ref{LRM}. Proofs of the results stated in this section, together with additional theoretical results on sample path properties of $\zeta$, including total variation, quadratic variation, modulus of continuity, Hölder regularity, and the non-semimartingale property, are provided in Appendix~B. Appendix~B also contains additional results on non-stationarity and, for weight functions satisfying \(1\asymp f\) on compact intervals, lack of intrinsic stationarity, as well as further results on self-similarity and short-time asymptotics.

For clarity, throughout this article, stationarity refers to weak
stationarity. A second-order process \(X=(X_t)_{t\geq0}\) is weakly
stationary if \(\mathbb{E}(X_t)\) is constant and
\(\operatorname{Cov}(X_s,X_t)\) depends only on the lag \(t-s\).
It is intrinsically stationary if, for every admissible lag \(h\),
the mean and variance of \(X_{t+h}-X_t\) do not depend on \(t\).

The following result provides a sufficient condition on $f$ for \eqref{critical-wegithed-covariance} to be a covariance function.

\begin{thm} \label{Main-th1}
    Let $f:\mathbb{R}_{+}\to \mathbb{R}_{+}$ be a measurable function such that, for any $\delta>0$,
    \begin{equation}\label{Integrability-condition}
        \int_0^\delta f(u)\, du<\infty .
    \end{equation}
    Then, $K_f(s,t)$ given by \eqref{critical-wegithed-covariance} is a covariance function.
\end{thm}




\begin{rem}
(i) The covariance kernel \eqref{critical-wegithed-covariance}, with the choice $f(r)\equiv 1$, appeared in \cite{BGT2007,BGT2008}.

(ii) For the family $\mathcal{C}_1:=\{f_{\alpha}(u)=u^{\alpha}: u>0,\ \alpha>-1\}$, the covariance function is
\begin{equation}\label{Cf-polinomial}
C_1(s,t):=2\int_0^{s\wedge t} u^\alpha \Big[(s+t-2u)\log(s+t-2u)-(s-u)\log(s-u)-(t-u)\log(t-u)\Big]\,du.
\end{equation}
Moreover, for $-1<\alpha<0$, the covariance function \eqref{Cf-polinomial} was obtained in the study of rescaled occupation-time fluctuations of a $(d,a,1,\gamma)$-branching particle system with $d=a$ and $\alpha=\gamma-1$; see Theorem~2.1(ii) in \cite{LMR2024}.

(iii) For the family $\mathcal{C}_2:=\{f_{\alpha}(u)=e^{\alpha u}: u>0,\ \alpha\in\mathbb{R}\}$, the covariance function is
\begin{equation}\label{Cf-exponential}
C_2(s,t):=2\int_0^{s\wedge t} e^{\alpha u}\Big[(s+t-2u)\log(s+t-2u)-(s-u)\log(s-u)-(t-u)\log(t-u)\Big]\,du.
\end{equation}
As will be seen in Section~\ref{Application}, the family $\mathcal{C}_2$ is used to model bat telemetry data.
\end{rem}

When dealing with LRD, the notion has no unique definition
\citet{samo}. In this article, the LRD index refers to the
normalization exponent in the increment-covariance criterion of
\citet[Theorem~2.4(5)]{BGT2007}, where the normalization is
\(T^{1-b}\). More precisely, for fixed \(0<r<\nu\) and \(0<s<t\),
a second-order stochastic process \(X=(X_t)_{t\geq0}\) is said to
exhibit LRD with index \(\kappa\) if
\[
\lim_{T\to\infty}
T^{\kappa}
\operatorname{Cov}\!\left(
X_{\nu}-X_r,\,
X_{t+T}-X_{s+T}
\right)
=
C(r,\nu,s,t),
\qquad
0<\left|C(r,\nu,s,t)\right|.
\]
Thus, the notation of \citet[Theorem~2.4(5)]{BGT2007} corresponds
to \(\kappa=1-b\). Several Gaussian processes satisfying this LRD
criterion have appeared in the study of rescaled occupation-time
fluctuations of branching particle systems; see, for example,
\cite{BGT-SPL2004,BGT2007,BGT2008,LMR2024}. The next result deals
with the memory properties of the model.
\begin{prop}\label{LRM}  The process $\zeta$ satisfies the following properties. 
\begin{itemize}
\item[(i)]  For $0<r<s$ it holds that 

\begin{equation}\label{long-memory}
    \lim_{T\rightarrow\infty}\frac{\E(\zeta_r \zeta_{s+T})}{\log(T)} =2\int_0^r f(u)(r-u) du.  
\end{equation}
\item[(ii)] (Long-range memory) For $0<r\leq \nu$ and $0<s\leq t$ it holds that 
\begin{eqnarray}\label{long-range-dependence}
  &&\lefteqn{ \lim_{T\rightarrow\infty} T\E\left[(\zeta_\nu-\zeta_r)(\zeta_{t+T}-\zeta_{s+T})\right]}\nonumber\\&&=2(t-s)\left[\int_0^\nu f(u)(\nu-u)du-\int_0^r f(u)(r-u)du\right].
\end{eqnarray}

\end{itemize}
\end{prop}

Observe that the limit in \eqref{long-memory} does not depend on \(s\). This is because \(s\) is kept fixed while \(T\to\infty\), so that the shift \(s+T\) is asymptotically equivalent to \(T\) at logarithmic scale. Consequently, the contribution of \(s\) is absorbed into lower-order terms, and the leading term is determined by the earlier time \(r\).

Whenever the right-hand side of
\eqref{long-range-dependence} is nonzero,
Proposition~\ref{LRM}(ii) shows that the process \(\zeta\)
satisfies the adopted LRD criterion with index \(\kappa=1\).
Moreover, by equation~(5) in \cite{Gorostiza2015}, the process
\(\zeta\) belongs to the long-memory regime.

\section{Application to telemetry data}\label{Application}

In this section, we apply the  process $\zeta$ to animal telemetry data. We begin with preliminaries and a review of previous telemetry models based on the integral of a fractional Ornstein--Uhlenbeck process. We then derive the likelihood function from the finite-dimensional distributions. Finally, we present the application to bat telemetry data and the comparison with several competing models.

\subsection{Preliminaries and previous work}\label{Preliminaries}

As outlined in the Introduction, previous models for telemetry data
have been based on position processes obtained by integrating
Ornstein--Uhlenbeck-type velocity processes. More precisely, the
fractional modeling framework considered in \cite{RCMS24} and
\cite{JY} describes the velocity dynamics through the stochastic
differential equation
\begin{equation}\label{velocity-sde}
    dv_H(t)
    =
    -\beta v_H(t)\,dt
    +
    \sigma\,dW_H(t),
    \qquad
    \beta,\sigma>0,
\end{equation}
where $W_H$ is a fractional Brownian motion with Hurst parameter
\[
H\in(0,1)\setminus\left\{\frac12\right\},
\]
and the associated position process is given by
\begin{equation}\label{position-process}
    \mu_H(t)
    =
    \mu_H(0)
    +
    \int_0^t v_H(s)\,ds.
\end{equation}
The process $v_H$ is defined as the path-by-path solution of
\eqref{velocity-sde}. A key feature of this model is the memory
behavior of the position process. Its covariance function is
\begin{equation}\label{QMH}
    R(s,t)
    :=
    \operatorname{Cov}\bigl(\mu_H(s),\mu_H(t)\bigr)
    =
    \sigma^2
    \int_0^t
    \int_0^s
    e^{-\beta v}
    c_H(t-v,s-u)
    e^{-\beta u}
    \,du\,dv,
    \qquad
    s,t\geq0,
\end{equation}
where
\[
    c_H(s,t)
    :=
    \frac12
    \left(
        t^{2H}
        +
        s^{2H}
        -
        |t-s|^{2H}
    \right),
    \qquad
    H\in(0,1).
\]

By Proposition~2.7 of \cite{RCMS24}, for
$H\in(0,1)\setminus\{\frac12\}$ and fixed
$0<r<\nu$ and $0<s<t$,
\begin{equation}\label{long-range-dependencemu}
\begin{aligned}
    &\lim_{T\to\infty}
    T^{2(1-H)}
    \mathbb{E}
    \Big[
        \bigl(\mu_H(\nu)-\mu_H(r)\bigr)
        \bigl(\mu_H(t+T)-\mu_H(s+T)\bigr)
    \Big]
    \\
    &\qquad
    =
    H(2H-1)(t-s)\,
    h_{\sigma^2,\beta}(\nu,r),
\end{aligned}
\end{equation}
where
\[
    h_{\sigma^2,\beta}(\nu,r)
    :=
    \frac{\sigma^2}{\beta^3}
    \left[
        e^{-\beta r}(\beta r+1)
        -
        e^{-\beta\nu}(\beta\nu+1)
    \right].
\]
Since the function
\[
    x\longmapsto e^{-\beta x}(\beta x+1)
\]
is strictly decreasing on $(0,\infty)$, one has
$h_{\sigma^2,\beta}(\nu,r)>0$ whenever $0<r<\nu$.
Consequently, for
$H\in(0,1)\setminus\{\frac12\}$, the limit in
\eqref{long-range-dependencemu} is finite and nonzero, and therefore
the process $\mu_H$ satisfies the LRD criterion adopted in this
article. Moreover, it belongs to the long-memory regime when
$H>\frac12$ and to the short-memory regime when $H<\frac12$; see
equation~(5) in \cite{Gorostiza2015}.

In the particular case $H=\frac12$, the process $W_H$ reduces to
standard Brownian motion, and the model becomes the integrated
Ornstein--Uhlenbeck process considered in \cite{John}. This Brownian
case is not covered by the preceding LRD conclusion: the covariance
between separated increments decays exponentially, and hence
$\mu_{1/2}$ does not satisfy the LRD criterion adopted in this
article.

Among the admissible weight functions covered by
Theorem~\ref{Main-th1}, we consider
\(f_{\sigma^2,\beta}(u)=\sigma^2e^{-\beta u}\), with
\(\sigma^2>0,\beta\ge0\). This choice is motivated by the covariance
structure of the position process \(\mu_H\) in \eqref{QMH}. In that
model, \(\sigma^2\) is a covariance scale parameter, whereas the factors
\(e^{-\beta u}\) and \(e^{-\beta v}\) describe the exponential
attenuation induced by the velocity dynamics. The specification
\(f_{\sigma^2,\beta}(u)=\sigma^2e^{-\beta u}\) preserves these two
parameter roles within the covariance kernel
\eqref{critical-wegithed-covariance}: \(\sigma^2\) determines the overall
covariance scale and \(\beta\) controls the exponential weighting of past
times. It therefore provides a direct two-parameter analogue of the
exponential structure appearing in \eqref{QMH}.

The exponential specification also belongs to a broader class of
admissible weights. Proposition~B.3 of Appendix~B shows that every
completely monotone and locally integrable weight can be approximated on
finite intervals by finite positive sums of exponential functions; the
associated covariance kernels converge uniformly on compact time
squares, and the corresponding Gaussian processes converge in
finite-dimensional distributions. In particular, Remark~B.4 considers
the power-law weights \(f(u)=a u^{-\gamma}\), with \(a\geq0\) and
\(0<\gamma<1\), which are completely monotone and admit a positive
mixture representation in terms of exponential weights. Equivalently,
these weights correspond to the family \(f_\alpha(u)=u^\alpha\) with
\(-1<\alpha<0\). For this range of \(\alpha\), the resulting covariance
kernel was first obtained from the rescaled occupation-time fluctuations
of a \((d,a,1,\gamma)\)-branching particle system with \(d=a\) and
\(\alpha=\gamma-1\); see Theorem~2.1(ii) of \cite{LMR2024}.

For the present application, the process $\zeta$ therefore provides a natural alternative to the model $\mu_H$. More generally, one may also consider broader families of weight functions satisfying the assumptions of Theorem~\ref{Main-th1}.

\subsection{The likelihood function}\label{subsec:likelihood}

In this subsection, we derive the likelihood function associated with the
finite-dimensional distributions of the process \(\zeta\) under the
parametric specification introduced above. For \(s,t\geq0\), its covariance
function is
\begin{equation}\label{covexp}
\begin{aligned}
\operatorname{Cov}(\zeta_s,\zeta_t)
={}&
2\sigma^2
\int_0^{s\wedge t}
e^{-\beta u}
\Big[
(s+t-2u)\log(s+t-2u) \\
&\hspace{3.2cm}
-(s-u)\log(s-u)
-(t-u)\log(t-u)
\Big]\,du .
\end{aligned}
\end{equation}

Let \(0<t_1<\cdots<t_n\) be fixed observation times, and define
\[
\bm\zeta_{\bm t}
:=
\big(\zeta_{t_1},\ldots,\zeta_{t_n}\big)^\top,
\qquad
\bm t:=(t_1,\ldots,t_n).
\]
Since \(\zeta\) is a centered Gaussian process,
\[
\bm\zeta_{\bm t}
\sim
N\!\left(
\bm 0,
\bm\Sigma_{\sigma^2,\beta,\bm t}
\right),
\]
where
\[
\bm\Sigma_{\sigma^2,\beta,\bm t}(i,j)
:=
\operatorname{Cov}(\zeta_{t_i},\zeta_{t_j}),
\qquad i,j=1,\ldots,n,
\]
and the covariance on the right-hand side is given by
\eqref{covexp}.

For an observed realization
\(\bm\zeta_{\bm t}=\bm z\), the likelihood function is therefore
\begin{equation}\label{eq:likelihood_zeta}
\begin{aligned}
L(\sigma^2,\beta\mid\bm z)
={}&
(2\pi)^{-n/2}
\left|
\bm\Sigma_{\sigma^2,\beta,\bm t}
\right|^{-1/2} \\
&\times
\exp\!\left\{
-\frac12
\bm z^\top
\bm\Sigma_{\sigma^2,\beta,\bm t}^{-1}
\bm z
\right\}.
\end{aligned}
\end{equation}

It follows from \eqref{covexp} that \(\sigma^2\) enters the covariance
matrix as a scale parameter. More precisely,
\[
\bm\Sigma_{\sigma^2,\beta,\bm t}
=
\sigma^2
\bm\Sigma_{1,\beta,\bm t},
\]
where \(\bm\Sigma_{1,\beta,\bm t}\) denotes the covariance matrix
obtained by setting \(\sigma^2=1\) in \eqref{covexp}. Consequently,
\[
\bm\Sigma_{\sigma^2,\beta,\bm t}^{-1}
=
\frac{1}{\sigma^2}
\bm\Sigma_{1,\beta,\bm t}^{-1},
\qquad
\left|
\bm\Sigma_{\sigma^2,\beta,\bm t}
\right|
=
(\sigma^2)^n
\left|
\bm\Sigma_{1,\beta,\bm t}
\right|.
\]

Hence, the log-likelihood can be written as
\begin{equation}\label{eq:loglik_scaled}
\ell(\sigma^2,\beta\mid\bm z)
=
-\frac{n}{2}\log(2\pi)
-\frac{n}{2}\log(\sigma^2)
-\frac12
\log\left|
\bm\Sigma_{1,\beta,\bm t}
\right|
-\frac{1}{2\sigma^2}
\bm z^\top
\bm\Sigma_{1,\beta,\bm t}^{-1}
\bm z.
\end{equation}

For convenience, define
\[
Q_\beta(\bm z)
:=
\bm z^\top
\bm\Sigma_{1,\beta,\bm t}^{-1}
\bm z.
\]
Differentiating \eqref{eq:loglik_scaled} with respect to the scale
parameter \(\sigma^2\) gives
\begin{equation}\label{eq:dloglik_dsigma2}
\frac{\partial}{\partial(\sigma^2)}
\ell(\sigma^2,\beta\mid\bm z)
=
-\frac{n}{2\sigma^2}
+
\frac{Q_\beta(\bm z)}{2(\sigma^2)^2}.
\end{equation}
Therefore, for each fixed \(\beta\), the maximum likelihood estimator
of \(\sigma^2\) is
\begin{equation}\label{eq:sigmahat}
\widehat{\sigma}^{\,2}(\beta)
=
\frac{Q_\beta(\bm z)}{n}.
\end{equation}

Substituting \eqref{eq:sigmahat} into
\eqref{eq:loglik_scaled} yields the profile log-likelihood
\begin{equation}\label{eq:profile_beta}
\ell_p(\beta\mid\bm z)
=
-\frac{n}{2}
\left[
1+\log(2\pi)
+\log\!\left(
\frac{Q_\beta(\bm z)}{n}
\right)
\right]
-\frac12
\log\left|
\bm\Sigma_{1,\beta,\bm t}
\right|.
\end{equation}
It follows that
\begin{equation}\label{eq:betahat}
\widehat{\beta}
\in
\operatorname*{arg\,max}_{\beta\in\mathcal B}
\ell_p(\beta\mid\bm z),
\qquad
\widehat{\sigma}^{\,2}
=
\widehat{\sigma}^{\,2}(\widehat{\beta}),
\end{equation}
where \(\mathcal B\subseteq[0,\infty)\) denotes the parameter set used
in the numerical maximization. The value \(\beta=0\) corresponds to the
boundary model with constant weight.
\newpage
An analogous profiling argument applies to the process \(\mu_H\); see
\citet[Section~4]{JY}. Indeed, its covariance function in \eqref{QMH}
also depends linearly on the scale parameter~\(\sigma^2\).

\begin{rem}
The preceding likelihood formulation extends directly to a deterministic
initial condition \(\zeta_0\in\mathbb R\). In this case,
\[
\big(
\zeta_{t_1}-\zeta_0,\ldots,\zeta_{t_n}-\zeta_0
\big)^\top
\sim
N\!\left(
\bm 0,
\bm\Sigma_{\sigma^2,\beta,\bm t}
\right),
\]
because subtraction of a deterministic constant changes the mean but
does not alter the covariance matrix. Consequently, the likelihood is
evaluated using the shifted observation vector
\[
\big(
\zeta_{t_1}-\zeta_0,\ldots,\zeta_{t_n}-\zeta_0
\big)^\top.
\]
\end{rem}



\subsection{Inference and application to telemetry data}\label{bat}

In \cite{mara}, the authors recorded the three-dimensional migratory movements of five bats in Germany during 2016–2017. Individual flights varied both within and among bats, indicating that flight decisions are influenced by factors such as wind, landscape, navigational conditions, or other yet-unknown factors. The diagnostics reported in Appendix~C provide evidence against weak
stationarity of the raw telemetry series. For the first-differenced series,
the ADF and KPSS results are mixed: both procedures are consistent with
stationarity in four of the ten cases, whereas at least one of them indicates
non-stationarity in each of the remaining six cases. Autocorrelation remaining
after first differencing indicates serial dependence, but does not by itself
imply lack of intrinsic stationarity. The DFA analysis provides complementary evidence of large-scale dependence:
over the selected scaling ranges, the estimated exponents exceed one for all
ten telemetry series. Thus, stationary first differences are supported for a
non-negligible subset of the data, although not uniformly across all series.
Accordingly, the comparison below includes both stationary-increment models,
such as fractional Brownian motion, and models without stationary increments.

For the class \(f_{\sigma^2,\beta}(u)=\sigma^2e^{-\beta u}\), the process
\(\zeta\) is non-stationary and is not intrinsically stationary on compact
time intervals; see Appendix~B, Proposition~B.13. In addition,
Proposition~\ref{LRM}(ii) establishes its LRD property. These properties
motivate the consideration of \(\zeta\) for modeling the bat trajectories.

\begin{figure}[h]
\centering
\includegraphics[width=100mm]{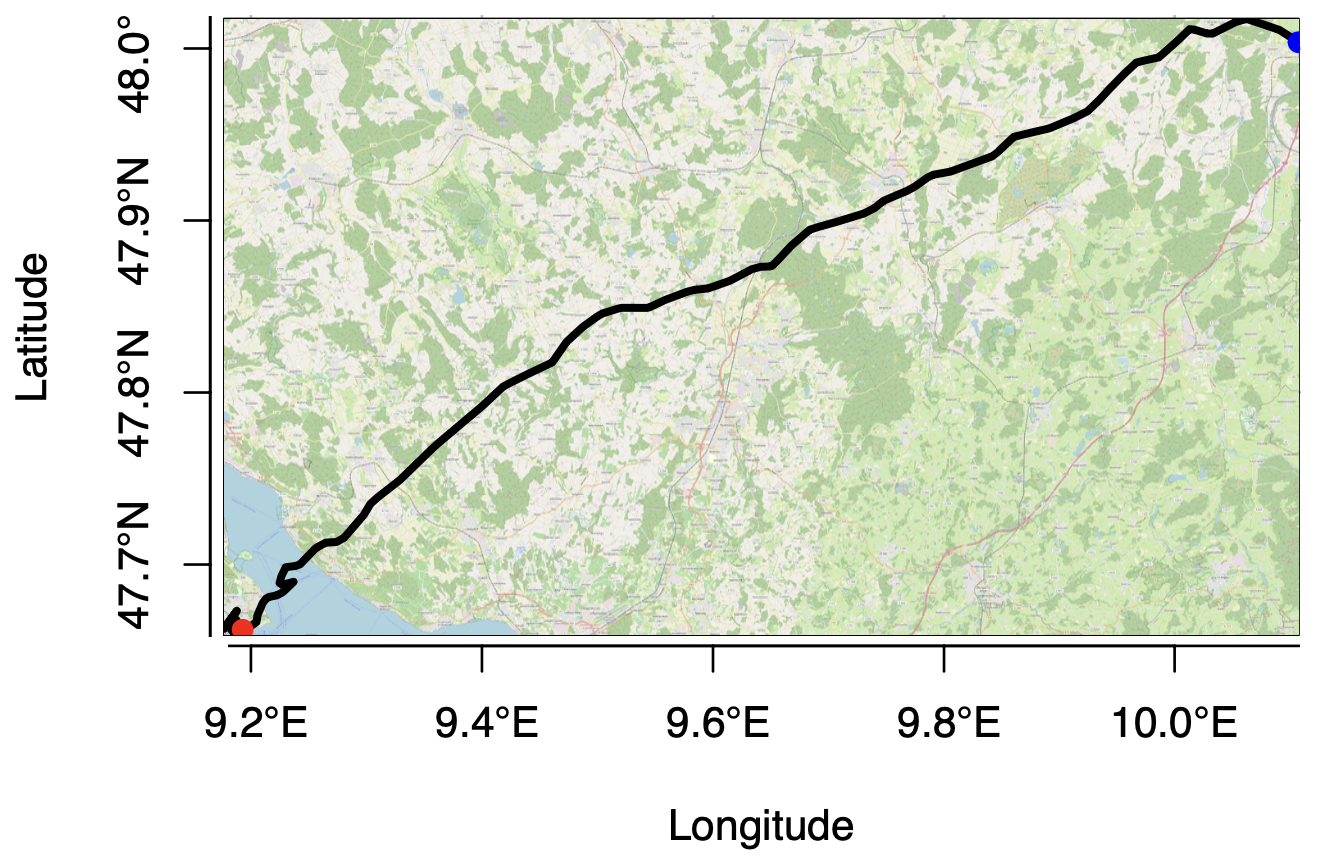}
\caption{Trajectory of bat ID $\#4$ showing longitude and latitude. The red point indicates the initial position, and the  blue point indicates the final position.}
 \label{mapa}
\end{figure}

Assuming independence between the coordinate processes, we fit
\((\zeta_t)_{t\geq0}\), with \(f_{\sigma^2,\beta}\) as above, separately
to the longitude and latitude telemetry series. We compare the fitted model
with the animal telemetry model \(\mu_H\) using the Akaike information
criterion (AIC). This dataset is openly available in the repository: \url{https://doi.org/10.5441/001/1.5d736bf0}. Among the five telemetry trajectories, altitude information is missing in two cases, limiting our analysis to latitude and longitude data. We focus on the telemetry of bat ID $\#4$ to describe the inference. To manage computational costs, we generate a periodically sampled trajectory by averaging the location data every 0.5 minutes. Figure \ref{mapa} illustrates the resulting trajectory.

\begin{figure}[H]
\centering
\includegraphics[width=0.98\textwidth]{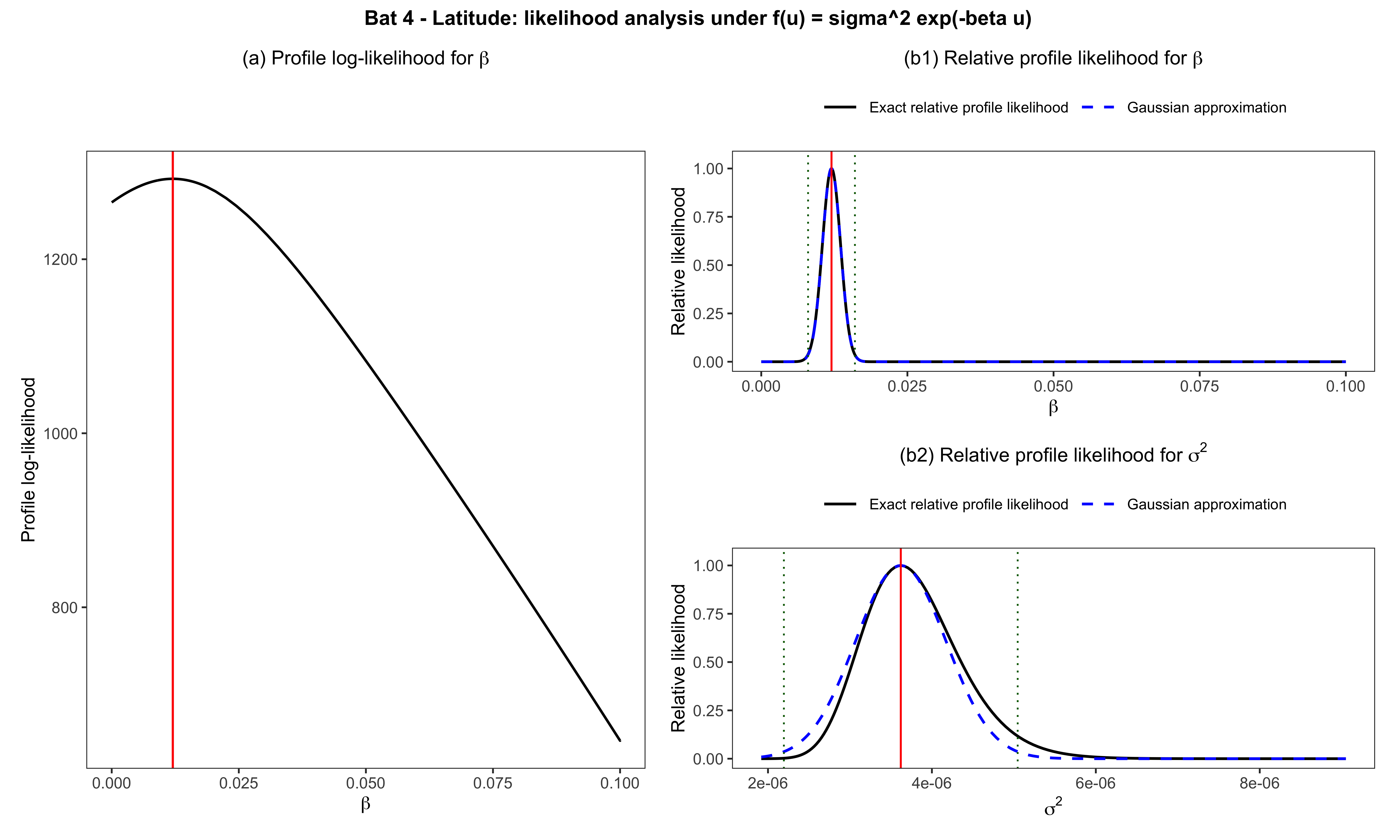}
\caption{Likelihood analysis for the model \((\zeta_{t,f_{\sigma^2,\beta}})_{t\geq0}\), with \(f_{\sigma^2,\beta}(u)=\sigma^2e^{-\beta u}\), fitted to the latitude telemetry data of Bat ID~\#4. Panel~(a) shows the profile log-likelihood for \(\beta\). Panels~(b1) and~(b2) show the relative profile likelihoods for \(\beta\) and \(\sigma^2\), respectively, together with their local Gaussian approximations. Each relative likelihood is normalized to attain its maximum value of one. The red vertical lines indicate the maximum likelihood estimates, and the green dotted lines indicate the approximate 99\% confidence limits obtained from the Gaussian approximations.}
\label{fig:profile_latitude}
\end{figure}

Inference was conducted by maximum likelihood using \texttt{R}, with the corresponding code available at \url{https://github.com/joseramirezgonzalez/f-wsfBm}. For the latitude telemetry data, under the model with \(f_{\sigma^2,\beta}(u)=\sigma^2e^{-\beta u}\), the maximum likelihood estimates are \((\widehat{\beta},\widehat{\sigma}^{\,2})=(0.0120167,3.620324\times10^{-6})\), with AIC value \(-2580.865\). A Gaussian approximation based on the local quadratic behavior of the profile log-likelihoods yields the confidence intervals shown in Figure~\ref{fig:profile_latitude}; see, for example, \cite{Pawitan2001}. At the \(99\%\) confidence level, the resulting intervals are \((0.008009,0.016024)\) for \(\beta\) and \((2.193960\times10^{-6},5.046688\times10^{-6})\) for \(\sigma^2\).

An analogous analysis was carried out for the longitude telemetry data. The profile log-likelihood for \(\beta\) attained its maximum at the boundary value \(\beta=0\). We therefore considered the boundary model associated with the constant weight \(f_{\sigma^2,0}(u)\equiv\sigma^2\). For this model, the maximum likelihood estimate is \(\widehat{\sigma}^{\,2}=5.699848\times10^{-6}\), with maximized log-likelihood \(1123.186\) and AIC value \(-2244.372\). The corresponding \(99\%\) confidence interval for \(\sigma^2\), obtained from the Gaussian approximation, is \((4.370499\times10^{-6},\,7.029196\times10^{-6})\).

\begin{figure}[H]
\centering
\includegraphics[width=0.98\textwidth]
{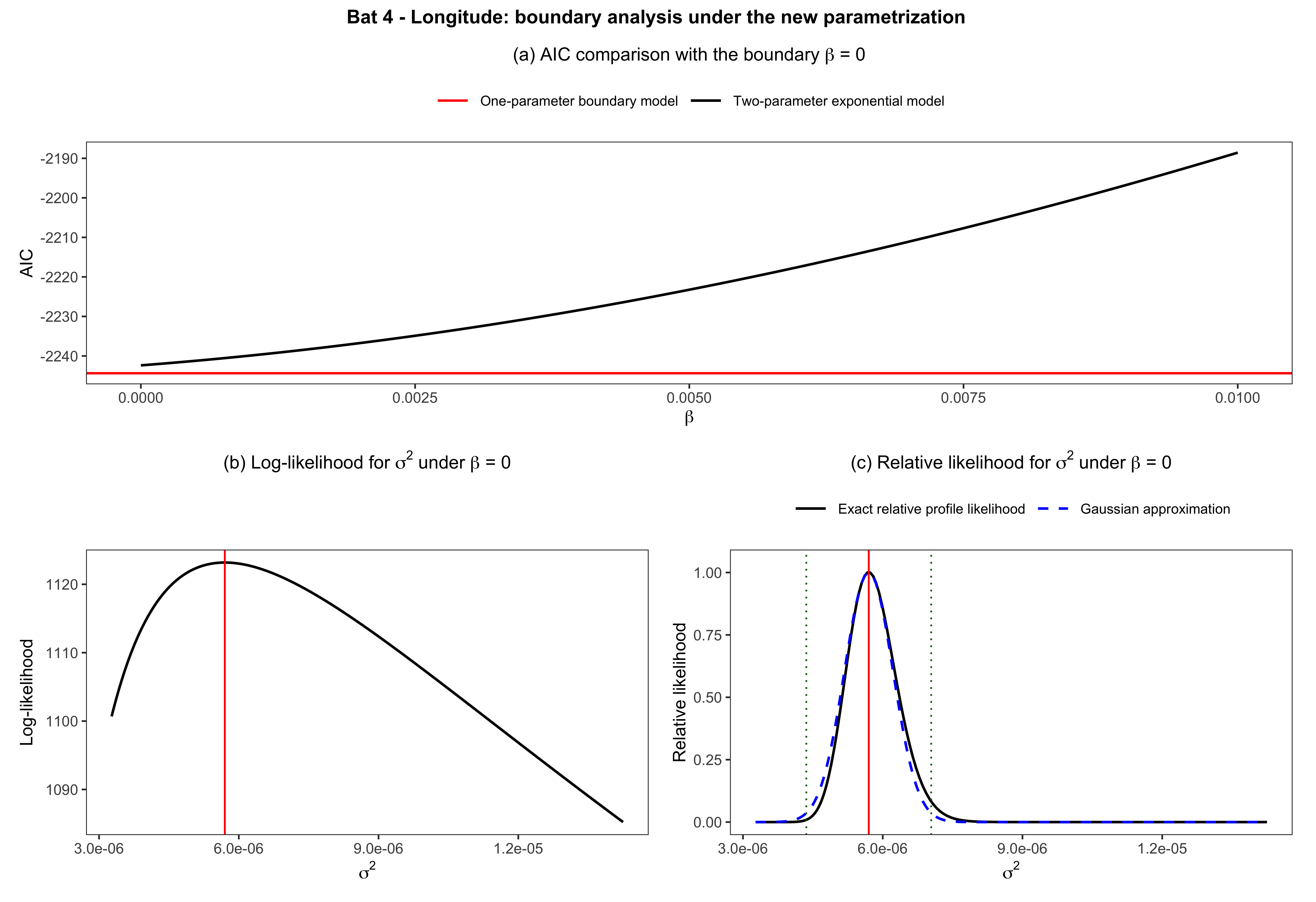}
\caption{Boundary analysis for the model
\(\left(\zeta_{t,f_{\sigma^2,\beta}}\right)_{t\geq0}\), with
\(f_{\sigma^2,\beta}(u)=\sigma^2e^{-\beta u}\), fitted to the longitude
telemetry data of Bat ID~\#4. Panel~(a) compares the AIC of the
two-parameter exponential model for \(\beta>0\) with that of the
one-parameter boundary model \(\beta=0\), for which
\(f_{\sigma^2,0}(u)\equiv\sigma^2\). Panel~(b) shows the log-likelihood
for \(\sigma^2\) under the boundary model. Panel~(c) shows the
corresponding relative likelihood, normalized to attain a maximum value
of one, together with its local Gaussian approximation. The red vertical
lines indicate the maximum likelihood estimates, and the green dotted
lines indicate the approximate 99\% confidence limits obtained from the
Gaussian approximation.}
\label{fig:profile_longitude_boundary}
\end{figure}

Figure~\ref{fig:profile_longitude_boundary} summarizes the results. Panel~(a) compares the AIC values of the exponential model with \(\beta>0\) and the boundary model with \(\beta=0\), with the latter providing the better fit. Panels~(b) and~(c) show, respectively, the log-likelihood for \(\sigma^2\) and the corresponding relative likelihood together with its Gaussian approximation.

Inference for the remaining four bat trajectories, in both latitude and
longitude, was carried out in the same manner. Table~\ref{T1} also reports
comparisons with several competing models having different stationarity,
dependence, and covariance-growth structures. Specifically, we compare the
proposed non-stationary and non-intrinsically stationary LRD model
\(\zeta_{t,f_{\sigma^2,\beta}}\), including the boundary specification
\(\zeta_{t,f_{\sigma^2,0}}\) for the longitude trajectory of Bat~\#4, with
the fractional Ornstein--Uhlenbeck position model \(\mu_H(t)\), for
\(H\in(0,1)\setminus\{\frac12\}\); the integrated Ornstein--Uhlenbeck
benchmark \(\mu_{1/2}(t)\); the fractional Brownian motion model
\(\sigma W_H(t)\); the stationary confluent hypergeometric model
\(A_{\sigma^2,\eta,\alpha,\beta}(t)\); and two stationary temporal
covariance models derived from Equation~(11) of \citet{Stein2005}.

The covariance of \(A_{\sigma^2,\eta,\alpha,\beta}(t)\) is
\[
\operatorname{Cov}\!\left(
A_{\sigma^2,\eta,\alpha,\beta}(t),
A_{\sigma^2,\eta,\alpha,\beta}(s)
\right)
=
\frac{\sigma^2\Gamma(\eta+\alpha)}{\Gamma(\eta)}
U\!\left(
\alpha,1-\eta,
\eta\left(\frac{|t-s|}{\beta}\right)^2
\right),
\qquad s,t\geq0,
\]
where \(U(a,b,z)\) denotes Tricomi's confluent hypergeometric function.
As shown in Appendix~D, Proposition~D.3, this is a valid stationary
covariance belonging to the confluent hypergeometric class of
\citet{MaBhadra2023}. Its local mean-square regularity is determined by
\(\eta\), its covariance has the polynomial tail
\(R_A(h)\sim c_A|h|^{-2\alpha}\), and it satisfies the LRD criterion adopted
in this article with normalization \(T^{2\alpha+2}\).

For the comparison with \citet{Stein2005}, we use the fixed-site temporal
marginal of Equation~(11). This is the appropriate temporal reduction
because each coordinate is modeled as a scalar response indexed by time;
the observed latitude and longitude values are not external spatial sites
of a space--time random field. Let \(\varepsilon\) denote the advection
parameter, let \(d\) be the spatial dimension appearing in Equation~(11),
and let \(\nu>0\) and \(\rho\geq0\) denote the parameters governing the
order of its temporal marginal. Proposition~D.1 of Appendix~D shows
that this marginal satisfies the adopted LRD criterion if and only if
\(\varepsilon=0\) and \(\rho>0\). Under precisely these conditions, its
covariance is
\[
R_{d,\nu,\rho}(h)
=
\sigma^2
\frac{\Gamma(\nu+d/2)}{\Gamma(\nu)}
\frac{\Gamma(\nu+\rho|h|)}
     {\Gamma(\nu+\rho|h|+d/2)},
\qquad h\in\mathbb{R}.
\]
The other fixed-site temporal branches of Equation~(11) are either
constant or exponentially decaying and do not satisfy the LRD criterion
used here.

For the numerical comparison, we considered the LRD-compatible fixed-site
temporal covariance \(R_{d,\nu,\rho}\) for \(d=1,2,3\), under the
parametrization
\[
\rho=\nu\lambda,\qquad \lambda>0.
\]
In the likelihood maximization for the ten coordinate series, the cases
\(d=1\) and \(d=3\) did not yield interior maximizers in \(\nu\). For
\(d=1\), the likelihood profiles approached their suprema as
\(\nu\downarrow0\). By Corollary~D.2 of Appendix~D, the corresponding
limiting covariance is
\[
\frac{\sigma^2}{1+\lambda|h|}.
\]
For \(d=2\), the gamma recurrence formula gives
\[
R_{2,\nu,\nu\lambda}(h)
=
\frac{\sigma^2}{1+\lambda|h|},
\]
so the covariance does not depend on \(\nu\). Hence, the boundary case
obtained from \(d=1\) coincides with the exact \(d=2\) reduction and does
not define an additional identifiable family.

For \(d=3\), the likelihood profiles approached their suprema as
\(\nu\to\infty\). Again by Corollary~D.2 of Appendix~D, the
corresponding limiting covariance is
\[
\frac{\sigma^2}{(1+\lambda|h|)^{3/2}}.
\]
Consequently, the comparison reported in Table~\ref{T1} includes the two
distinct Stein-type covariance models
\[
\operatorname{Cov}\!\left(
S_{\sigma^2,\lambda}^{(d)}(t),
S_{\sigma^2,\lambda}^{(d)}(s)
\right)
=
\frac{\sigma^2}
     {(1+\lambda|t-s|)^{d/2}},
\qquad d\in\{2,3\}.
\]
Here, \(S_{\sigma^2,\lambda}^{(2)}\) is the exact \(d=2\) fixed-site
temporal marginal of Equation~(11) of \citet{Stein2005}, whereas
\(S_{\sigma^2,\lambda}^{(3)}\) is the limiting \(d=3\) model obtained as
\(\nu\to\infty\). Both are stationary Gaussian models satisfying the LRD
criterion adopted in this article; see Corollary~D.2 of Appendix~D.

The proposed \(\zeta\)-type model is neither stationary nor intrinsically
stationary and exhibits LRD. For
\(H\in(0,1)\setminus\{\frac12\}\), the process \(\mu_H(t)\) is likewise
neither stationary nor intrinsically stationary and exhibits LRD. Its
Brownian counterpart \(\mu_{1/2}(t)\), corresponding to the integrated
Ornstein--Uhlenbeck process, is neither stationary nor intrinsically
stationary, but it does not satisfy the LRD criterion adopted in this
article because its separated-increment covariance decays exponentially.

The process \(\sigma W_H(t)\) is non-stationary and intrinsically
stationary, belonging to the short-memory regime when \(H<1/2\) and to
the long-memory regime when \(H>1/2\); when \(H=1/2\), it is standard
Brownian motion and has independent increments; see
\cite{Gorostiza2015}. By contrast,
\(A_{\sigma^2,\eta,\alpha,\beta}(t)\),
\(S_{\sigma^2,\lambda}^{(2)}(t)\), and
\(S_{\sigma^2,\lambda}^{(3)}(t)\) are stationary Gaussian processes with
LRD.

For some telemetry datasets, the profile likelihood of \(\beta\) under
the model \(\mu_H(t)\) was flat and increasing over the numerically
admissible range, so that a finite maximum likelihood estimate of
\(\beta\) could not be obtained. In those cases, the numerical search
interval was \([0.001,400]\) and, writing \(\ell_{\mathrm p}\) for the
profile log-likelihood, we defined \(\beta^*\) as the smallest value
\(\beta\in[0.001,400]\) satisfying
\[
\ell_{\mathrm p}(400)-\ell_{\mathrm p}(\beta^*)\leq 10^{-2},
\]
up to numerical tolerance. Thus, \(\beta^*\) is a representative of the
likelihood plateau rather than a finite maximum likelihood estimate of
\(\beta\). These cases are denoted by \(\beta^*\) in Table~\ref{T1}.

The AIC values in Table~\ref{T1} show that the proposed
\(\zeta\)-type model attains the smallest AIC in three of the ten
datasets, fractional Brownian motion \(\sigma W_H(t)\) does so in five,
\(\mu_H(t)\) does so in one, and the stationary confluent hypergeometric
model \(A_{\sigma^2,\eta,\alpha,\beta}(t)\) does so in one. Neither
\(\mu_{1/2}(t)\), \(S_{\sigma^2,\lambda}^{(2)}(t)\), nor
\(S_{\sigma^2,\lambda}^{(3)}(t)\) attains the smallest AIC in any
of the ten datasets.

These results are consistent with the empirical diagnostics reported
in Appendix~C, which provide evidence against weak stationarity
and show that stationarity of the first-differenced series is not
uniformly supported across the ten datasets. Both ADF and KPSS are
simultaneously consistent with stationary first differences in four cases,
and KPSS does not reject stationarity in two additional cases. This makes a
stationary-increment description plausible for a substantial subset of the
data and is compatible with the strong performance of fractional Brownian
motion, without implying a one-to-one correspondence between the diagnostic
tests and the full Gaussian likelihood comparison. All estimates obtained
under \(\sigma W_H(t)\) and \(\mu_H(t)\) satisfy
\(\widehat{H}>1/2\), and hence correspond to the long-memory regimes of
those models. Moreover, seven of the \(\sigma W_H(t)\) fits yield
\(\widehat{H}>0.75\). In four of those seven cases, the proposed
\(\zeta\)-type model has a smaller AIC than \(\sigma W_H(t)\).
Thus, even among trajectories for which the fractional Brownian motion
fits indicate pronounced long memory, the \(\zeta\)-type model remains
competitive. This empirical finding complements the theoretical contrast
established above between the logarithmic covariance growth of the proposed
model and the polynomial covariance growth of \(\mu_H(t)\).
\begin{table*}[p]
\centering
\rotatebox{90}{%
\scalebox{-1}[-1]{%
\begin{minipage}{\textheight}
\centering
\tiny
\setlength{\tabcolsep}{2pt}
\renewcommand{\arraystretch}{1.15}

\resizebox{0.97\textheight}{!}{%
\begin{tabular}{|c|c|
>{\centering\arraybackslash}p{2.10cm}|
>{\centering\arraybackslash}p{2.55cm}|
>{\centering\arraybackslash}p{2.10cm}|
>{\centering\arraybackslash}p{2.00cm}|
>{\centering\arraybackslash}p{3.35cm}|
>{\centering\arraybackslash}p{2.15cm}|
>{\centering\arraybackslash}p{2.25cm}|}
\hline
\multicolumn{2}{|c|}{Dataset}
& $\zeta_{t,f_{\sigma^2,\beta}}$
& $\mu_H(t)$
& $\mu_{\frac12}(t)$
& $\sigma W_H(t)$
& $A_{\sigma^2,\eta,\alpha,\beta}(t)$
& $S_{\sigma^2,\lambda}^{(2)}(t)$
& $S^{(3)}_{\sigma^2,\lambda}(t)$
\\ \hline

\multirow{3}{*}{\shortstack[c]{Latitude\\Bat \#1}}
& AIC
& $-263.7006$
& $-282.6614$
& $-253.4788$
& $\mathbf{-284.7042}$
& $-276.8744$
& $-246.9013$
& $-247.0100$
\\
& Parameters
& $(\beta,\sigma^2)$
& $(\beta^{*},\sigma,H)$
& $(\beta,\sigma)$
& $(H,\sigma)$
& $(\sigma^2,\eta,\alpha,\beta)$
& $(\sigma^2,\lambda)$
& $(\sigma^2,\lambda)$
\\
& Estimates
& $(0.006,8.496\times10^{-6})$
& $(219.167,0.859,0.779)$
& $(1.149,0.010)$
& $(0.779,0.004)$
& $(0.721,0.774,9.226,10940.590)$
& $(0.318,7.85\times10^{-5})$
& $(0.285,5.82\times10^{-5})$
\\ \hline

\multirow{3}{*}{\shortstack[c]{Longitude\\Bat \#1}}
& AIC
& $-231.4647$
& $-245.3374$
& $-230.0733$
& $\mathbf{-247.3626}$
& $-240.4568$
& $-220.4685$
& $-220.5866$
\\
& Parameters
& $(\beta,\sigma^2)$
& $(\beta^{*},\sigma,H)$
& $(\beta,\sigma)$
& $(H,\sigma)$
& $(\sigma^2,\eta,\alpha,\beta)$
& $(\sigma^2,\lambda)$
& $(\sigma^2,\lambda)$
\\
& Estimates
& $(0.004,1.521\times10^{-5})$
& $(118.714,0.722,0.703)$
& $(1.682,0.017)$
& $(0.703,0.006)$
& $(0.546,0.699,9.301,9128.046)$
& $(0.403,1.11\times10^{-4})$
& $(0.359,8.30\times10^{-5})$
\\ \hline

\multirow{3}{*}{\shortstack[c]{Latitude\\Bat \#2}}
& AIC
& $-277.5199$
& $-284.9209$
& $-272.5946$
& $\mathbf{-286.9427}$
& $-281.2260$
& $-261.6505$
& $-261.7662$
\\
& Parameters
& $(\beta,\sigma^2)$
& $(\beta^{*},\sigma,H)$
& $(\beta,\sigma)$
& $(H,\sigma)$
& $(\sigma^2,\eta,\alpha,\beta)$
& $(\sigma^2,\lambda)$
& $(\sigma^2,\lambda)$
\\
& Estimates
& $(0.015,1.557\times10^{-5})$
& $(99.147,0.494,0.694)$
& $(2.070,0.016)$
& $(0.694,0.005)$
& $(0.167,0.688,9.312,11433.090)$
& $(0.114,2.07\times10^{-4})$
& $(0.102,1.54\times10^{-4})$
\\ \hline

\multirow{3}{*}{\shortstack[c]{Longitude\\Bat \#2}}
& AIC
& $-252.4751$
& $-266.6972$
& $-257.4578$
& $\mathbf{-268.7201}$
& $-263.1141$
& $-247.8731$
& $-247.9854$
\\
& Parameters
& $(\beta,\sigma^2)$
& $(\beta^{*},\sigma,H)$
& $(\beta,\sigma)$
& $(H,\sigma)$
& $(\sigma^2,\eta,\alpha,\beta)$
& $(\sigma^2,\lambda)$
& $(\sigma^2,\lambda)$
\\
& Estimates
& $(0.0098,2.049\times10^{-5})$
& $(101.543,0.625,0.682)$
& $(2.328,0.021)$
& $(0.682,0.006)$
& $(0.159,0.675,9.325,9014.932)$
& $(0.134,2.41\times10^{-4})$
& $(0.120,1.79\times10^{-4})$
\\ \hline

\multirow{3}{*}{\shortstack[c]{Latitude\\Bat \#3}}
& AIC
& $-583.5546$
& $\mathbf{-586.5635}$
& $-585.7510$
& $-580.2283$
& $-579.6776$
& $-536.6979$
& $-536.8127$
\\
& Parameters
& $(\beta,\sigma^2)$
& $(\beta,\sigma,H)$
& $(\beta,\sigma)$
& $(H,\sigma)$
& $(\sigma^2,\eta,\alpha,\beta)$
& $(\sigma^2,\lambda)$
& $(\sigma^2,\lambda)$
\\
& Estimates
& $(0.014,5.021\times10^{-6})$
& $(1.755,0.006,0.682)$
& $(0.853,0.004)$
& $(0.860,0.003)$
& $(0.005,0.954,1.305,136.496)$
& $(0.009,4.18\times10^{-4})$
& $(0.008,3.12\times10^{-4})$
\\ \hline

\multirow{3}{*}{\shortstack[c]{Longitude\\Bat \#3}}
& AIC
& $\mathbf{-564.5836}$
& $-559.9343$
& $-533.2260$
& $-561.9265$
& $-556.2895$
& $-467.6996$
& $-467.8157$
\\
& Parameters
& $(\beta,\sigma^2)$
& $(\beta,\sigma,H)$
& $(\beta,\sigma)$
& $(H,\sigma)$
& $(\sigma^2,\eta,\alpha,\beta)$
& $(\sigma^2,\lambda)$
& $(\sigma^2,\lambda)$
\\
& Estimates
& $(0.036,2.412\times10^{-5})$
& $(193.950,0.670,0.871)$
& $(0.632,0.006)$
& $(0.872,0.003)$
& $(0.093,0.893,9.107,1974.061)$
& $(0.069,1.73\times10^{-4})$
& $(0.062,1.29\times10^{-4})$
\\ \hline

\multirow{3}{*}{\shortstack[c]{Latitude\\Bat \#4}}
& AIC
& $\mathbf{-2580.8650}$
& $-2534.8409$
& $-2473.3950$
& $-2536.6267$
& $-2531.1563$
& $-2298.2396$
& $-2298.3567$
\\
& Parameters
& $(\beta,\sigma^2)$
& $(\beta,\sigma,H)$
& $(\beta,\sigma)$
& $(H,\sigma)$
& $(\sigma^2,\eta,\alpha,\beta)$
& $(\sigma^2,\lambda)$
& $(\sigma^2,\lambda)$
\\
& Estimates
& $(0.012,3.620\times10^{-6})$
& $(19.766,0.049,0.844)$
& $(1.099,0.005)$
& $(0.857,0.002)$
& $(0.099,0.868,9.132,3336.100)$
& $(0.087,3.61\times10^{-5})$
& $(0.077,2.70\times10^{-5})$
\\ \hline

\multirow{3}{*}{\shortstack[c]{Longitude\\Bat \#4}}
& AIC
& $-2244.3720$
& $-2259.0606$
& $-2146.9770$
& $\mathbf{-2261.0896}$
& $-2250.8079$
& $-1913.0073$
& $-1913.1250$
\\
& Parameters
& $(\sigma^2)$
& $(\beta^{*},\sigma,H)$
& $(\beta,\sigma)$
& $(H,\sigma)$
& $(\sigma^2,\eta,\alpha,\beta)$
& $(\sigma^2,\lambda)$
& $(\sigma^2,\lambda)$
\\
& Estimates
& $(5.699\times10^{-6})$
& $(224.903,0.975,0.849)$
& $(0.762,0.009)$
& $(0.849,0.004)$
& $(3.867,0.085,9.148,16400.040)$
& $(0.624,2.50\times10^{-5})$
& $(0.555,1.87\times10^{-5})$
\\ \hline

\multirow{3}{*}{\shortstack[c]{Latitude\\Bat \#5}}
& AIC
& $-1648.9230$
& $-1649.6038$
& $-1590.8140$
& $-1651.6334$
& $\mathbf{-1652.2620}$
& $-1488.9655$
& $-1489.1822$
\\
& Parameters
& $(\beta,\sigma^2)$
& $(\beta^{*},\sigma,H)$
& $(\beta,\sigma)$
& $(H,\sigma)$
& $(\sigma^2,\eta,\alpha,\beta)$
& $(\sigma^2,\lambda)$
& $(\sigma^2,\lambda)$
\\
& Estimates
& $(0.004,3.134\times10^{-6})$
& $(231.861,0.797,0.888)$
& $(1.284,0.007)$
& $(0.888,0.003)$
& $(0.005,0.901,9.099,406.667)$
& $(0.004,1.21\times10^{-3})$
& $(0.003,9.86\times10^{-4})$
\\ \hline

\multirow{3}{*}{\shortstack[c]{Longitude\\Bat \#5}}
& AIC
& $\mathbf{-1585.8437}$
& $-1581.2897$
& $-1535.6170$
& $-1583.3100$
& $-1584.8542$
& $-1392.8529$
& $-1393.0755$
\\
& Parameters
& $(\beta,\sigma^2)$
& $(\beta^{*},\sigma,H)$
& $(\beta,\sigma)$
& $(H,\sigma)$
& $(\sigma^2,\eta,\alpha,\beta)$
& $(\sigma^2,\lambda)$
& $(\sigma^2,\lambda)$
\\
& Estimates
& $(0.001,4.166\times10^{-6})$
& $(168.604,0.999,0.945)$
& $(0.856,0.007)$
& $(0.945,0.006)$
& $(0.008,0.097,9.029,310.546)$
& $(0.006,1.43\times10^{-3})$
& $(0.005,1.17\times10^{-3})$
\\ \hline

\end{tabular}%
}

\captionof{table}{Comparison of the fitted models for the latitude and
longitude telemetry data from five bat trajectories. For each dataset,
the table reports the AIC value, parameter vector, and corresponding
estimates; the smallest AIC is shown in bold. For Bat \#4 Longitude,
the \(\zeta\)-type model corresponds to the boundary case \(\beta=0\),
for which \(f_{\sigma^2,0}(u)=\sigma^2\) and only \(\sigma^2\) is
estimated. Here, \(\beta^{*}\) denotes the smallest value
\(\beta\in[0.001,400]\) satisfying
\(\ell_{\mathrm p}(400)-\ell_{\mathrm p}(\beta^{*})\leq 10^{-2}\),
up to numerical tolerance, according to the rule described in the text.
The model \(S^{(3)}_{\sigma^2,\lambda}\) is the identifiable
two-parameter boundary model obtained from the \(d=3\) temporal
covariance family as \(\nu\to\infty\).}
\label{T1}

\end{minipage}%
}}
\end{table*}

\subsection{On the R Shiny interactive app}

To facilitate simulation, estimation, and prediction for the process \(\zeta\) studied in this article, we developed an interactive R Shiny application, publicly available at \url{https://jhramgon.shinyapps.io/aplication_r_log_w/}. The application is specifically tailored to the covariance structure considered here and is organized into three modules. In the simulation modules, trajectories are generated over the grid \((0,\frac{T}{n},2\frac{T}{n},\dots,T)\), while a common header displays the covariance function associated with the process.

\begin{figure}[H]
    \centering
    \includegraphics[width=0.8\textwidth]{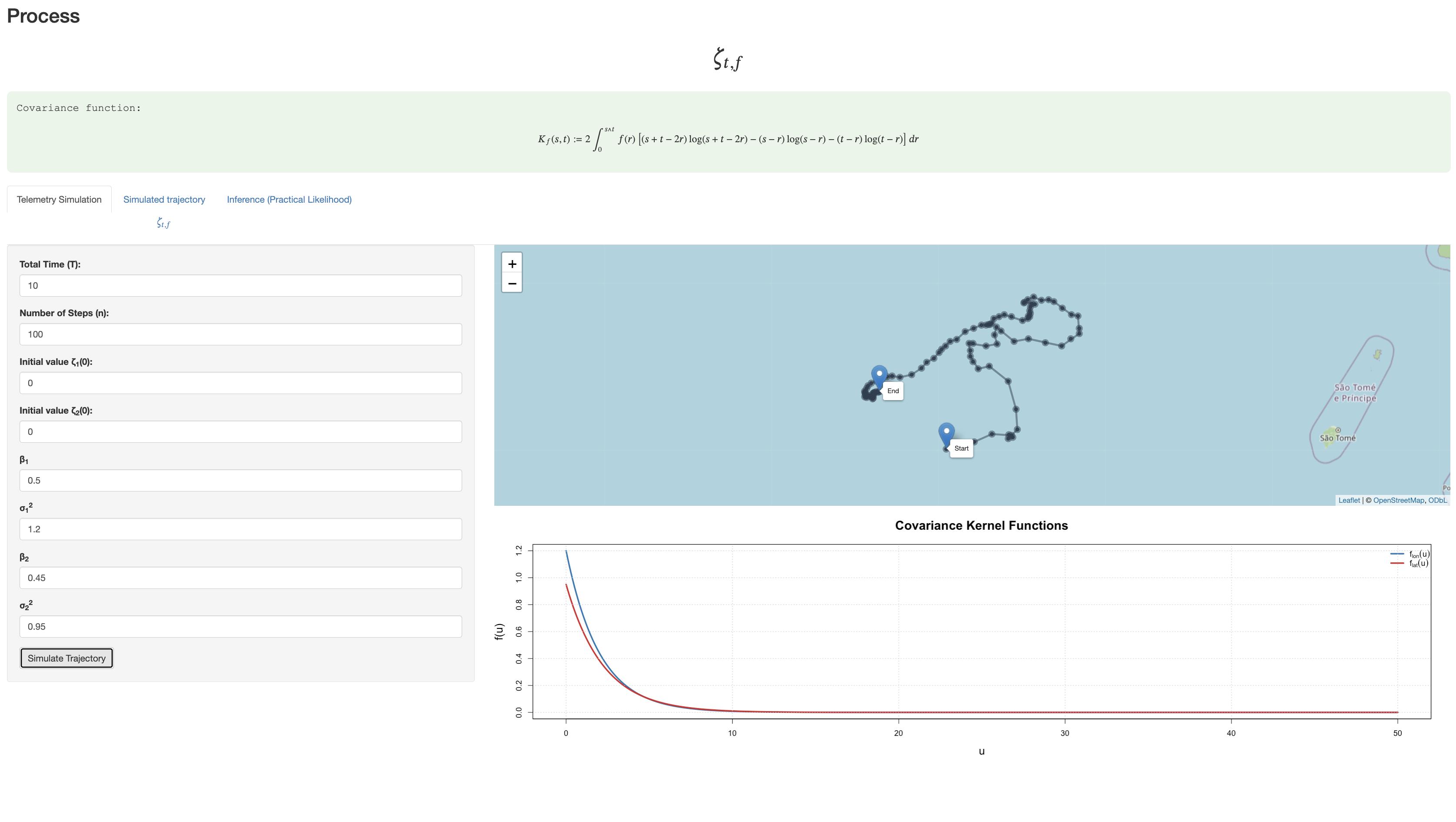}
    \caption{Output of the \textit{Telemetry Simulation} module. A realization of a movement path generated from two independent Gaussian coordinates associated with the process \(\zeta\) is displayed on a geographic map. The lower panel shows the kernel functions used for the longitude and latitude coordinates.}
    \label{fig:telemetry_sim}
\end{figure}

The first module, \textit{Telemetry Simulation}, is designed to generate simulated movement trajectories on a geographic map. In this module, the longitude and latitude coordinates are modeled as two independent Gaussian processes with covariance kernels of the form given in~(\ref{critical-wegithed-covariance}), using the functions \(f_{\sigma_i^2,\beta_i}(u)=\sigma_i^2e^{-\beta_i u}\), \(i=1,2\). The user specifies the time horizon, the number of discretization steps, the initial coordinates, and the parameters \((\beta_1,\sigma_1^2,\beta_2,\sigma_2^2)\). Once the \textit{Simulate Trajectory} button is activated, the application generates a realization of the corresponding bivariate path, displays it interactively on a geographic map, and simultaneously plots the kernel functions associated with the longitude and latitude coordinates. This module thus provides a direct visual connection between the chosen covariance specification and the resulting trajectory; see Figure~\ref{fig:telemetry_sim}.

\begin{figure}[H]
    \centering
    \includegraphics[width=0.92\textwidth]{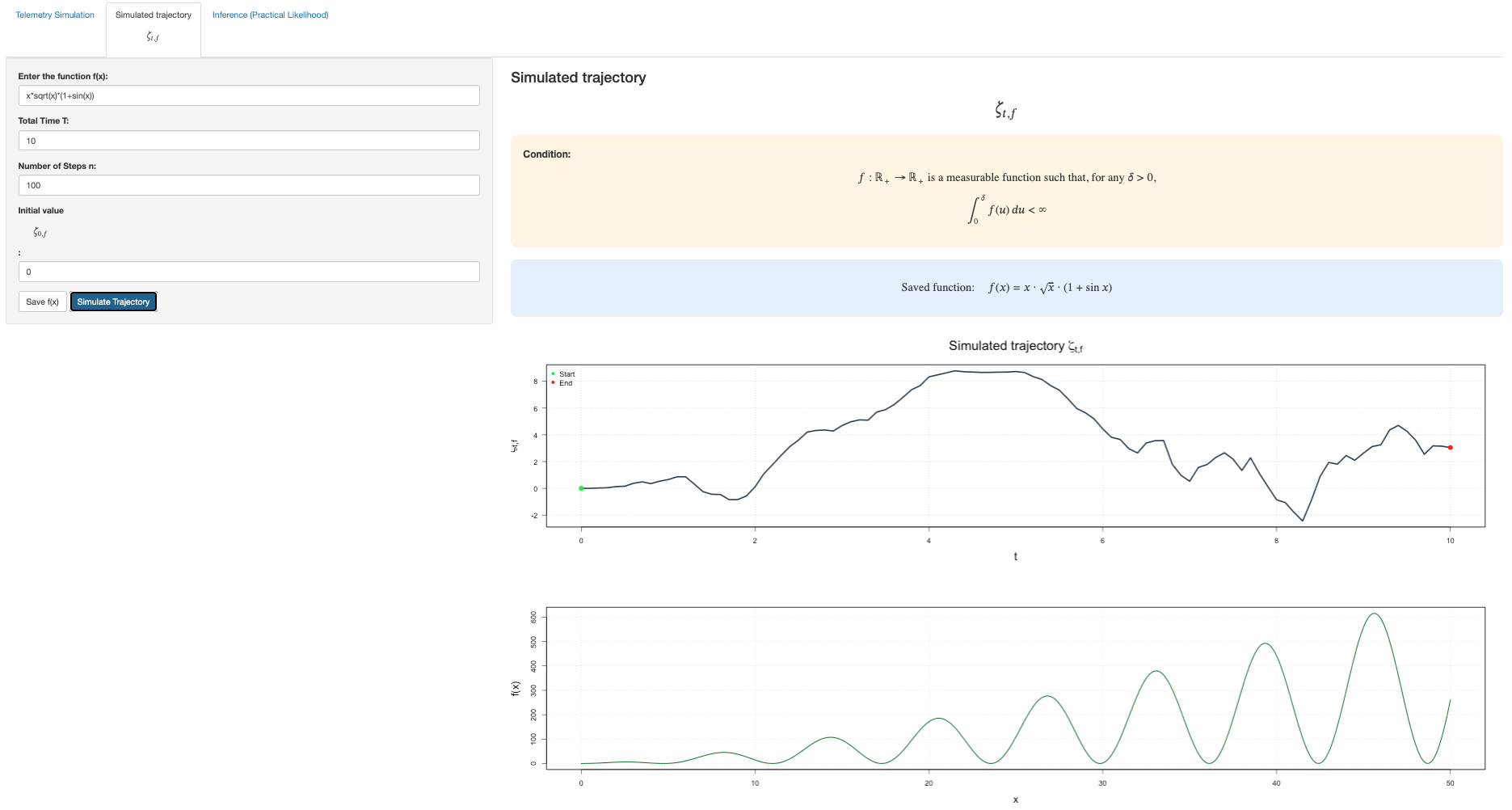}
    \caption{Output of the \textit{Simulated trajectory \(\zeta_{t,f}\)} module. The panel displays the integrability condition imposed on \(f\), the user-defined function entered in the application, a realization of the process \(\zeta_{t,f}\), and the graph of the selected function \(f\) used to construct the covariance matrix in~(\ref{critical-wegithed-covariance}).}
    \label{fig:sim_trajectory_fx}
\end{figure}

The second module, \textit{Simulated trajectory \(\zeta_{t,f}\)}, allows the user to investigate the model under a user-specified kernel function \(f\). The function must be entered in standard R syntax and is required to satisfy the integrability condition in~(\ref{Integrability-condition}). Once the expression has been saved, the application displays both the admissibility condition and the mathematical representation of the selected function. It then constructs the covariance matrix associated with~(\ref{critical-wegithed-covariance}) and generates a realization of the Gaussian process \(\zeta_{t,f}\). The output consists of the simulated trajectory together with the graph of the chosen function \(f\), thereby illustrating how the specification of the kernel influences the resulting sample path. An example of this module is shown in Figure~\ref{fig:sim_trajectory_fx}.

The third module, \textit{Inference (Practical Likelihood)}, implements the estimation and prediction procedure developed in this article for the \(\zeta\)-process. It is restricted to the parametric specification adopted in the inferential part of the paper, namely the kernel family \(f_{\sigma^2,\beta}(u)=\sigma^2e^{-\beta u}\). The user uploads a telemetry data set in Excel format, organized in consecutive three-column blocks corresponding to time, longitude, and latitude for each animal. Once the animal of interest has been selected, the corresponding block is extracted and the recorded observation times are converted into integer indices through \(\mathrm{round}(t/\Delta)\), where \(\Delta\) is specified by the user. The interface also allows separate search intervals for \(\beta\) in the longitude and latitude coordinates, the specification of the grid size for profile-likelihood evaluation, and the choice of the prediction settings.

After the \textit{Run inference} button is activated, the module estimates \(\widehat{\beta}\) and \(\widehat{\sigma^2}\) separately for longitude and latitude by maximizing the practical log-likelihood under the fitted model. For each coordinate, the observed trajectory is expressed relative to its initial value, the corresponding covariance matrix is evaluated at the observation times, and numerical optimization over the prescribed \(\beta\)-interval yields the parameter estimates. The resulting values of \(\widehat{\beta}\), \(\widehat{\sigma^2}\), and the associated AIC criterion are then reported. When requested, the profile log-likelihood is also evaluated on a user-defined grid and the corresponding curve is displayed for each coordinate.

\begin{figure}[H]
    \centering
    \includegraphics[width=0.92\textwidth]{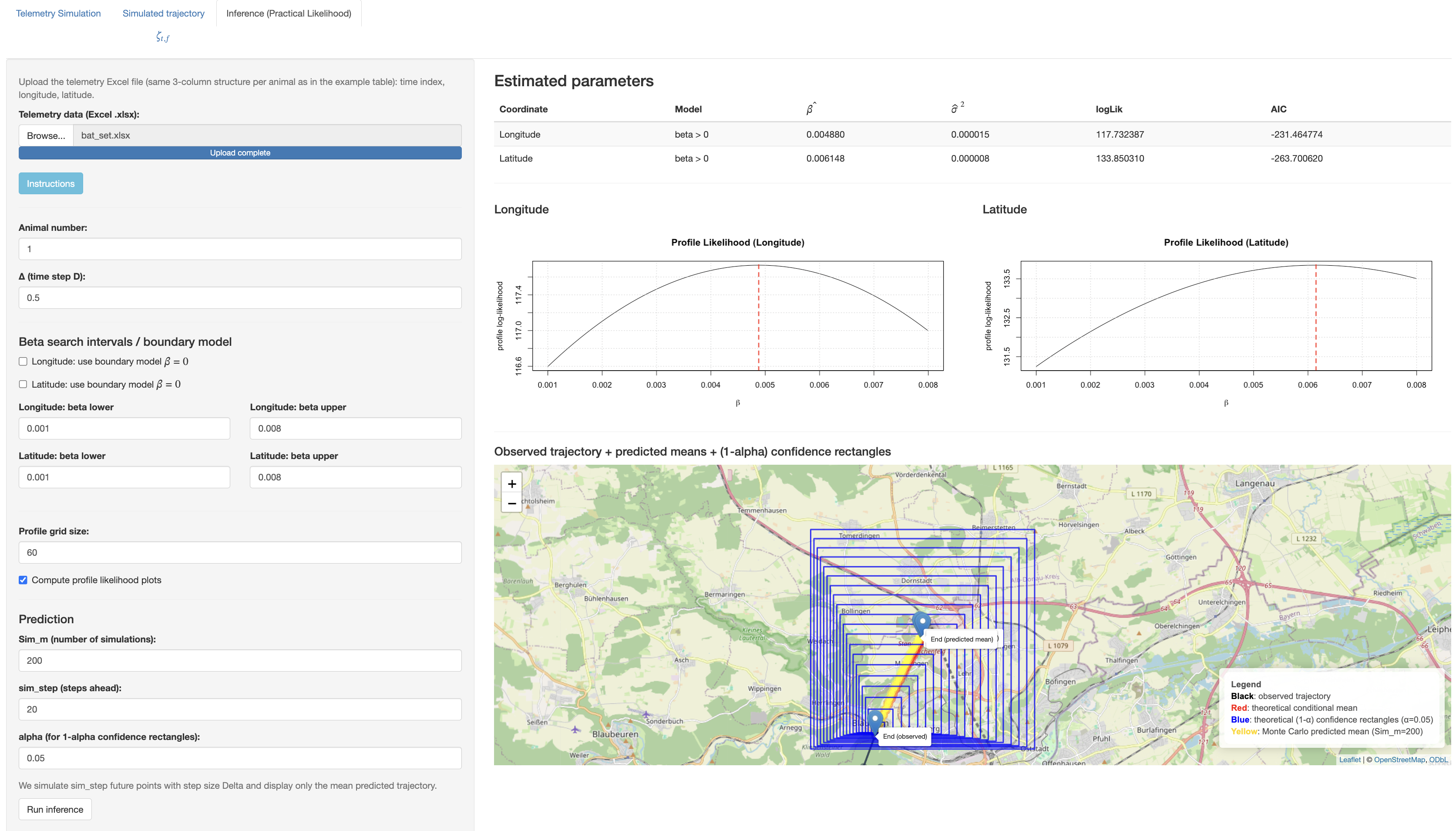}
    \caption{Output of the \textit{Inference (Practical Likelihood)} module for the \(\zeta\)-process considered in this article. The left panel contains the data input and tuning parameters, including the animal identifier, the time-step \(\Delta\), the search intervals for \(\beta\), the profile-likelihood settings, and the prediction controls. The right panel reports the estimates \(\widehat{\beta}\) and \(\widehat{\sigma^2}\), the corresponding AIC values, the profile log-likelihood curves for longitude and latitude, and an interactive map displaying the observed trajectory together with the theoretical conditional mean path, the coordinatewise marginal prediction rectangles, and the Monte Carlo mean trajectory.}
    \label{fig:inference_app}
\end{figure}

The module also includes a prediction component based on the fitted model. For each coordinate separately, it computes the conditional mean trajectory over the selected number of future steps. This output is used to display the theoretical conditional mean path and, at each future step, the corresponding coordinatewise marginal \((1-\alpha)\) prediction intervals for longitude and latitude. The blue rectangles shown on the map are obtained by combining these marginal intervals coordinatewise and should therefore be interpreted as Cartesian products of separate marginal intervals rather than as joint bivariate \((1-\alpha)\) prediction regions. In addition, when Monte Carlo prediction is requested, future continuations for longitude and latitude are simulated under the fitted model, and their average continuation is displayed on the map. These outputs are superimposed on the observed trajectory, thereby allowing estimation and short-term prediction to be examined jointly; see Figure~\ref{fig:inference_app}. In particular, for the \(\zeta\)-process considered in this article, the inferential analysis reported in Table~\ref{T1} can be reproduced in a direct and transparent manner within the application.

For practical use, the inference module also includes an \textit{Instructions} button that opens a guidance window describing the required structure of the telemetry file, the role of each estimation and prediction parameter, and the interpretation of the numerical and graphical outputs. This additional feature makes the module operational and self-contained from the user’s perspective.

\section{Conclusions}\label{Discu}

We have introduced a family of Gaussian processes and focused our
investigation on the trajectory and asymptotic properties of two specific
classes, denoted by \(\mathcal{C}_1\) and \(\mathcal{C}_2\). These classes
are characterized by long-range dependence and covariance functions with
logarithmic asymptotic growth. In Section~\ref{bat}, we presented an
application to bat telemetry data from Germany. The proposed
\(\mathcal{C}_2\)-type model attained the smallest AIC among the candidate
models in three of the ten coordinate series. Fractional Brownian motion
attained the smallest AIC in five series, the fractional
Ornstein--Uhlenbeck position model \(\mu_H\) did so in one, and the
stationary confluent hypergeometric model did so in one. Thus, the
preferred models span stationary specifications, non-stationary models
with stationary increments, and non-stationary models without stationary
increments. As shown in Proposition~\ref{LRM}, the covariance of the
proposed model with remote future positions grows at logarithmic order.

In contrast, the fractional Ornstein--Uhlenbeck position model
\(\mu_H\) considered in \cite{JY}, with \(H\in(1/2,1)\), has polynomial
asymptotic covariance growth of order \(T^{2H-1}\), which is regularly
varying. Our results suggest that the proposed model provides a useful
alternative in settings where polynomial growth may overstate the
influence of remote past positions. In such cases, covariance growth
governed by a slowly varying function, such as the logarithm, may offer
a more suitable description of long-range dependence.

In Section~\ref{Application}, we adopted the common assumption that
telemetry data are uncorrelated across coordinate axes. Although this
assumption simplifies the analysis, it may be inadequate when dependence
across coordinates is non-negligible. A natural extension of the present
work is therefore to develop a multidimensional version of the proposed
model capable of capturing such cross-coordinate dependence.
Additionally, in Appendix A, we present simulations of finite-dimensional distributions of the process \(\zeta\) and examine the empirical performance of the MLEs for the parametric family of weight functions
\[
\mathcal{C}:=\left\{ f_{\sigma^2,\beta}(u)=\sigma^2e^{-\beta u} : \sigma^2>0,\beta\ge0 \right\}.
\]
Currently, no asymptotic results are available for parameter estimation within the families \(\mathcal{C}_1\) and \(\mathcal{C}_2\). Developing asymptotic estimators tailored to these classes constitutes a promising direction for future research.

Finally, we believe that the proposed model may also be applicable to model (and predict) the trajectories of other species with robust memory capabilities, such as  chimpanzees and  elephants, among others (\cite{Janmaat2013,Polansky2015}).

\noindent{\bf Acknowledgments.}
This study was supported by the King Abdullah University of Science and Technology (KAUST). J.H. Ram\'irez-Gonzalez and Ying Sun gratefully acknowledge this support. AMS thanks the University of Guanajuato for granting a sabbatical leave during which this work was completed.

\medskip

\noindent{\bf Funding Statement.}
No funding was received for this manuscript.

\medskip

\noindent{\bf Conflict of Interest Statement.}
The authors declare no conflict of interest.

\medskip

\noindent{\bf Data Availability Statement.}
The bat telemetry data used in this study are openly available at Movebank: \url{https://doi.org/10.5441/001/1.5d736bf0}. The code used for the inference is publicly available at \url{https://github.com/joseramirezgonzalez/f-wsfBm}. An interactive R Shiny application developed to explore the proposed stochastic models is available at \url{https://jhramgon.shinyapps.io/aplication_r_log_w/}.

\medskip

\noindent{\bf Code Availability.}
The code used in this study is publicly available at \url{https://github.com/joseramirezgonzalez/f-wsfBm}. An interactive R Shiny application for exploring the proposed stochastic models is available at \url{https://jhramgon.shinyapps.io/aplication_r_log_w/}.

\textbf{Corresponding author:} josehermenegildo.ramirezgonzalez@kaust.edu.sa.
\clearpage

\let\thm\suppthm
\let\endthm\endsuppthm
\let\lem\supplem
\let\endlem\endsupplem
\let\prop\suppprop
\let\endprop\endsuppprop
\let\cor\suppcor
\let\endcor\endsuppcor
\let\definition\suppdefinition
\let\enddefinition\endsuppdefinition
\let\example\suppexample
\let\endexample\endsuppexample
\let\rem\supprem
\let\endrem\endsupprem

\appendix
\renewcommand{\thesection}{\Alph{section}}

\section{Simulation Study}\label{sa}

\subsection{Simulation of sample paths}

\begin{figure}[H]
\centering
\subfigure[$f(u):=u^{-0.93}$]{\includegraphics[width=53mm,height=45mm]{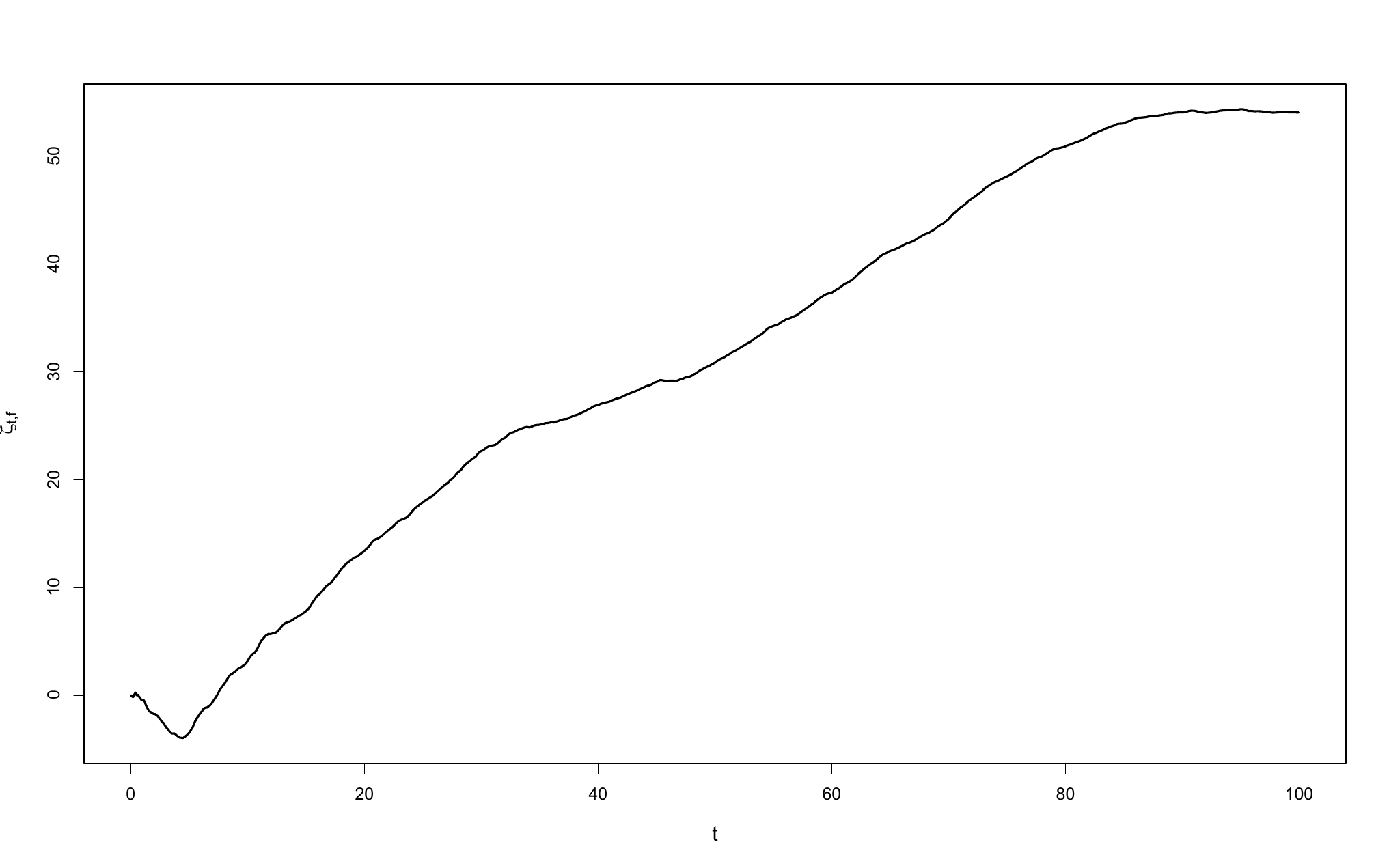}}
\subfigure[$f(u):=e^{-0.8u}$]{\includegraphics[width=53mm,height=45mm]{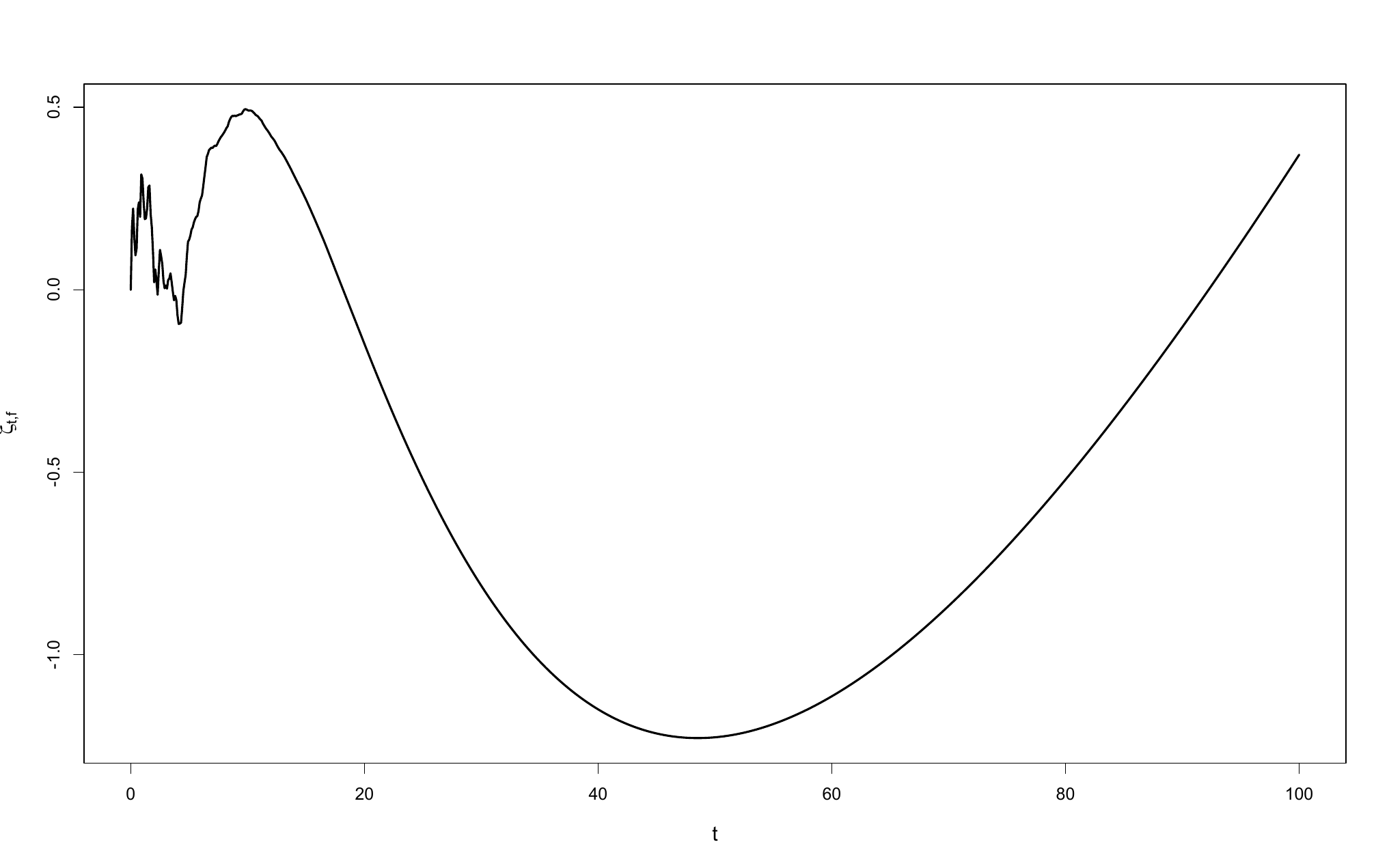}}
\subfigure[$f\equiv 1$]{\includegraphics[width=53mm,height=45mm]{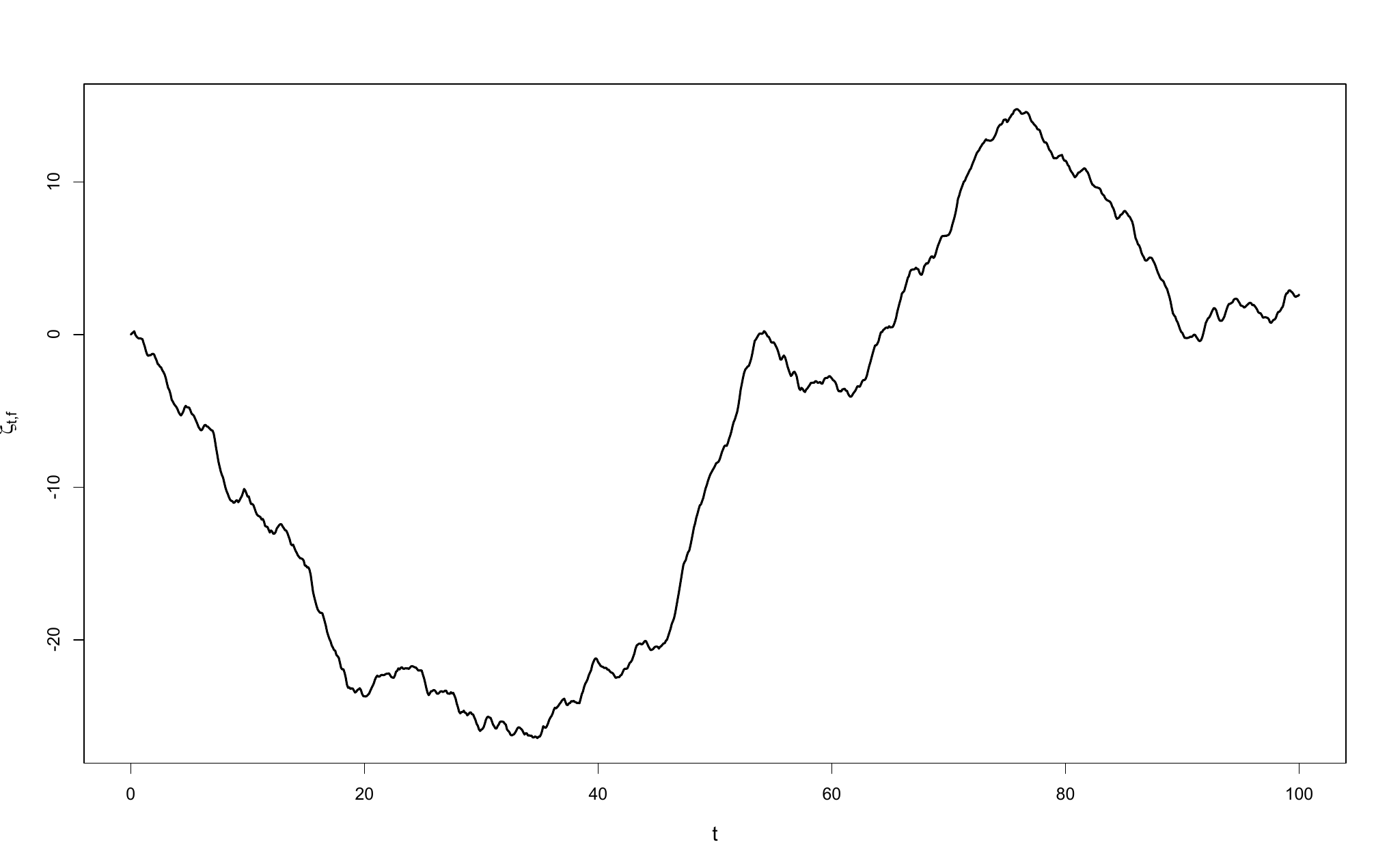}}
\subfigure[$f(u):=u^{10.6}$]{\includegraphics[width=53mm,height=45mm]{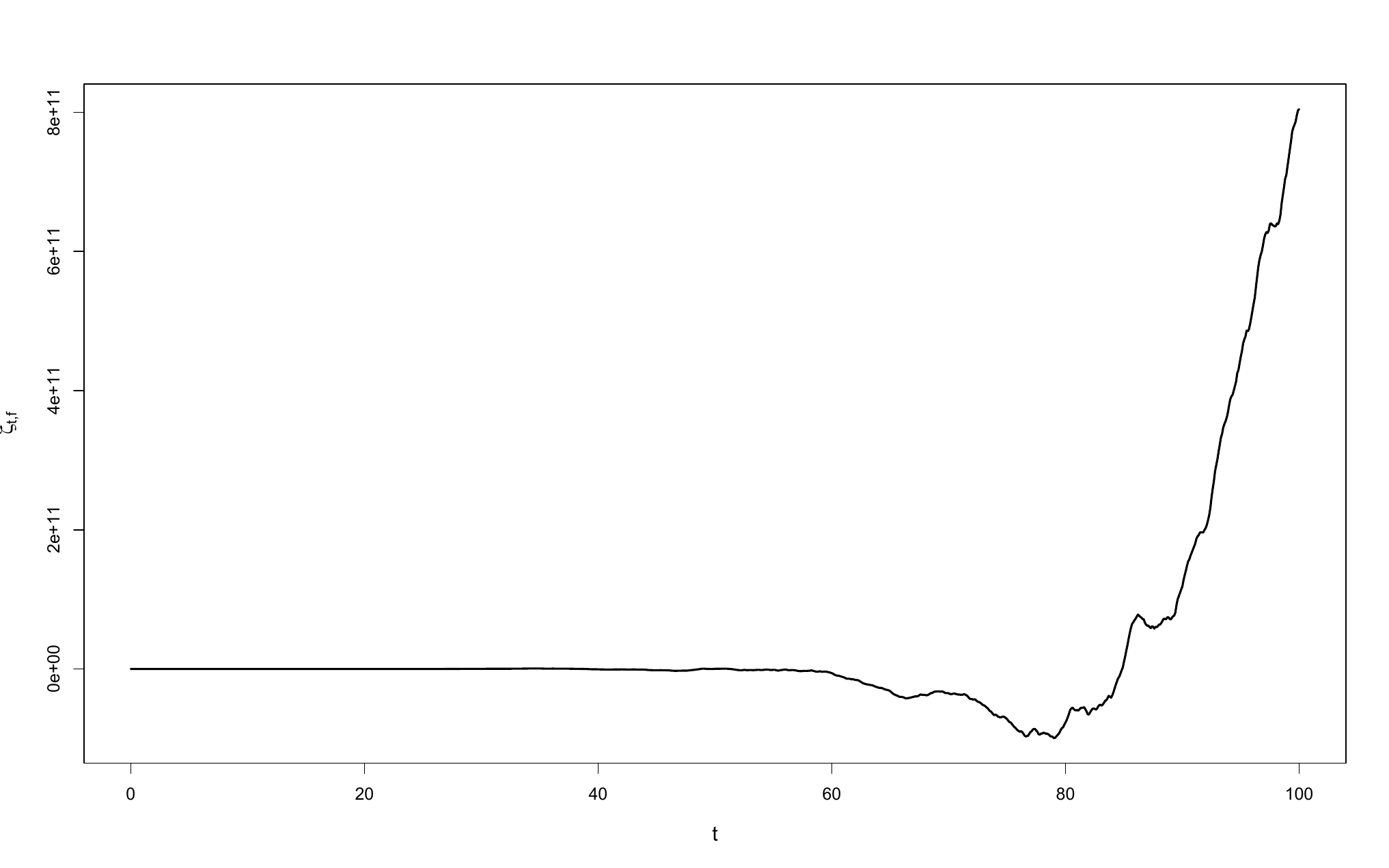}}
\subfigure[$f(u)=\cos\!\left(\frac{u}{1+u}\right)$]{\includegraphics[width=53mm,height=45mm]{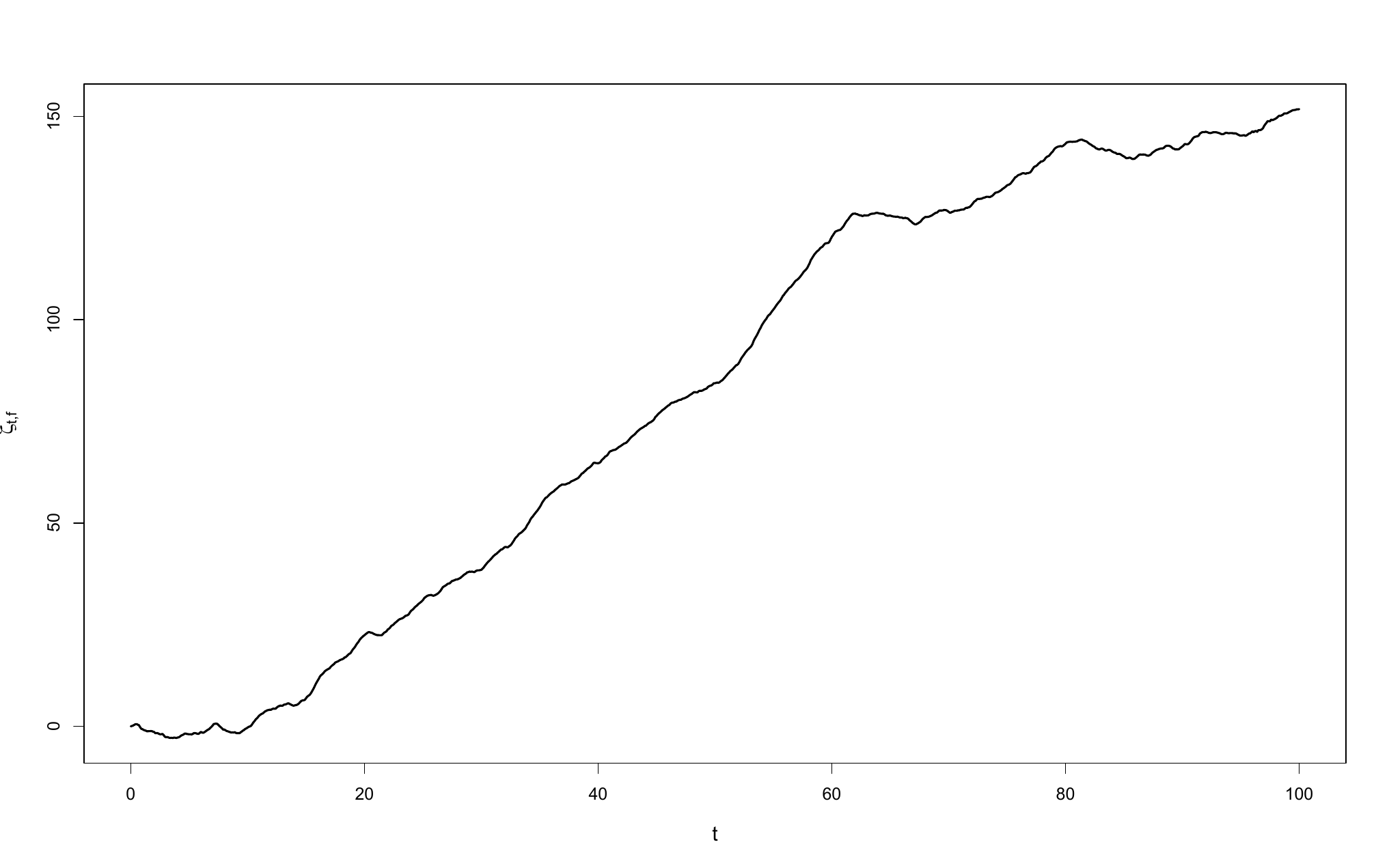}}
\subfigure[$f(u)=1-\cos(u)$]{\includegraphics[width=53mm,height=45mm]{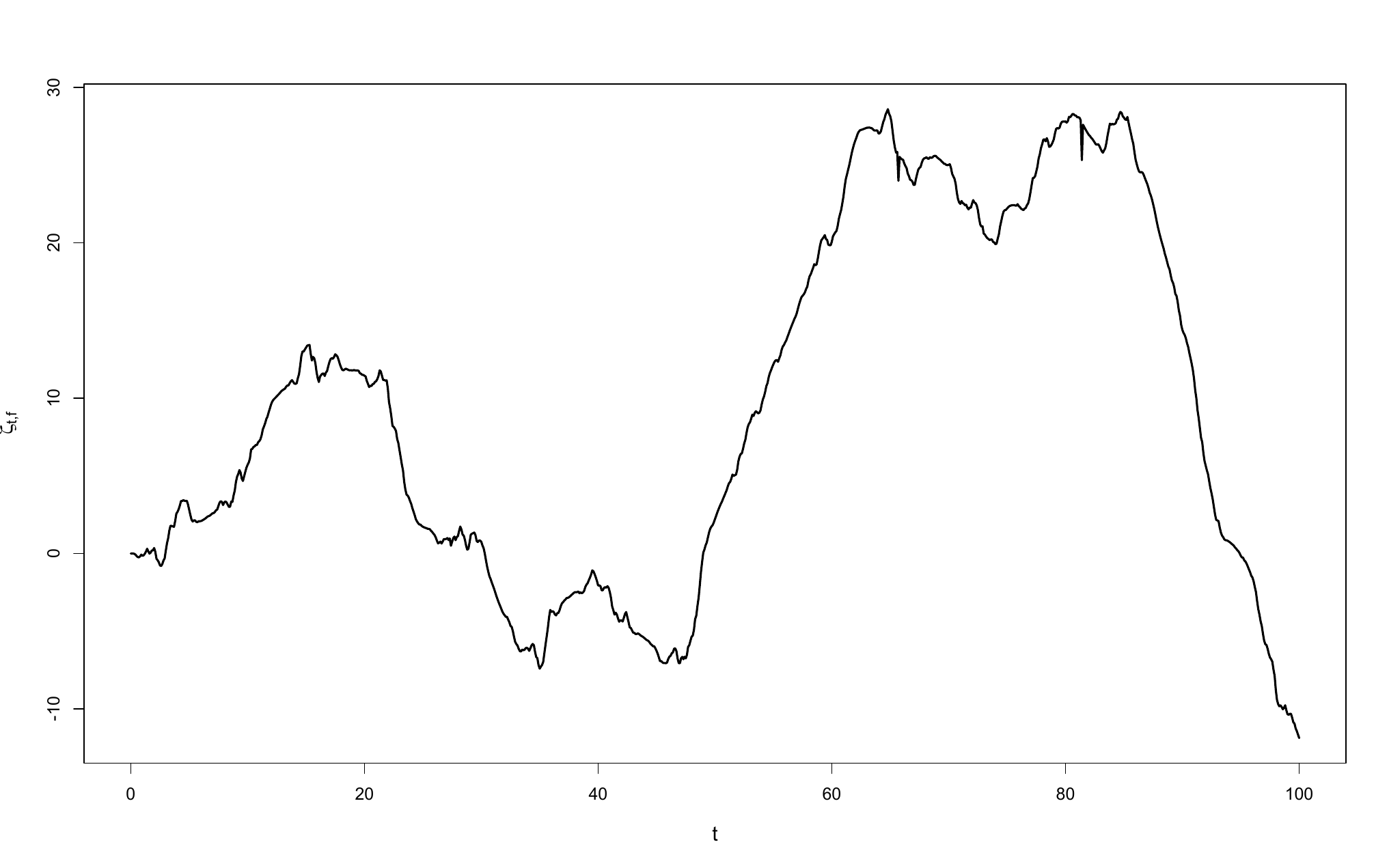}}
\subfigure[$f(u)=1-\sin(u)$]{\includegraphics[width=53mm,height=45mm]{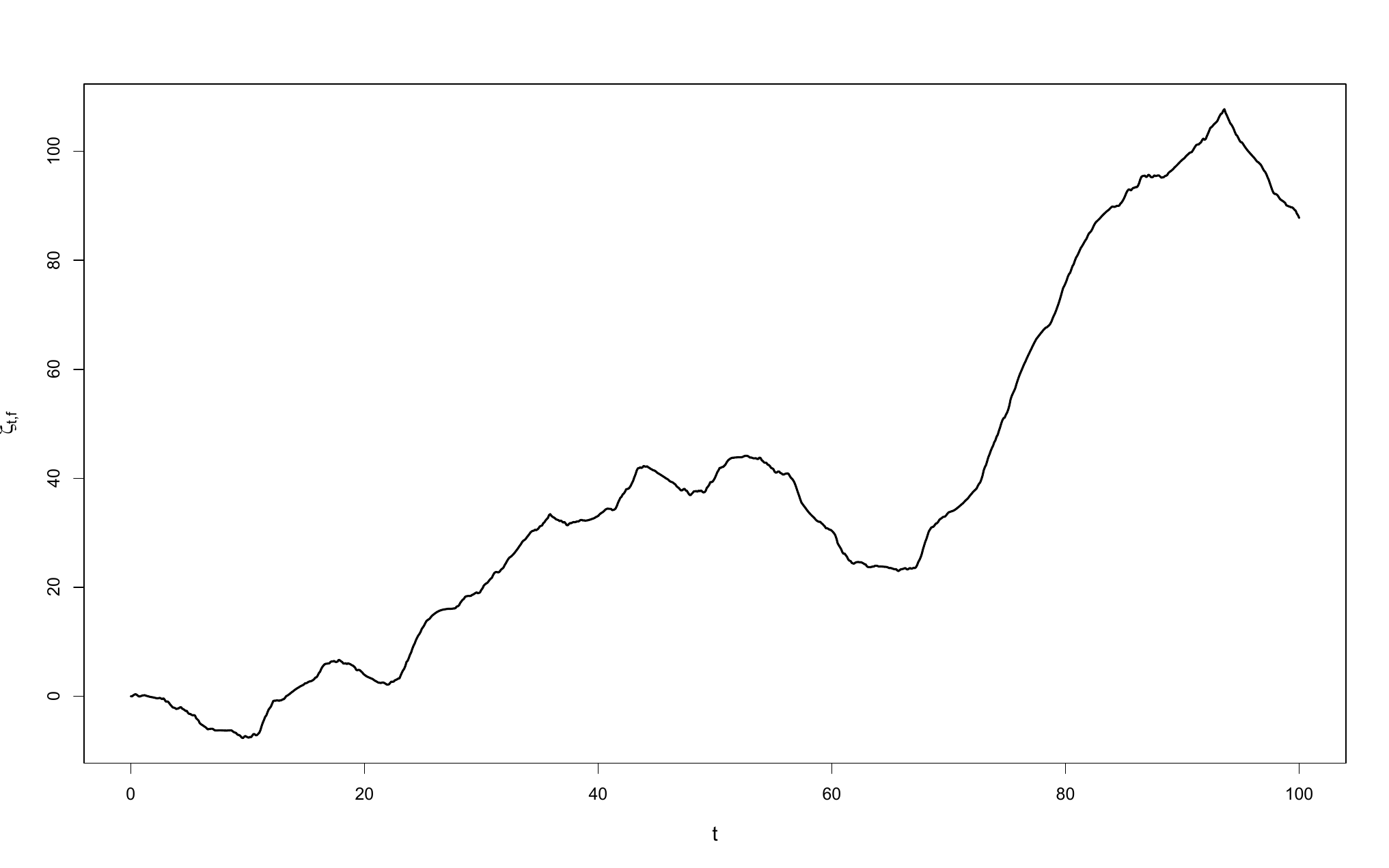}}
\subfigure[$f(u)=(1-\sin(u))(1-\cos(u))$]{\includegraphics[width=53mm,height=45mm]{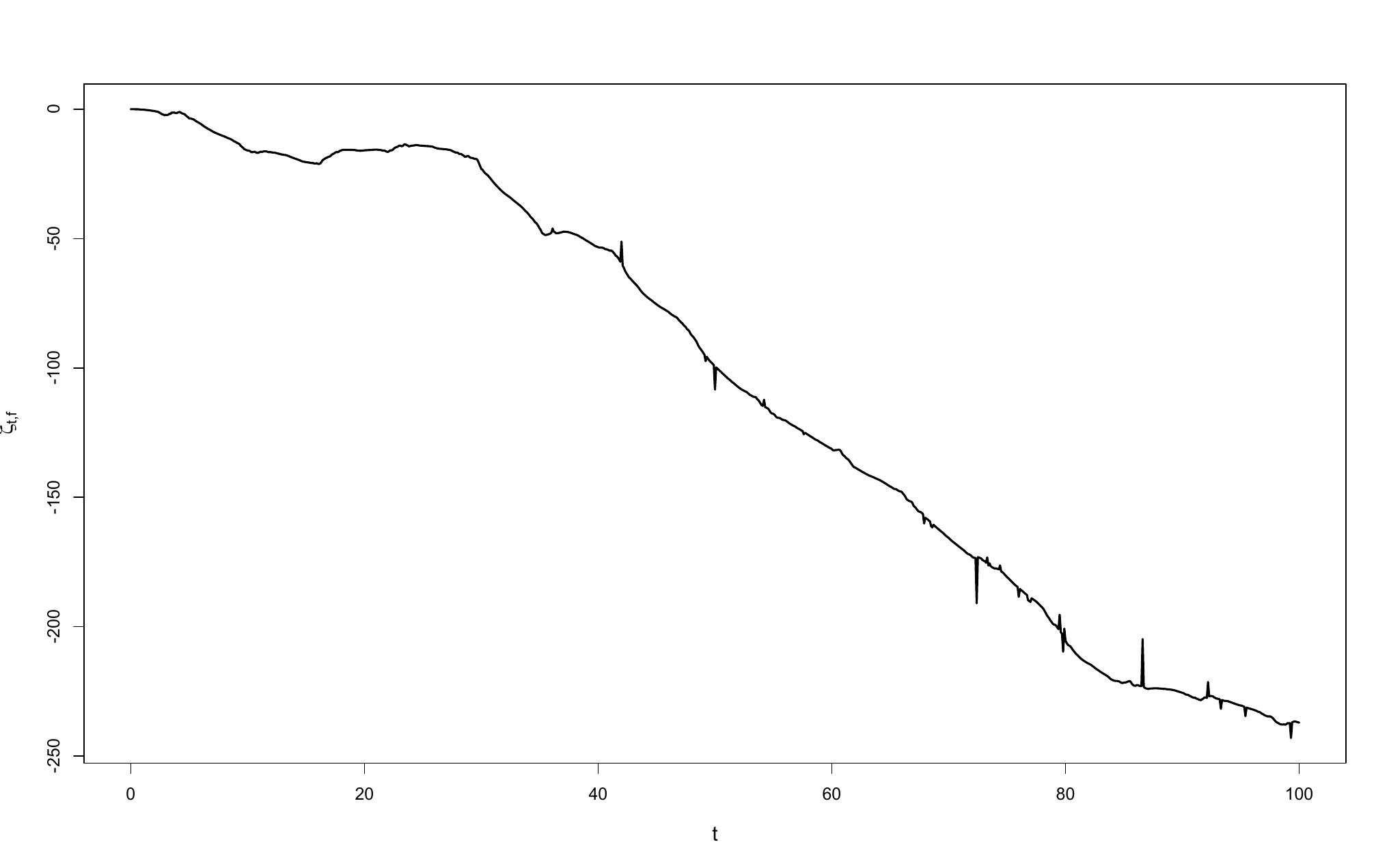}}
\subfigure[$f(u)=\frac{u}{1+u}$]{\includegraphics[width=53mm,height=45mm]{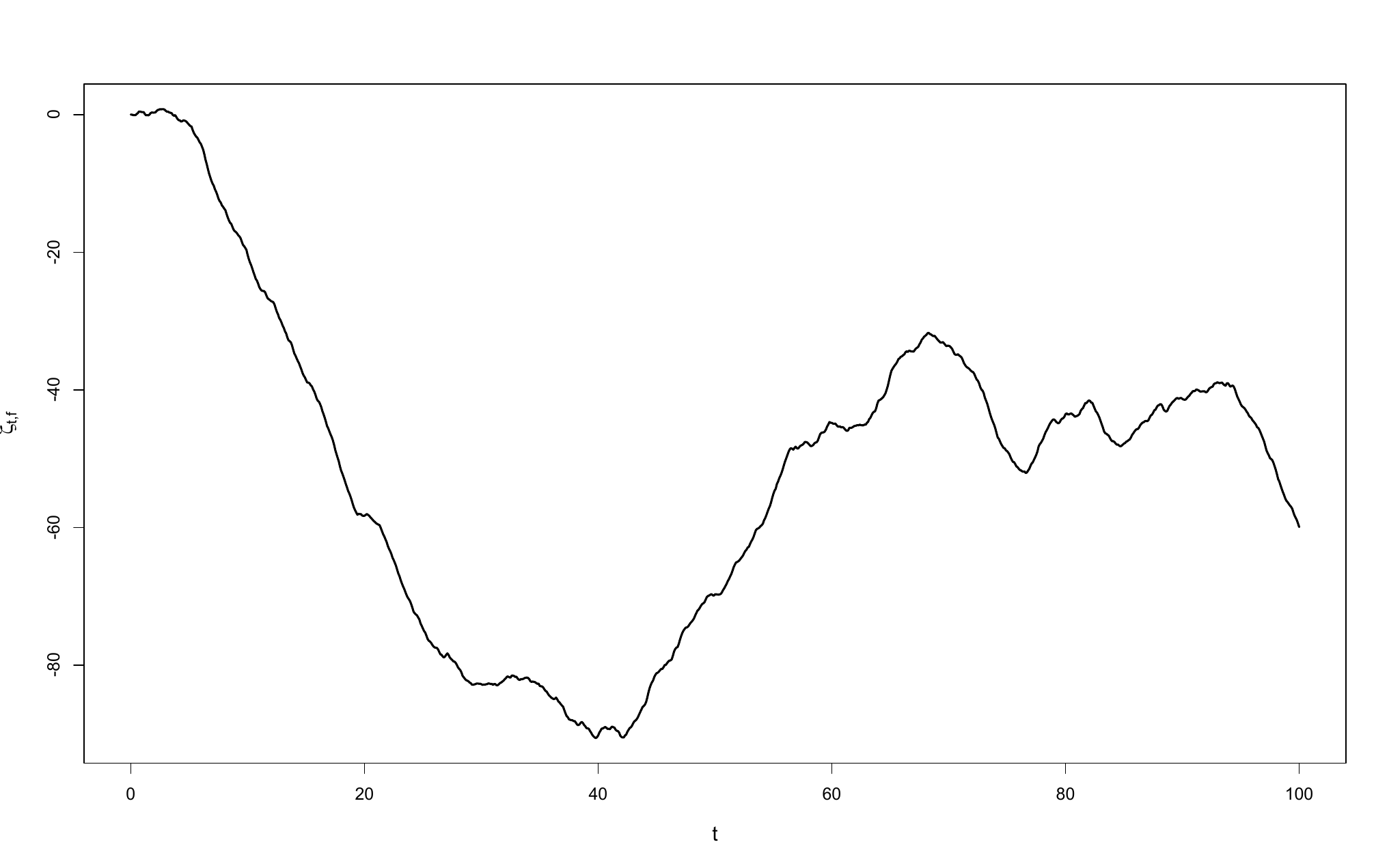}}
\subfigure[$f(u)=(1-\cos(u))e^{-0.5u}$]{\includegraphics[width=53mm,height=45mm]{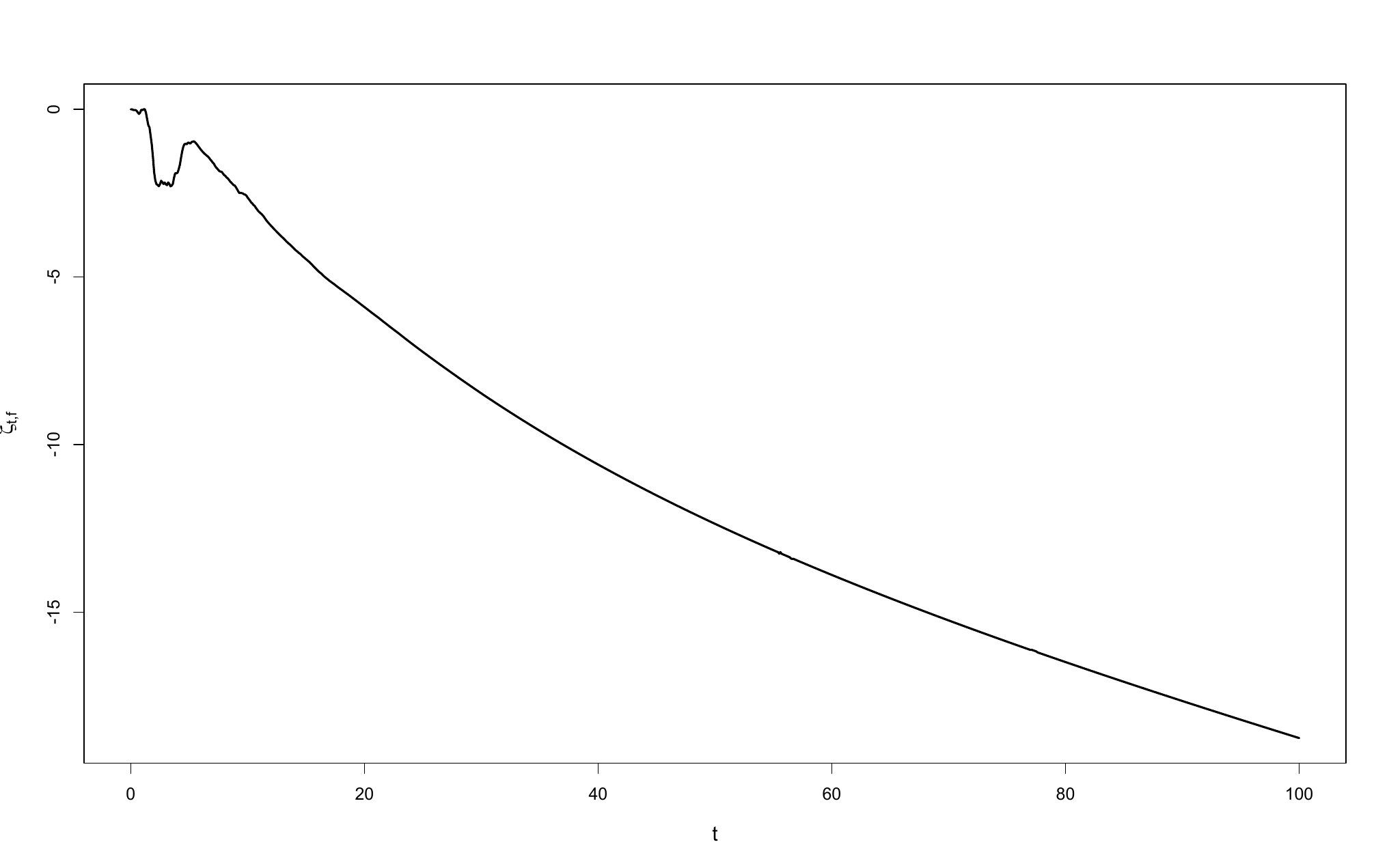}}
\subfigure[$f(u)=(1-\sin(u))e^{-0.5u}$]{\includegraphics[width=53mm,height=45mm]{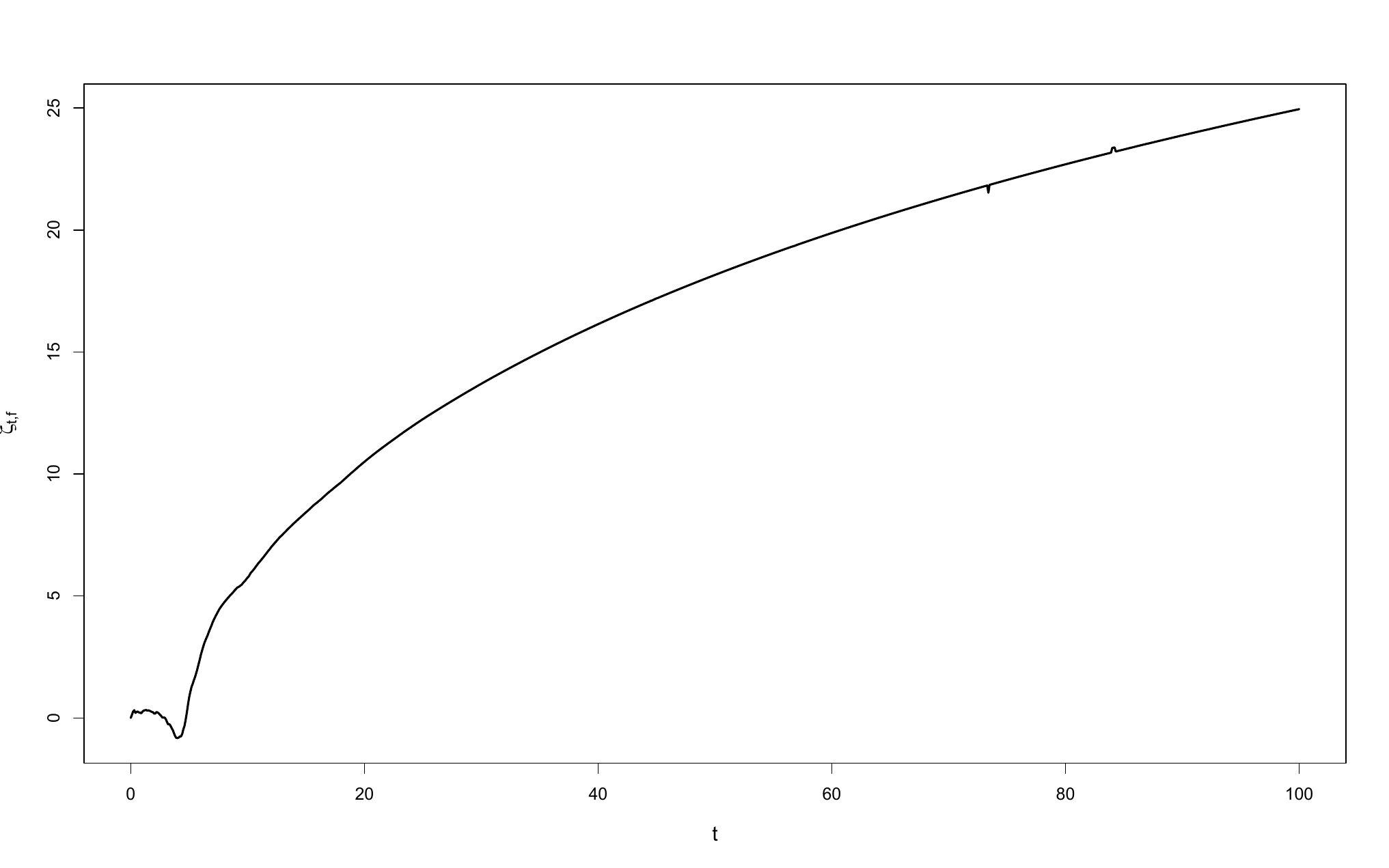}}
\subfigure[$f(u)=1-\cos(u)\sin(u)$]{\includegraphics[width=53mm,height=45mm]{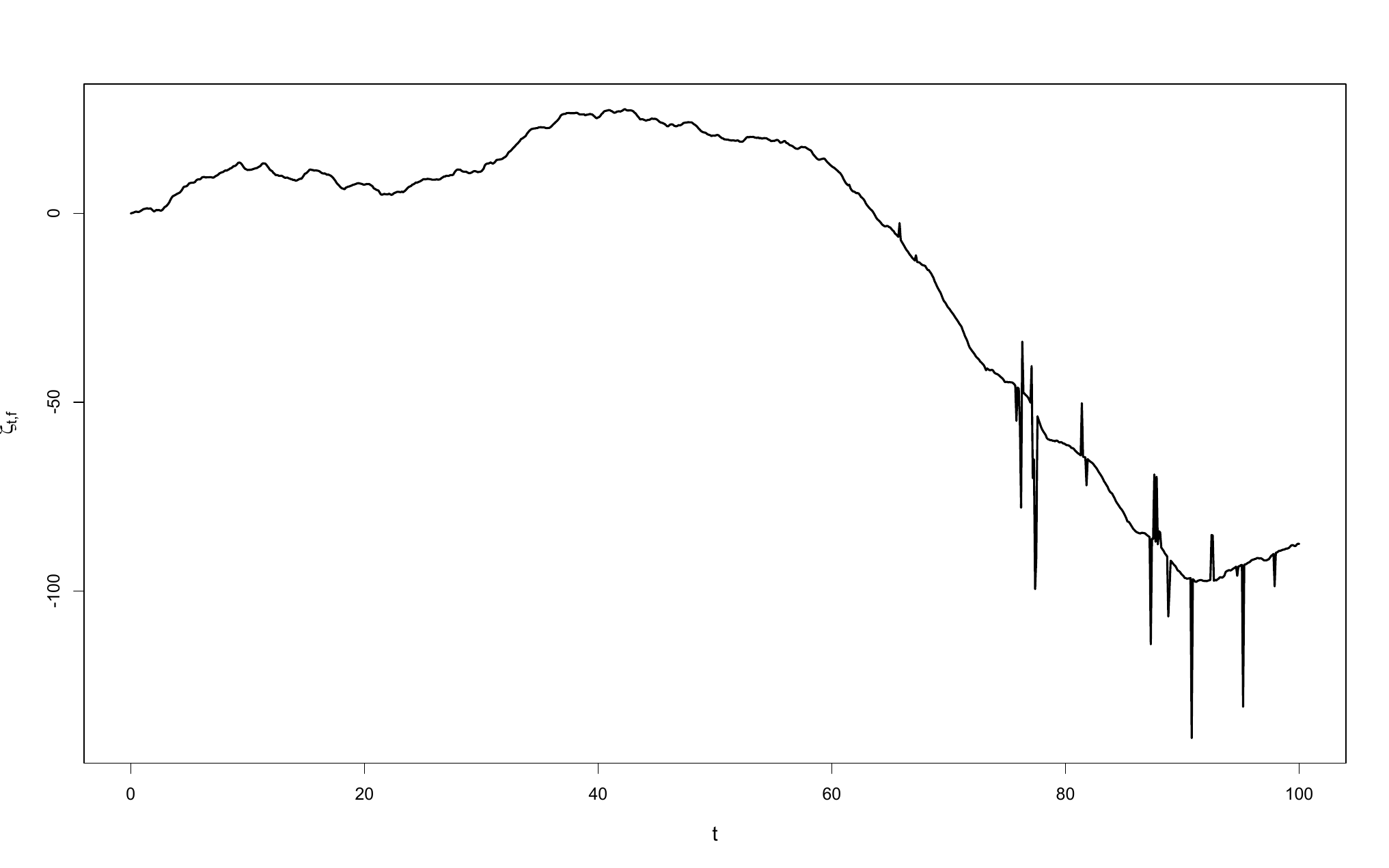}}
\caption{Simulated sample paths of the process $(\zeta_{t,f})_{t\geq0}$ for several representative choices of $f$. The simulations were performed on a uniform grid of 1000 time points over $[0,100]$, with initial condition $\zeta_{0,f}=0$.}
\label{F1}
\end{figure}

Let $\bar{t}=(t_0,t_1,\dots,t_n)$ denote a discrete set of time points. Simulation and prediction for the process $\zeta_f:=(\zeta_{t,f})_{t\geq0}$ are based on the finite-dimensional Gaussian distributions associated with the covariance function $K_f(s,t)$ defined in equation~(1). Figure~\ref{F1} presents simulated sample paths of $\zeta_f$ for several representative choices of the function $f$ satisfying condition~(2). In all cases, we consider a uniform discretization $t_i=i\Delta$ with $\Delta=1/10$ over the interval $[0,T]$, where $T=100$. The figure shows that the qualitative behavior of the sample paths substantially depends on the choice of $f$.
\newpage

\subsection{Inference from simulated data}

In this subsection, we study likelihood-based inference for the finite-dimensional distributions of the process \(\zeta_f\) under the exponential class
\[
\mathcal{C}=\left\{f_{\sigma^2,\beta}(u)=\sigma^2 e^{-\beta u}:\sigma^2>0,\ \beta\ge0\right\}.
\]
This is the same parametric class used in the analysis of the real data in the main text. Maximum likelihood estimation was implemented in \texttt{R}; the corresponding code is available at \url{https://github.com/joseramirezgonzalez/f-wsfBm}.

Let \(T>0\), \(n\in\mathbb{N}\), and \(\bar{t}=(t_0:=0,\ldots,t_n:=T)\) be a partition of \([0,T]\) into \(n\) subintervals, with \(\Delta_k=t_k-t_{k-1}\) for \(k=1,\ldots,n\). Throughout this subsection, we consider a uniform partition, so that \(\Delta_k=T/n\) for all \(k=1,\ldots,n\).

For the simulation experiment, we set \(T=120\), \(n=400\), \(\sigma^2=32.8409\), and \(\beta=0.044\), and generated data from the finite-dimensional distributions induced by \(K_f(s,t)\) on the grid \(\bar{t}\). The first \(90\%\) of the simulated observations were used for parameter estimation, while the remaining \(10\%\), corresponding to the final \(40\) observations, were reserved for prediction. The likelihood was profiled as described in Section~3.2 of the main text: for each fixed \(\beta\), the scale parameter is estimated by
\(
\widehat{\sigma}^{\,2}(\beta)=\frac{Q_\beta}{n},
\)
and the resulting profile log-likelihood is maximized over \(\beta\).

The maximum likelihood estimates obtained from the estimation segment were \(\widehat{\sigma}^{\,2}=33.0625\) and \(\widehat{\beta}=0.04366\). The corresponding approximate \(95\%\) confidence intervals, obtained from the local Gaussian approximation to the profile likelihoods, were \((28.2325,37.8925)\) for \(\sigma^2\) and \((0.03921,0.04810)\) for \(\beta\). Thus, in this simulated trajectory, the true parameter values used to generate the data lie within the corresponding approximate confidence intervals.

\begin{figure}[H]
\centering
\subfigure[Relative profile likelihood for \(\sigma^2\) under \(\mathcal{C}\).]{\includegraphics[width=0.92\textwidth]{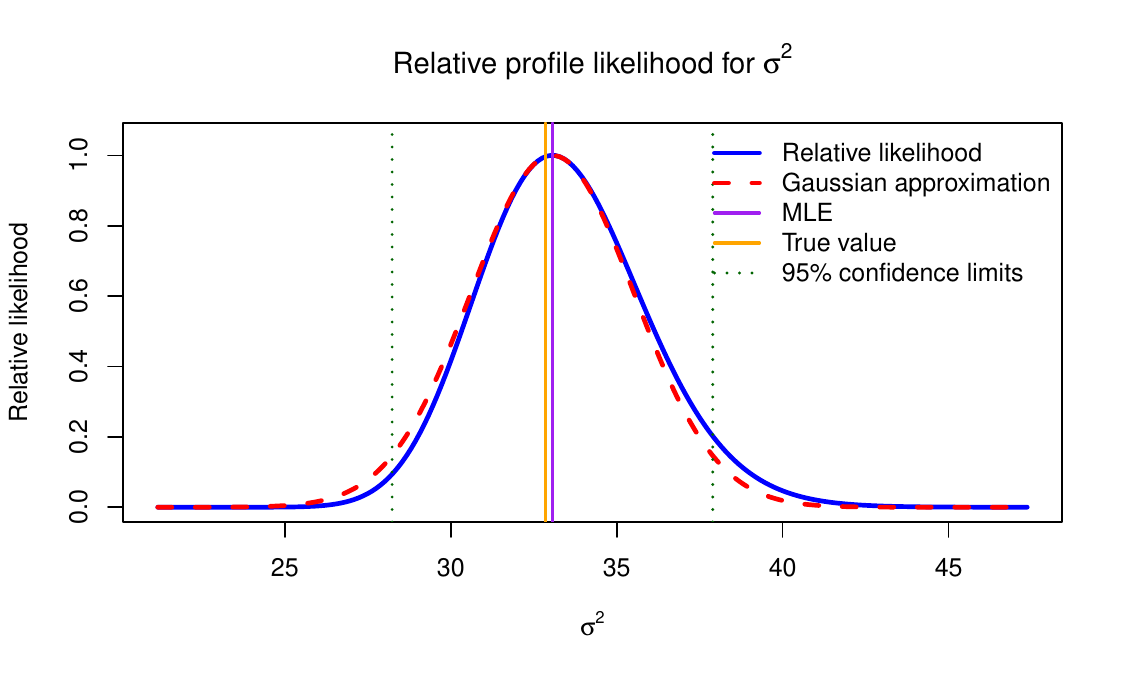}}
\subfigure[Relative profile likelihood for \(\beta\) under \(\mathcal{C}\).]{\includegraphics[width=0.92\textwidth]{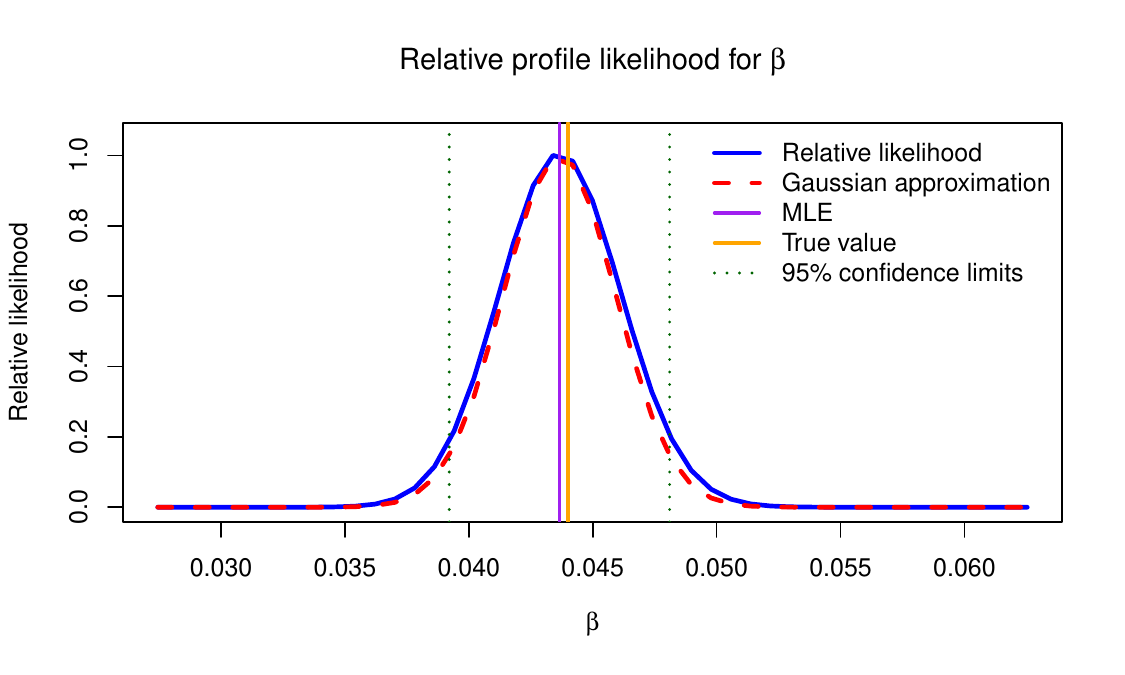}}
\caption{Relative profile likelihoods for \(\sigma^2\) and \(\beta\) under the exponential class \(\mathcal{C}=\{f_{\sigma^2,\beta}(u)=\sigma^2e^{-\beta u}:\sigma^2>0,\ \beta\ge0\}\), together with their local Gaussian approximations. The estimation stage is based on the first \(90\%\) of a simulated trajectory generated with \(T=120\), \(n=400\), \(\sigma^2=32.8409\), and \(\beta=0.044\). Each relative likelihood is normalized to attain a maximum value of one.}
\label{Ratiolikelihood_exp}
\end{figure}

Figure~\ref{Ratiolikelihood_exp} displays the relative profile likelihoods for \(\sigma^2\) and \(\beta\), together with their local Gaussian approximations. Each relative likelihood is normalized to attain a maximum value of one. The true parameter values are indicated in the plots, together with the maximum likelihood estimates and the corresponding approximate \(95\%\) confidence limits.

\begin{figure}[H]
\centering
\includegraphics[width=0.92\textwidth]{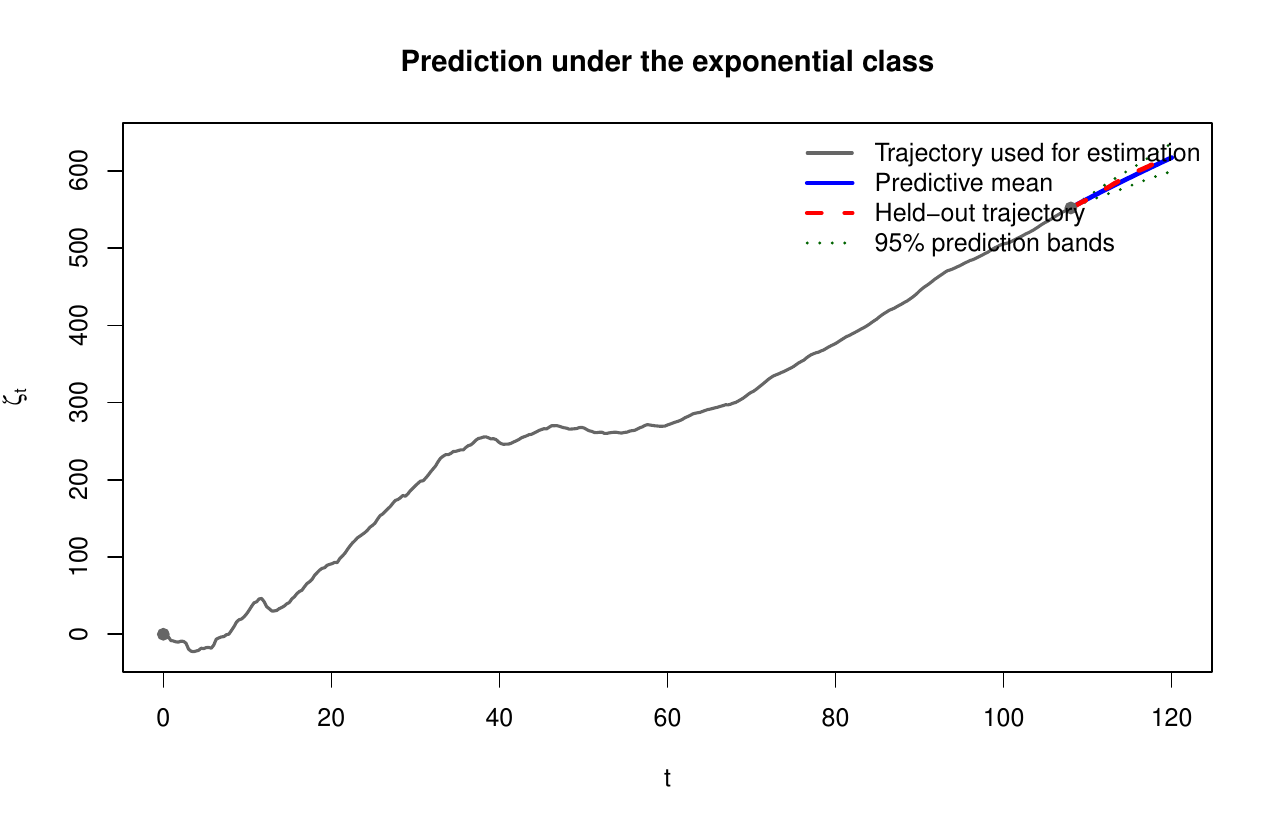}
\caption{Prediction of the final \(40\) time points under the exponential class \(\mathcal{C}\), based on the maximum likelihood estimates obtained from the first \(90\%\) of the simulated observations. The predictive mean and the prediction bands are computed from 1000 conditional simulations.}
\label{Ratiolikelihood_exp_pred2}
\end{figure}

Prediction was carried out for the remaining \(10\%\) of the time points, corresponding to the final \(40\) observations, which were not used for parameter estimation. Figure~\ref{Ratiolikelihood_exp_pred2} shows the predicted continuation of the simulated trajectory. We generated 1000 conditional simulations under the fitted model and used the corresponding predictive mean trajectory for comparison with these reserved observations. The prediction error, computed on this reserved segment using the relative mean squared error, was \(1.7484\times10^{-5}\). This calculation illustrates the implementation of likelihood-based estimation and conditional prediction under the class \(\mathcal{C}\).

\section{Theoretical results}

This appendix contains detailed proofs of the theoretical results stated in the main text, together with additional  results that complement the analysis of the proposed process. In particular, the material presented here covers conditions ensuring the positive definiteness of the covariance kernel, path properties, and long-range dependence, as well as further results on self-similarity, short-time asymptotics, and continuity properties.

\subsection*{Proof of Theorem 1} 
It is straightforward to verify that \( K_f(s,t) \) is symmetric; therefore, it suffices to prove that it is nonnegative definite. Let \( a > 0 \) and \( s, t \geq 0 \). We first prove that
\small
\begin{eqnarray*}
&&\lefteqn{Q_a(t,s) :=} \\
&&2\left[(s+t)\log(s+t) - (t+a)\log(t+a) - (s+a)\log(s+a) + 2a\log(2a)\right] \cdot 1_{\{a \leq s\}} 1_{\{a \leq t\}},
\end{eqnarray*}
\normalsize
is nonnegative definite.  A direct computation shows that
\small
\begin{equation}\label{PDQ_a}
 Q_a(t,s) = \int_a^s \int_a^t \int_0^\infty 2 e^{-(r + r')x} \, dx \, dr \, dr' =\int_a^s \int_a^t h(r,r') \, dr \, dr'
\end{equation}
\normalsize
where \( h(r,r') := 2\int_0^\infty e^{-(r + r')x} \, dx = \frac{2}{r+r'} \). Moreover, given \( m \in \mathbb{N} \),  \( x_1, \ldots, x_m \in [0,\infty) \) and  \( y_1, \ldots, y_m \in \mathbb{R} \), we get
\begin{equation}\label{h_nnd}
   \sum_{i=1}^m \sum_{j=1}^m h(x_i, x_j) y_i y_j = 2\int_0^\infty \sum_{i=1}^m \sum_{j=1}^m e^{-(x_i + x_j)x} y_i y_j \, dx = 2\int_0^\infty \left( \sum_{i=1}^{m} y_i e^{-x_i x} \right)^2 dx \geq 0,
\end{equation}
which shows that \( h(r,r') \) is nonnegative definite.

Now, let \( b > a \) and let \( g \in L^2([a,b]) \). By a change in the order of integration, we obtain
\small
\begin{equation}\label{PDQ_a2}
\begin{aligned}
\int_{a}^{b}\int_{a}^{b} g(s)\, Q_a(s,t)\, g(t)\, ds\, dt 
&= \int_{a}^{b}\int_{a}^{b} g(s)g(t)\left( \int_a^s \int_a^t h(r,r') \, dr \, dr' \right) ds\, dt \\
&= \int_{a}^{b}\int_{a}^{b} \left( \int_r^b g(s)\, ds \right) h(r,r') \left( \int_{r'}^b g(t)\, dt \right) dr\, dr' \\
&= \int_{a}^{b}\int_{a}^{b} \psi_b(r)\, h(r,r')\, \psi_b(r') \, dr\, dr',
\end{aligned}
\end{equation}
\normalsize
where \( \psi_b(r) := \int_r^b g(s)\, ds \) for all \( r \in [a,b] \).

Since \( \psi_b(r) \) is continuous on \([a,b]\) and \( h(r,r') \) is continuous on \([a,b]^2 \), the integrand \( \psi_b(r)\, h(r,r')\, \psi_b(r') \) is continuous on \([a,b]^2 \). Therefore, by expression~\eqref{PDQ_a2} and the convergence of Riemann sums, we have that
\small
\begin{equation}\label{h_integral}
\begin{aligned}
&\int_{a}^{b}\int_{a}^{b}g(s)Q_a(s,t)g(t)\,ds\,dt =\int_{a}^{b}\int_{a}^{b}\psi_b(r)h(r,r')\psi_b(r^{'}) \, dr \, dr'\\
&= \lim_{n\to \infty}\sum_{i=1}^n \sum_{j=1}^n \Delta_n^{a,b}\psi_b(a + \left(i - \tfrac{1}{2}\right)\Delta_n^{a,b})h\left(a + \left(i - \tfrac{1}{2}\right)\Delta_n^{a,b},\, a + \left(j - \tfrac{1}{2}\right)\Delta_n^{a,b}\right)\psi_b(a + \left(j - \tfrac{1}{2}\right)\Delta_n^{a,b})\Delta_n^{a,b}\\
&=\lim_{n\to \infty}\sum_{i=1}^n \sum_{j=1}^n w_{n,i,2}^{a,b}h(w_{n,i,1}^{a,b},w_{n,j,1}^{a,b})w_{n,j,2}^{a,b}
\end{aligned}
\end{equation}
\normalsize
where $\Delta_n^{a,b}:=\frac{b - a}{n}$, $w_{n,i,1}^{a,b}:=a + \left(i - \tfrac{1}{2}\right)\Delta_n^{a,b}$, $w_{n,i,2}^{a,b}=\Delta_n^{a,b}\psi_b(a + \left(i - \tfrac{1}{2}\right)\Delta_n^{a,b})$ for $i=1,...,n$. Therefore, from (\ref{h_nnd}) and (\ref{h_integral}), we conclude that
\[
\int_{a}^{b}\int_{a}^{b} g(s)\, Q_a(s,t)\, g(t)\, ds\, dt \geq 0,
\]
and, according to Section~6.3 of \cite{Fuk10}, this implies that \(Q_a(\cdot,\cdot)\) is a nonnegative-definite kernel. Furthermore, since \( \lim_{a \downarrow 0} 2a\log(2a) = 0 \), we deduce that the pointwise limit
\begin{equation}\label{Q-Kernel}
Q(s,t) := \lim_{a \downarrow 0} Q_a(s,t) = 2\left[(s+t)\log(s+t) - t\log(t) - s\log(s)\right] \mathbf{1}_{\{s \geq 0\}} \mathbf{1}_{\{t \geq 0\}}, \quad s, t > 0,
\end{equation}
is also nonnegative definite. 


Assume now that the function \( f: [0,\infty) \to [0,\infty) \) satisfies the integrability condition
\begin{equation}\label{Integrability-condition-appendix}
\int_0^\delta f(u) \, du < \infty, \text{ for every } \delta > 0.
\end{equation}
Under this assumption, it is immediate to check that \( |K_f(s,t)| < \infty \) for all \( s,t \geq 0 \). Furthermore,
\begin{eqnarray*}
\sum_{i=1}^m \sum_{j=1}^m K_f(x_i, x_j) y_i y_j &=& \int_0^\infty f(r) \sum_{i=1}^m \sum_{j=1}^m Q(x_i - r, x_j - r) y_i 1_{\{0 \leq r \leq x_i\}} y_j 1_{\{0 \leq r \leq x_j\}} \, dr \\
&=& \int_0^\infty f(r) \sum_{i=1}^m \sum_{j=1}^m Q(x_i - r, x_j - r) z_i(r) z_j(r) \, dr,
\end{eqnarray*}
where \( z_i (r):= y_i 1_{\{0 \leq r \leq x_i\}} \). Since \( f \geq 0 \) and \(Q(\cdot,\cdot)\) is a nonnegative-definite kernel, the function \( K_f (s,t)\) is also nonnegative definite. This completes the proof.

\hfill \(\square\)
\begin{rem}
    As an alternative argument to show that \( Q(s,t) \) is nonnegative definite, let us consider the function  
\[
Q_b(s,t) := \frac{1}{1-b}\left[s^b + t^b - (s + t)^b\right], \quad s, t \geq 0,
\]
defined for \( b \in [0,1) \cup (1,2] \). As shown in equation (2.4) of \cite{BojParti}, \( Q_b(s,t) \) is a covariance function for each such \( b \). 
Moreover, using L’Hôpital’s rule, it is straightforward to verify that
\[
\lim_{b \to 1} 2 Q_b(s,t) = Q(s,t),
\]
for all \( s, t \geq 0 \). Since the set of covariance functions is closed under pointwise limits, it follows that \( Q(s,t) \) is nonnegative definite.
\end{rem}

\subsection{Sample path properties}\label{memoryp}
We now discuss several path properties of the process \(\zeta\), including total variation, quadratic variation, and modulus of continuity. We also show that \(\zeta\) is not a semimartingale; consequently, the standard It\^o stochastic calculus for semimartingales does not apply. For background on this topic, particularly in connection with processes of finite quadratic variation, see, for example, \cite{ruso}.

By setting \(f\equiv 1\) in \((1)\), we obtain
\begin{equation}\label{logcov}
   {K}_1(s,t)= -\Bigg(s^{2}\log(s)+t^{2}\log(t)-\frac{1}{2}\Bigg[(t+s)^{2}\log(s+t)+|s-t|^{2}\log(|s-t|)\Bigg]\Bigg).
\end{equation}
In \cite{BGT2007} the process $\eta=(\eta_t)_{t\geq 0}$ with covariance \eqref{logcov} was considered.

For the results below, we introduce the following two classes of weight functions:
\[
\mathcal{C}_1:=\{f_{\alpha}(u)=u^{\alpha}: u>0,\ \alpha>-1\},
\qquad
\mathcal{C}_2:=\{f_{\alpha}(u)=e^{\alpha u}: u>0,\ \alpha\in\mathbb{R}\}.
\]
Throughout this subsection, when we write that \(f\in\mathcal{C}_1\) or \(f\in\mathcal{C}_2\), we mean that the covariance kernel \(K_f\) in \((1)\) is obtained by substituting the corresponding weight function into the general formula.

In the following result we establish the self-similarity property of the process $\zeta$ for weights in the class $\mathcal{C}_1$.

\begin{prop}\label{auto}
For \(f(u)=u^\alpha\in\mathcal{C}_1\), the process \(\zeta\) is self-similar with index \((\alpha+2)/2\). More precisely, for every \(c>0\),
\[
(\zeta_{ct})_{t\geq 0}\overset{d}= \bigl(c^{(\alpha+2)/2}\zeta_t\bigr)_{t\geq 0}.
\]
\end{prop}

\begin{proof}
Let \(\alpha>-1\). For every \(c>0\), a change of variables together with the identity \(\log(ab)=\log(a)+\log(b)\) shows that
\[
K_f(cs,ct)=c^{\alpha+2}K_f(s,t), \qquad s,t>0,
\]
where \(f(u)=u^\alpha\). The conclusion follows from the Gaussianity of \(\zeta\). \qed
\end{proof}

\begin{prop}
\label{prop:completely-monotone-approximation}
Let \(f:(0,\infty)\to[0,\infty)\) be completely monotone, that is,
\(f\in C^\infty(0,\infty)\) and
\[
(-1)^k f^{(k)}(u)\geq0,
\qquad
u>0,\quad k=0,1,2,\ldots,
\]
and assume that
\[
\int_0^T f(u)\,du<\infty
\qquad
\text{for every }T>0.
\]
Define the family of finite nonnegative exponential sums
\[
\mathcal{E}
:=
\left\{
u\longmapsto\sum_{j=1}^{m}a_j e^{-b_j u}
:
m\in\mathbb N,\quad a_j\geq0,\quad b_j\geq0
\right\}.
\]
Then the following statements hold.

\begin{enumerate}
\item[(i)]
There exists a sequence
\[
f_n(u)
=
\sum_{j=1}^{m_n}a_{n,j}e^{-b_{n,j}u}
\in\mathcal{E}
\]
such that
\[
\int_0^T|f_n(u)-f(u)|\,du
\longrightarrow0
\]
for every \(T>0\).

\item[(ii)]
Let \(K_h\) denote the covariance kernel associated in Theorem~1 with
an admissible weight \(h\). Then, for every \(T>0\),
\[
\sup_{0\leq s,t\leq T}
|K_{f_n}(s,t)-K_f(s,t)|
\longrightarrow0.
\]

\item[(iii)]
For each \(n\), let
\(X^{(n,1)},\ldots,X^{(n,m_n)}\) be independent centered Gaussian
processes satisfying
\[
\operatorname{Cov}
\left(
X^{(n,j)}_s,X^{(n,j)}_t
\right)
=
K_{e^{-b_{n,j}\,\cdot}}(s,t),
\]
and define
\[
Y^{(n)}_t
:=
\sum_{j=1}^{m_n}
\sqrt{a_{n,j}}\,X^{(n,j)}_t.
\]
If \(X^f\) denotes a centered Gaussian process with covariance \(K_f\),
then
\[
Y^{(n)}
\overset{\mathrm{f.d.d.}}{\longrightarrow}
X^f.
\]
Thus, every process generated by a completely monotone and locally
integrable weight is an f.d.d.\ limit of finite sums of independent
Gaussian processes generated by exponential weights.
\end{enumerate}
\end{prop}

\begin{proof}
By Bernstein's representation theorem
\citep[Theorem~1.4]{SchillingSongVondracek2010}, there exists a unique
positive Borel measure \(\mu\) on \([0,\infty)\) such that
\[
f(u)
=
\int_{[0,\infty)}e^{-bu}\,\mu(db),
\qquad u>0.
\]
Fix \(T>0\) and set
\[
w_T(b)
:=
\int_0^T e^{-bu}\,du
=
\begin{cases}
\dfrac{1-e^{-bT}}{b}, & b>0,\\[5pt]
T, & b=0.
\end{cases}
\]
Tonelli's theorem gives
\[
\int_{[0,\infty)}w_T(b)\,\mu(db)
=
\int_0^T f(u)\,du
<\infty.
\]

Given \(\varepsilon>0\), choose \(M>0\) such that
\[
\int_{(M,\infty)}w_T(b)\,\mu(db)<\frac{\varepsilon}{2}.
\]
Since \(w_T(b)\geq w_T(M)>0\) on \([0,M]\), one has
\(\mu([0,M])<\infty\). Moreover, for \(b,c\geq0\),
\[
\int_0^T|e^{-bu}-e^{-cu}|\,du
\leq
\frac{T^2}{2}|b-c|.
\]
Partition \([0,M]\) into disjoint intervals
\(I_1,\ldots,I_m\) of diameter at most \(\delta\), choose
\(b_j\in I_j\), and set
\[
a_j:=\mu(I_j),
\qquad
g_{\delta,M}(u)
:=
\sum_{j=1}^m a_j e^{-b_j u}.
\]
Then
\[
\begin{aligned}
\|f-g_{\delta,M}\|_{L^1(0,T)}
&\leq
\int_{(M,\infty)}w_T(b)\,\mu(db)
+
\frac{T^2}{2}\delta\,\mu([0,M]).
\end{aligned}
\]
Choosing \(\delta\) sufficiently small makes the right-hand side less
than \(\varepsilon\). Thus, \(\mathcal E\) is dense in \(L^1(0,T)\)
within this class of weights. For each \(n\in\mathbb N\), choose
\(f_n\in\mathcal E\) such that
\[
\int_0^n|f_n(u)-f(u)|\,du<\frac1n.
\]
Then, for every fixed \(T>0\) and all sufficiently large \(n\),
\[
\int_0^T|f_n(u)-f(u)|\,du
\leq\frac1n,
\]
which proves~(i).

For~(ii), write
\[
K_h(s,t)
=
2\int_0^{s\wedge t}h(r)H_{s,t}(r)\,dr,
\]
where
\[
H_{s,t}(r)
:=
(s+t-2r)\log(s+t-2r)
-(s-r)\log(s-r)
-(t-r)\log(t-r),
\]
with \(0\log0:=0\). For fixed \(T>0\), the function
\((s,t,r)\mapsto H_{s,t}(r)\) is continuous on
\[
D_T
:=
\left\{
(s,t,r)\in[0,T]^3:
0\leq r\leq s\wedge t
\right\},
\]
and hence
\[
C_T
:=
\sup_{(s,t,r)\in D_T}|H_{s,t}(r)|
<\infty.
\]
Therefore,
\[
\begin{aligned}
|K_{f_n}(s,t)-K_f(s,t)|
&\leq
2C_T\int_0^T|f_n(r)-f(r)|\,dr
\end{aligned}
\]
for \(0\leq s,t\leq T\). Taking the supremum and applying~(i)
proves~(ii).

Finally, independence and the linearity of \(h\mapsto K_h\) yield
\[
\begin{aligned}
\operatorname{Cov}(Y^{(n)}_s,Y^{(n)}_t)
&=
\sum_{j=1}^{m_n}
a_{n,j}K_{e^{-b_{n,j}\,\cdot}}(s,t)
=
K_{f_n}(s,t).
\end{aligned}
\]
For fixed \(t_1,\ldots,t_k\geq0\), let
\[
(\Sigma_n)_{i\ell}
:=
K_{f_n}(t_i,t_\ell),
\qquad
\Sigma_{i\ell}
:=
K_f(t_i,t_\ell).
\]
By~(ii), \(\Sigma_n\to\Sigma\) entrywise. Hence, for every
\(\theta\in\mathbb R^k\),
\[
\exp\left\{-\frac12\theta^\top\Sigma_n\theta\right\}
\longrightarrow
\exp\left\{-\frac12\theta^\top\Sigma\theta\right\}.
\]
Convergence of the corresponding Gaussian characteristic functions
proves~(iii).
\end{proof}

\begin{rem}
As a particular example, let
\[
f(u)=a u^{-\gamma},
\qquad
a\geq0,\quad 0<\gamma<1.
\]
Then
\[
(-1)^k f^{(k)}(u)
=
a\frac{\Gamma(\gamma+k)}{\Gamma(\gamma)}
u^{-\gamma-k}
\geq0
\]
for every \(k\geq0\), and
\[
\int_0^T f(u)\,du
=
a\frac{T^{1-\gamma}}{1-\gamma}
<\infty.
\]
Moreover,
\[
a u^{-\gamma}
=
\frac{a}{\Gamma(\gamma)}
\int_0^\infty b^{\gamma-1}e^{-bu}\,db.
\]
Thus, power-law weights with \(0<\gamma<1\) satisfy all the conclusions
of Proposition~\ref{prop:completely-monotone-approximation}.
\end{rem}

\begin{thm} \label{not-semimartingale}
    For \(f \in \mathcal{C}_1 \cup \mathcal{C}_2\), the process $\zeta$ is not a semimartingale. 
\end{thm}
The proof of Theorem~\ref{not-semimartingale} is based on showing that the process \( \zeta \) has infinite variation and zero quadratic variation on the interval \([0,1]\).

The following result concerns the modulus of continuity of \( \zeta \).
\begin{thm}\label{teo3}
(i)    Suppose that there exists a positive constant $c_{f}$ such that $0\leq f\leq c_f$ on $[0,T]$ for some $T>0$. Then,  for any $\kappa\in(0,1)$ there exists a positive constant $c_{f,\kappa}(T)$ such that
\begin{equation*}
        \E(\zeta_{t}-\zeta_{s})^2\leq  c_{f,\kappa}(T)|t-s|^{2-\kappa},\,\,  0\leq s,t\leq T.
\end{equation*}
In particular, the process $\zeta$ has Hölder continuous paths of index $\delta$, for any $0<\delta<1$.

(ii) Let $\alpha \in (-1,0)$. Assume that there exists a positive constant $c_{f,\alpha}$ such that $0\leq f\leq c_{f,\alpha}u^{\alpha}$ on $[0,T]$ for some $T>0$.
\\(ii.a) Then,  for any $\kappa\in(0,1)$ such that $1+\alpha-\kappa>0$ there exists a positive constant $c_{f,\kappa,\alpha}(T)$ such that:
\begin{equation*}\label{continuity2}
        \E(\zeta_{t}-\zeta_{s})^2\leq  c_{f,\kappa,\alpha}(T)(s\wedge t)^{\alpha}|t-s|^{2+\alpha-\kappa},\,\, 0< s, t\leq T.
\end{equation*}
In particular, the process $\zeta$ has locally Hölder continuous paths on $(0,\infty)$ of index $\delta$, for any $0<\delta<1+\alpha/2$.
\\(ii.b) Let $p\geq 1$. There exists a positive constant $c_{f,\alpha,p}(T)$ such that:
\begin{equation}\label{continuity2a}
        \E|\zeta_{t}-\zeta_{s}|^p\leq  c_{f,\alpha,p}(T)|t-s|^{\frac{p}{2}},\,\, 0\leq s,t\leq T.
\end{equation}
In particular, the process $\zeta$ has Hölder continuous paths of index $\delta$, for any $0<\delta<1/2$.
\end{thm}


\subsection*{Proof of Theorem~\ref{not-semimartingale}:} 
We begin by establishing several auxiliary results.

\begin{prop}\label{propA}
Suppose that \( f \) satisfies the condition
\begin{equation}\label{inf-var-condition}
   \lim_{n\rightarrow\infty} \frac{1}{n}\sum_{k=1}^{n-1} \sqrt{\int_0^{\frac{k}{n}} f(u)\frac{1}{\frac{k+1}{n}-u}\,du} = \infty.
\end{equation}
Then \( (\zeta_t)_{t \geq 0} \) has infinite variation almost surely on the interval \([0,1]\).
\end{prop}

\noindent\textbf{Proof.} A direct computation shows that, for \(0\le s<t\le 1\),
\begin{equation}\label{sec2}
        \E(\zeta_t-\zeta_s)^2 =4\log(2)\int_s^t f(u)(t-u)du+4\int_0^s f(u)\int_{s}^{t}\int_{s}^{t}\frac{1}{r+r^{'}-2u}drdr^{'}du.
   \end{equation}
Consider \(s = \frac{k}{n}\) and \(t = \frac{k+1}{n}\), for \(k \in \{0, 1, \ldots, n-1\}\). In this case,

\small
 \begin{eqnarray}\label{sec}
  \E\left(\zeta_{\frac{k+1}{n}}-\zeta_{\frac{k}{n}}\right)^2&=& 4\log(2)\int_{\frac{k}{n}}^{\frac{k+1}{n}} f(u)\left(\frac{k+1}{n}-u\right)du\\ \nonumber
       &+&4\int_0^{\frac{k}{n}} f(u)\int_{\frac{k}{n}}^{\frac{k+1}{n}}\int_{\frac{k}{n}}^{\frac{k+1}{n}}\frac{1}{r+r^{'}-2u}drdr^{'}du\\  \nonumber
       &\geq& \frac{2}{n^2} \int_0^\frac{k}{n} f(u)\frac{1}{\frac{k+1}{n}-u}du.
 \end{eqnarray}
Using \eqref{sec} and the absolute moment formula for Gaussian distributions, it follows that

\begin{equation}\label{infinite-variation-consequence}
        \E\Bigg(\Bigg|\zeta_{\frac{k+1}{n}}-\zeta_{\frac{k}{n}}\Bigg|\Bigg)
        \geq \sqrt{\frac{2}{\pi}}\frac{1}{n} \sqrt{\int_0^\frac{k}{n} f(u)\frac{1}{\frac{k+1}{n}-u}du}.
\end{equation}
Let \( |\zeta|_1 \) denote the total variation of \( \zeta \) on the interval \([0,1]\). From \eqref{infinite-variation-consequence} and \eqref{inf-var-condition}, we obtain \( \mathbb{E}(|\zeta|_1) = \infty \). According to the 0–1 law for Gaussian processes (\cite{Camb-Raj1973}, Theorem 2), the quantity \( |\zeta|_1 \) is either finite or infinite almost surely. If it were finite, Fernique's result (\cite{Fernique}) would imply that \( \mathbb{E}(|\zeta|_1) < \infty \), leading to a contradiction. This contradiction proves the result. \qed


As a consequence of Proposition~\ref{propA}, we derive the following result concerning the total variation of the process $\zeta$ for the function classes $\mathcal{C}_i$, $i = 1, 2$. More precisely,

\begin{cor}\label{corA}
Let \( f \in \mathcal{C}_1 \) with \( -1 < \alpha \leq 0 \), or \( f \in \mathcal{C}_2 \) with \( \alpha \in \mathbb{R} \). Then, the process \( \zeta \) has infinite variation almost surely on the interval \( [0,1] \).
\end{cor}

\smallskip
\noindent\textbf{Proof.} For the families of functions \( \mathcal{C}_1 \) and \( \mathcal{C}_2 \), note the following: if \( f(u) = u^\alpha \) with \( -1 < \alpha \leq 0 \), then \( f(u) \geq 1 \) for all \( u \in (0,1] \); likewise, if \( f(u) = e^{\alpha u} \) with \( \alpha \in \mathbb{R} \), then \( f(u) \geq \min\{1, e^{\alpha}\} \) for all \( u \in [0,1] \). Hence,
\[
\sqrt{\int_0^{\frac{k}{n}} f(u)\frac{1}{\frac{k+1}{n} - u} \, du} \geq 
C \sqrt{\int_0^{\frac{k}{n}} \frac{1}{\frac{k+1}{n} - u} \, du} 
= C \sqrt{\log(k+1)},
\]
for some constant \( C > 0 \). Therefore, condition~\eqref{inf-var-condition} is fulfilled, and the conclusion follows from Proposition~\ref{propA}.
\hfill\(\square\)


\begin{rem}\label{remA}
 (i) Consider the case \( f \equiv 1 \). In this case, we have
\[
\frac{1}{n}\sum_{k=1}^{n-1} \sqrt{\int_0^\frac{k}{n} f(u)\frac{1}{\frac{k+1}{n}-u}du}
= \frac{1}{n}\sum_{k=1}^{n-1}(\log(k+1))^{\frac{1}{2}} \xrightarrow[n \to \infty]{} \infty,
\]
which shows that condition~\eqref{inf-var-condition} is satisfied. Therefore, as a particular case of our criterion~\eqref{inf-var-condition}, we recover the result from \cite{BGT2007}.

(ii) If \( f \equiv c > 0 \), then by the same reasoning as in part (i), it follows that the process \( \zeta \) has infinite variation almost surely.
\end{rem}



The following proposition establishes a result concerning the quadratic variation of the process \( \zeta \).


\begin{prop}\label{propB}
Let \( f \) belong to the class \( \mathcal{C}_1 \) or \( \mathcal{C}_2 \). Then, the quadratic variation of the process \( \zeta \) on the interval \( [0,1] \) is equal to zero.
\end{prop}

\noindent\textbf{Proof.} The proof is structured in two parts:

\noindent\textbf{Part 1.} Consider the case \( f \in \mathcal{C}_1 \) with \( \alpha \geq 0 \) or \( f \in \mathcal{C}_2 \). Observe that \( u^{\alpha} \leq 1 \) for all \( u \in [0,1] \) when \( \alpha \geq 0 \), and \( e^{\alpha u} \leq \max\{1, e^{\alpha}\} \) for all \( u \in [0,1] \) when \( \alpha \in \mathbb{R} \). 

Substituting \( s = \frac{k}{n} \) and \( t = \frac{k+1}{n} \) into \eqref{sec2}, we obtain:
\small
\begin{eqnarray*}
\E\left(\zeta_{\frac{k+1}{n}} - \zeta_{\frac{k}{n}}\right)^2 
&=& 4\log(2) \int_{\frac{k}{n}}^{\frac{k+1}{n}} f(u)\left(\frac{k+1}{n} - u\right) du \\
&+& 4 \int_0^{\frac{k}{n}} f(u) \int_{\frac{k}{n}}^{\frac{k+1}{n}} \int_{\frac{k}{n}}^{\frac{k+1}{n}} \frac{1}{r + r' - 2u} \, dr \, dr' \, du \\
&\leq& C \left[ 4\log(2) \int_{\frac{k}{n}}^{\frac{k+1}{n}} \left(\frac{k+1}{n} - u\right) du + 4 \int_0^{\frac{k}{n}} \int_{\frac{k}{n}}^{\frac{k+1}{n}} \int_{\frac{k}{n}}^{\frac{k+1}{n}} \frac{1}{r + r' - 2u} \, dr \, dr' \, du \right] \\
&=& C \, \E\left(\eta_{\frac{k+1}{n}} - \eta_{\frac{k}{n}}\right)^2,
\end{eqnarray*}
\normalsize
for some constant \( C > 0 \). 

Using the convergence result
\begin{equation}\label{quadratic-variation-eta}
\lim_{n \rightarrow \infty} \sum_{k=0}^{n-1} \E\left[ \left( \eta_{\frac{k+1}{n}} - \eta_{\frac{k}{n}} \right)^2 \right] = 0,
\end{equation}
it follows that the quadratic variation of \( \zeta \) on \([0,1]\) is zero.\\

\textbf{Part 2.} Let us consider  $f\in\mathcal{C}_1$ with $-1<\alpha<0$. Note that
\begin{eqnarray}\label{sec2-1}
 \lefteqn{\int_0^{\frac{k}{n}} u^\alpha\int_{\frac{k}{n}}^{\frac{k+1}{n}}\int_{\frac{k}{n}}^{\frac{k+1}{n}}\frac{1}{r+r^{'}-2u}drdr^{'}du}\nonumber\\
 &=& \int_{\frac{k}{n}}^{\frac{k+1}{n}}\int_{\frac{k}{n}}^{\frac{k+1}{n}}\int_0^{\frac{k}{n}}\frac{u^{\alpha}}{r+r^{'}}\frac{1}{1-\frac{2u}{r+r^{'}}}dudrdr^{'}\nonumber\\
 &=& \sum_{i=0}^{\infty} \int_{\frac{k}{n}}^{\frac{k+1}{n}}\int_{\frac{k}{n}}^{\frac{k+1}{n}}\int_0^{\frac{k}{n}}\frac{u^{\alpha}}{r+r^{'}}\left(\frac{2u}{r+r^{'}}\right)^i dudrdr^{'}\\\nonumber
 &=& \left(\frac{k}{n}\right)^{\alpha}\sum_{i=0}^{\infty} \int_{\frac{k}{n}}^{\frac{k+1}{n}}\int_{\frac{k}{n}}^{\frac{k+1}{n}} 2^i \left(\frac{k}{n}\right)^{i+1} \frac{1}{(r+r^{'})^{i+1}} \cdot \frac{1}{1+i+\alpha} drdr^{'}.
\end{eqnarray}

Since $1+i+\alpha \geq (1+i)(1+\alpha) > 0$ for all $i\geq 0$, the left-hand side of \eqref{sec2-1} can be bounded as follows:
\begin{eqnarray}\label{sec2-2}
\lefteqn{\left(\frac{k}{n}\right)^{\alpha} \sum_{i=0}^{\infty} \int_{\frac{k}{n}}^{\frac{k+1}{n}} \int_{\frac{k}{n}}^{\frac{k+1}{n}} 2^i \left(\frac{k}{n}\right)^{i+1} \frac{1}{(r+r^{'})^{i+1}} \cdot \frac{1}{1+i+\alpha} drdr^{'}} \nonumber\\
&\leq& \frac{1}{1+\alpha} \left(\frac{k}{n}\right)^{\alpha} \sum_{i=0}^{\infty} \int_{\frac{k}{n}}^{\frac{k+1}{n}} \int_{\frac{k}{n}}^{\frac{k+1}{n}} 2^i \left(\frac{k}{n}\right)^{i+1} \frac{1}{(r+r^{'})^{i+1}} \cdot \frac{1}{1+i} drdr^{'} \nonumber\\
&=& \frac{1}{1+\alpha} \left(\frac{k}{n}\right)^{\alpha} \int_0^{\frac{k}{n}} \int_{\frac{k}{n}}^{\frac{k+1}{n}} \int_{\frac{k}{n}}^{\frac{k+1}{n}} \frac{1}{r+r^{'} - 2u} \, dr \, dr^{'} \, du.
\end{eqnarray}

On the other hand, observe that
\begin{equation}\label{sec2-3}
\int_{\frac{k}{n}}^{\frac{k+1}{n}} f(u)\left(\frac{k+1}{n}-u\right) du \leq \frac{1}{1+\alpha} \left(\frac{k}{n}\right)^{\alpha} \int_{\frac{k}{n}}^{\frac{k+1}{n}} \left(\frac{k+1}{n}-u\right) du.
\end{equation}

From \eqref{sec2-1}--\eqref{sec2-3}, we deduce that
\begin{equation}\label{key-quad-var-negative}
\E\left(\zeta_{\frac{k+1}{n}} - \zeta_{\frac{k}{n}}\right)^2 \leq C \left(\frac{k}{n}\right)^{\alpha} \E\left(\eta_{\frac{k+1}{n}} - \eta_{\frac{k}{n}}\right)^2 \leq C \left(\frac{k}{n}\right)^{\alpha} \cdot \frac{\log(k)}{n^2}, \quad k \geq e^4,
\end{equation}
for some constant \( C > 0 \), where the last inequality follows from
\begin{equation}\label{gorostiza-key}
\frac{\log(k)}{2n^2} \leq \E\left(\eta_{\frac{k+1}{n}} - \eta_{\frac{k}{n}}\right)^2 \leq \frac{\log(k)}{n^2}, \quad \text{for } k \geq e^4,
\end{equation}
as established in \cite{BGT2007}.

Finally, since \( -1 < \alpha < 0 \), \eqref{key-quad-var-negative} implies
\[
\lim_{n\rightarrow\infty} \sum_{k=0}^{n-1} \E\left(\zeta_{\frac{k+1}{n}} - \zeta_{\frac{k}{n}}\right)^2 = 0.
\]
Thus, the quadratic variation of the process \( \zeta \) on \( [0,1] \) is equal to zero almost surely.
\hfill\(\square\)


The following result completes the analysis regarding the total variation of the process~$\zeta$. More precisely, we establish the following:

\begin{prop}\label{propC}
Let \( f \in \mathcal{C}_1 \) with \( \alpha > 0 \). Then, the process \( \zeta \) has infinite variation almost surely on the interval \([0,1]\).
\end{prop}


\noindent\textbf{Proof.} We only provide a sketch of the proof, as it follows arguments similar to those employed in Part~2 of the proof of Proposition~\ref{propB}. Observe that the condition $\alpha > 0$ implies that $1 + i + \alpha \leq (1 + i)(1 + \alpha)$ for all $i \geq 0$. Hence, using the lower bound in~\eqref{gorostiza-key}, we obtain
\begin{equation}\label{infinite-variation-positive}
   \E\left(\zeta_{\frac{k+1}{n}}-\zeta_{\frac{k}{n}}\right)^2 \geq C  \left(\frac{k}{n}\right)^{\alpha}\frac{\log(k)}{2n^2},
\end{equation}
for some constant $C > 0$. Therefore, by Gaussianity, inequality~\eqref{infinite-variation-positive} implies that
\begin{equation*}
   \E\left|\zeta_{\frac{k+1}{n}}-\zeta_{\frac{k}{n}}\right| \geq C \frac{\sqrt{k^\alpha\log(k)}}{n^{1+\frac{\alpha}{2}}}.
\end{equation*}
It is easy to show that
\[
\frac{1}{n^{1 + \frac{\alpha}{2}}} \sum_{k = 1}^{n} \sqrt{k^{\alpha} \log(k)} \xrightarrow[n \to \infty]{} \infty.
\]
The result follows by the same arguments used to complete the proof of Proposition~\ref{propA}.
\hfill$\square$

Finally, we arrive at the desired result.  

\noindent\textbf{Proof of Theorem~\ref{not-semimartingale}.} The result follows directly from Corollary~\ref{corA}, Proposition~\ref{propB}, and Proposition~\ref{propC}.
\hfill$\square$

\subsection*{Proof of Theorem~\ref{teo3}:} The proof relies on the auxiliary result stated in Lemma~\ref{lemA}.
Throughout this argument, we adopt the following notation: given two positive real-valued functions \( f \) and \( g \), we write \( f \asymp g \) if there exist constants \( c_1, c_2 > 0 \) such that \( c_1 f \leq g \leq c_2 f \). Observe that if \( f \asymp g \) holds with constants \( c_1 \) and \( c_2 \), then \( g \asymp f \) also holds, with constants \( c_2^{-1} \) and \( c_1^{-1} \).


\begin{lem}\label{lemA}
(i) Let us assume that \( 1 \asymp f \) on \( [0,T] \), for some \( T>0 \), with constants \( c_{1,f} \) and \( c_{2,f} \). Then,
{\small
\begin{eqnarray}\label{zetad}
   c_{1,f}\frac{(t-s)^2}{4}\left[2\log(2s)+3-2\log(t-s)\right] &\leq& \E(\zeta_{t}-\zeta_{s})^2 \\
   &\leq& c_{2,f}\frac{(t-s)^2}{4}\left[2\log(2t)+3-2\log(t-s)\right].\nonumber
\end{eqnarray}}

(ii) Let us assume that \( f \asymp u^{\alpha} \) on \( [0,T] \), for some \( T>0 \) and \( \alpha \in (-1,\infty) \), with constants \( c_{1,f,\alpha} \) and \( c_{2,f,\alpha} \).

(ii.a) If \( \alpha \in (-1,0) \), then for all \( s,t \in [0,T] \),
{\small
\begin{eqnarray}\label{zetad2}
&& c_{1,f,\alpha}s^{\alpha}\frac{(t-s)^2}{4}\left[2\log(2s)+3-2\log(t-s)-8\log(2)\right] \leq \E(\zeta_{t}-\zeta_{s})^2 \\\nonumber
&& \leq \frac{c_{2,f,\alpha}}{1+\alpha}s^{\alpha}\frac{(t-s)^2}{4}\left[2\log(2t)+3-2\log(t-s)\right] + c_{2,f,\alpha}\frac{4\log(2)c_{\alpha}}{1+\alpha}(t-s)^{2+\alpha},
\end{eqnarray}}
where \( c_{\alpha} \) is a positive constant.

(ii.b) If \( \alpha \in [0,\infty) \), then for all \( s,t \in [0,T] \),
{\small
\begin{eqnarray}\label{zetad3}
&& \frac{c_{1,f,\alpha}}{1+\alpha}s^{\alpha}\frac{(t-s)^2}{4}\left[2\log(2s)+3-2\log(t-s)-8\log(2)\right] \leq \E(\zeta_{t}-\zeta_{s})^2 \\\nonumber
&& \leq c_{2,f,\alpha}s^{\alpha}\frac{(t-s)^2}{4}\left[2\log(2t)+3-2\log(t-s)\right] + 4c_{2,f,\alpha}\log(2)(t-s)^2.
\end{eqnarray}}
\end{lem}

\noindent\textbf{Proof.} (i) It is enough to consider \(0 \leq s < t \leq T\). Since \(1 \asymp f\), it follows from equation~\eqref{sec2} that there exist positive constants \(c_{1,f}\) and \(c_{2,f}\) such that
\[
c_{1,f} \, \E\left(\eta_{t}-\eta_{s}\right)^2 \leq \E\left(\zeta_{t}-\zeta_{s}\right)^2 \leq c_{2,f} \, \E\left(\eta_{t}-\eta_{s}\right)^2.
\]
From equation (3.7) in \cite{BGT2007}, and observing that \( \log(2s) \leq \log(2s + u + u') \leq \log(2t) \) for all \( 0 \leq u, u' \leq t - s \), the following estimates are deduced:
\begin{equation}\label{etad}
\frac{(t-s)^2}{4}\left[2\log(2s) + 3 - 2\log(t-s)\right] \leq \E(\eta_t - \eta_s)^2 \leq \frac{(t-s)^2}{4}\left[2\log(2t) + 3 - 2\log(t-s)\right].
\end{equation}
Therefore, inequality~\eqref{zetad} holds, which completes the proof of part~(i). Note that \( 2\log(2s) + 3 - 2\log(t - s) \geq 0 \) if and only if \( 2s e^{3/2} \geq t - s \).

(ii) The proof of part~(ii) follows by combining the arguments used in part~(i) with equations~\eqref{sec2}, \eqref{sec2-1}, and \eqref{etad}.
\hfill\(\square\)
\begin{prop}\label{prop:nonstat_intrinsic}
Assume that \(f:[0,\infty)\to[0,\infty)\) satisfies \eqref{Integrability-condition} and that there exists \(t_0>0\) such that
\[
\int_0^{t_0} f(u)\,du>0.
\]
Then the process \((\zeta_t)_{t\ge0}\) is non-stationary.

If, in addition, \(1\asymp f\) on \([0,T]\) for some \(T>0\), then \((\zeta_t)_{t\ge0}\) is not intrinsically stationary on \([0,T]\).
\end{prop}

\begin{proof}
We first prove that \((\zeta_t)_{t\ge0}\) is non-stationary. Since the process is centered, weak stationarity would imply that the variance is constant in time. By taking \(s=0\) in \eqref{sec2}, we obtain
\[
\E(\zeta_t^2)=4\log(2)\int_0^t f(u)(t-u)\,du,\qquad t\ge0.
\]
Let \(h>0\). Then
\[
\E(\zeta_{t_0+h}^2)-\E(\zeta_{t_0}^2)
=
4\log(2)\left[
h\int_0^{t_0} f(u)\,du+\int_{t_0}^{t_0+h} f(u)(t_0+h-u)\,du
\right].
\]
Since \(f\ge0\) and \(\int_0^{t_0} f(u)\,du>0\), the right-hand side is strictly positive. Hence
\[
\E(\zeta_{t_0+h}^2)>\E(\zeta_{t_0}^2),
\]
so the variance is not constant on \((0,\infty)\). Therefore, \((\zeta_t)_{t\ge0}\) is non-stationary.

We now prove that the process is not intrinsically stationary under the additional assumption \(1\asymp f\) on \([0,T]\). Then there exist constants \(c_1,c_2>0\) such that
\[
c_1\le f(u)\le c_2,\qquad u\in[0,T].
\]
For \(0<h<T\), taking \(s=0\) in \eqref{sec2}, we get
\[
\E(\zeta_h-\zeta_0)^2=\E(\zeta_h^2)=4\log(2)\int_0^h f(u)(h-u)\,du,
\]
and therefore
\[
2c_1\log(2)\,h^2
\le
\E(\zeta_h-\zeta_0)^2
\le
2c_2\log(2)\,h^2.
\]
Thus,
\[
\frac{\E(\zeta_h-\zeta_0)^2}{h^2}
\]
remains bounded as \(h\downarrow0\).

On the other hand, fix \(s\in(0,T)\). By \eqref{zetad}, with \(t=s+h\), we have for all sufficiently small \(h>0\),
\[
\E(\zeta_{s+h}-\zeta_s)^2
\ge
c_{1,f}\frac{h^2}{4}\bigl[2\log(2s)+3-2\log h\bigr].
\]
Hence
\[
\frac{\E(\zeta_{s+h}-\zeta_s)^2}{h^2}
\ge
c_{1,f}\frac{1}{4}\bigl[2\log(2s)+3-2\log h\bigr]
\longrightarrow \infty,
\qquad h\downarrow0.
\]
Therefore, for sufficiently small \(h>0\),
\[
\E(\zeta_h-\zeta_0)^2 \neq \E(\zeta_{s+h}-\zeta_s)^2.
\]
Hence the variance of the increment depends on the initial time, and the process is not intrinsically stationary on \([0,T]\). \qed
\end{proof}

\noindent{\bf Proof of Theorem~\ref{teo3}.} (i) Let \(\kappa \in (0,1)\). Proceeding analogously to the proof of Lemma~\ref{lemA}~(i), and using that the function \(x \mapsto -x^\kappa \log(x)\) is nonnegative and bounded on \([0,1]\), we obtain:
\begin{equation}\label{continuity}
\E(\zeta_{t}-\zeta_{s})^2 \leq c_f \frac{(t-s)^2}{2} \left[2\log(2(s \vee t)) + 3 - 2\log|t - s| \right] \leq c_{f,\kappa}(T) |t - s|^{2 - \kappa}, \quad \forall\, 0 \leq s, t \leq T,
\end{equation}
for some constant \( c_{f,\kappa}(T) > 0 \). Therefore, by the Kolmogorov continuity theorem, the sample paths of \(\zeta\) are Hölder continuous.

(ii.a) The result follows from arguments analogous to those in part~(i), combined with the estimates in Lemma~\ref{lemA}~(ii.a) and the Kolmogorov continuity theorem.

(ii.b) From equation~\eqref{sec2}, it follows that there exists a constant \( M_{f,\alpha}(T) > 0 \) such that
\[
\E(\zeta_t - \zeta_s)^2 \leq M_{f,\alpha}(T) |t - s|, \quad \text{for all } 0 \leq s, t \leq T.
\]
Since \(\zeta\) is a centered Gaussian process, the moment estimates for Gaussian increments and the Kolmogorov continuity theorem imply the Hölder continuity of the sample paths.
\hfill\(\square\)

\subsection{Memory properties}

In this subsection, we establish the asymptotic covariance behavior underlying the memory properties of the process $\zeta$. We first provide the proof of Proposition~1 from the main text, which yields the logarithmic asymptotic growth of the covariance and the corresponding long-range dependence property. We then present two additional complementary results concerning short-time asymptotics and the continuity behavior of the process at the origin.

We begin with the proof of Proposition~1 stated in the main text.

\subsection*{Proof of Proposition 1}

(i) For \( r, s > 0 \), we observe that
\begin{eqnarray*}
\frac{K_f(r, T+s)}{\log(T)} &=& 2\int_0^r \frac{f(u)}{\log(T)} \Big[(T+s+r - 2u)\log(T+s+r - 2u) \\
&& \quad - (r - u)\log(r - u) - (T + s - u)\log(T + s - u)\Big]\,du \\
&=& 2\int_0^r f(u) \left[\int_0^{r - u} \frac{\log(T + s - u + x) + 1}{\log(T)}\,dx - \frac{(r - u)\log(r - u)}{\log(T)}\right]\,du \\
&\xrightarrow[T \to \infty]{}& 2\int_0^r f(u)(r - u)\,du.
\end{eqnarray*}

(ii) Let \( 0 < r \leq \nu \) and \( 0 < s \leq t \). Then, it holds that
\[
K_f(r, T+t) - K_f(r, T+s) = 2\int_0^r f(u)\int_0^{r - u} \left[\log(T + t - u + x) - \log(T + s - u + x)\right]\,dx\,du.
\]
Applying L’Hôpital’s rule, the integrability condition~\eqref{Integrability-condition}, and the Dominated Convergence Theorem, we get
\begin{eqnarray*}
\frac{K_f(r, T+t) - K_f(r, T+s)}{T^{-1}} &=& 2\int_0^r f(u)\int_0^{r - u} \frac{\log(T + t - u + x) - \log(T + s - u + x)}{T^{-1}}\,dx\,du \\
&\xrightarrow[T \to \infty]{}& 2(t - s)\int_0^r f(u)(r - u)\,du.
\end{eqnarray*}
Therefore,
\begin{eqnarray*}
&&\frac{K_f(\nu, T + t) - K_f(\nu, T + s) - K_f(r, T + t) + K_f(r, T + s)}{T^{-1}} \\
&&\quad \xrightarrow[T \to \infty]{} 2(t - s) \left[\int_0^\nu f(u)(\nu - u)\,du - \int_0^r f(u)(r - u)\,du\right].
\end{eqnarray*}

This completes the proof. \hfill\(\square\)

The previous result describes the large-lag covariance behavior of the process. We next present two additional results that complement the local analysis of the sample paths.

\begin{cor}\label{corB}(Short-time asymptotics)
 Let \( t > 0 \). Then, the following assertions hold:

(i) Suppose that \( 1 \asymp f \), that is, there exist positive constants \( c_{1,f} \) and \( c_{2,f} \) such that \( c_{1,f} \leq f \leq c_{2,f} \). Then,
\begin{equation}\label{zetas}
\frac{c_{1,f}}{2} \leq \liminf_{\epsilon \downarrow 0} \frac{\E(\zeta_{t+\epsilon} - \zeta_{t})^2}{-\epsilon^2 \log(\epsilon)} \leq \limsup_{\epsilon \downarrow 0} \frac{\E(\zeta_{t+\epsilon} - \zeta_{t})^2}{-\epsilon^2 \log(\epsilon)} \leq \frac{c_{2,f}}{2}.
\end{equation}
In particular,
\begin{equation}\label{etas} 
\lim_{\epsilon \downarrow 0} \frac{\E(\eta_{t+\epsilon} - \eta_{t})^2}{-\epsilon^2 \log(\epsilon)} = \frac{1}{2}.
\end{equation}

(ii) Suppose that \( u^{\alpha} \asymp f \) on \([0,T]\), with \( \alpha \in (-1, \infty) \). That is, there exist positive constants \( c_{1,f,\alpha} \) and \( c_{2,f,\alpha} \) such that \( c_{1,f,\alpha} u^{\alpha} \leq f(u) \leq c_{2,f,\alpha} u^{\alpha} \).

(ii.a) If \( \alpha \in (-1,0) \), then
\begin{equation*}
\frac{c_{1,f,\alpha} t^{\alpha}}{2} \leq \liminf_{\epsilon \downarrow 0} \frac{\E(\zeta_{t+\epsilon} - \zeta_{t})^2}{-\epsilon^2 \log(\epsilon)} \leq \limsup_{\epsilon \downarrow 0} \frac{\E(\zeta_{t+\epsilon} - \zeta_{t})^2}{-\epsilon^2 \log(\epsilon)} \leq \frac{c_{2,f,\alpha} t^{\alpha}}{2(1 + \alpha)}.
\end{equation*}

(ii.b) If \( \alpha \in [0, \infty) \), then
\begin{equation*}
\frac{c_{1,f,\alpha} t^{\alpha}}{2(1 + \alpha)} \leq \liminf_{\epsilon \downarrow 0} \frac{\E(\zeta_{t+\epsilon} - \zeta_{t})^2}{-\epsilon^2 \log(\epsilon)} \leq \limsup_{\epsilon \downarrow 0} \frac{\E(\zeta_{t+\epsilon} - \zeta_{t})^2}{-\epsilon^2 \log(\epsilon)} \leq \frac{c_{2,f,\alpha} t^{\alpha}}{2}.
\end{equation*}
\end{cor}

\noindent\textbf{Proof.} (i) The identity in~\eqref{etas} follows directly from~\eqref{etad}, while the bounds in~\eqref{zetas} are a consequence of~\eqref{zetad}.

(ii) For the general case \( u^{\alpha} \asymp f \), observe that, by L’Hôpital’s rule, one can verify that
\[
\lim_{\epsilon \downarrow 0} \frac{4\log(2)}{\epsilon^2} \int_t^{t+\epsilon} u^{\alpha} (t+\epsilon - u) \, du = 2\log(2)\, t^{\alpha}.
\]
Hence,
\begin{equation}\label{lim_a}
\lim_{\epsilon \downarrow 0} \frac{4\log(2)}{-\epsilon^2 \log(\epsilon)} \int_t^{t+\epsilon} u^{\alpha} (t+\epsilon - u) \, du = 0.
\end{equation}

The desired result then follows by applying Lemma~\ref{lemA}~(ii), \eqref{sec2}, and the limit in~\eqref{lim_a}. \hfill\(\square\)

The following proposition concerns the path continuity of $\zeta$ at the origin.

\begin{prop}\label{propD} 
(i) Let \( f(u) = u^{\alpha} \in \mathcal{C}_1 \). Then,
\[
\lim_{\epsilon \downarrow 0} \frac{\E(\zeta_{\epsilon})^2}{\epsilon^{2+\alpha}} = \frac{2}{(\alpha+1)(\alpha+2)}.
\]
(ii) Let \( f(u) = e^{\alpha u} \in \mathcal{C}_2 \). Then,
\[
\lim_{\epsilon \downarrow 0} \frac{\E(\zeta_{\epsilon})^2}{\epsilon^{2}} = 1.
\]
\end{prop}
\noindent\textbf{Proof.} The conclusion follows by direct application of L’Hôpital’s rule to the integral representation of \(\E(\zeta_{\epsilon})^2\).
\hfill$\square$





\section{Empirical diagnostics for the telemetry data}

\begingroup
\setlength{\textfloatsep}{10pt plus 2pt minus 2pt}
\setlength{\floatsep}{8pt plus 2pt minus 2pt}
\setlength{\intextsep}{8pt plus 2pt minus 2pt}
\setlength{\abovecaptionskip}{4pt}
\setlength{\belowcaptionskip}{0pt}

\subsection{Empirical evidence on weak stationarity and stationarity of first differences}\label{appC:NS}

We examined the stationarity properties of the ten telemetry series obtained
from five bats, considering both longitude and latitude. Figure~\ref{fig:rolling_mean}
displays the rolling mean of the raw series, which varies over time, while
Figure~\ref{fig:rolling_variance} reports the corresponding rolling variance,
showing non-constant dispersion. The sample autocorrelation functions of the
raw series, shown in Figure~\ref{fig:acf_raw}, frequently exceed the 95\%
confidence bounds, indicating substantial temporal dependence. After first
differencing, the autocorrelations are reduced, but significant serial
dependence remains in several cases; see Figure~\ref{fig:acf_diff}. Thus,
first differencing does not consistently remove the serial dependence present
in the data. This residual autocorrelation does not, by itself, contradict
intrinsic stationarity, since intrinsically stationary increments may still be
serially correlated.

\begin{figure}[H]
    \centering
    \includegraphics[width=0.80\textwidth]{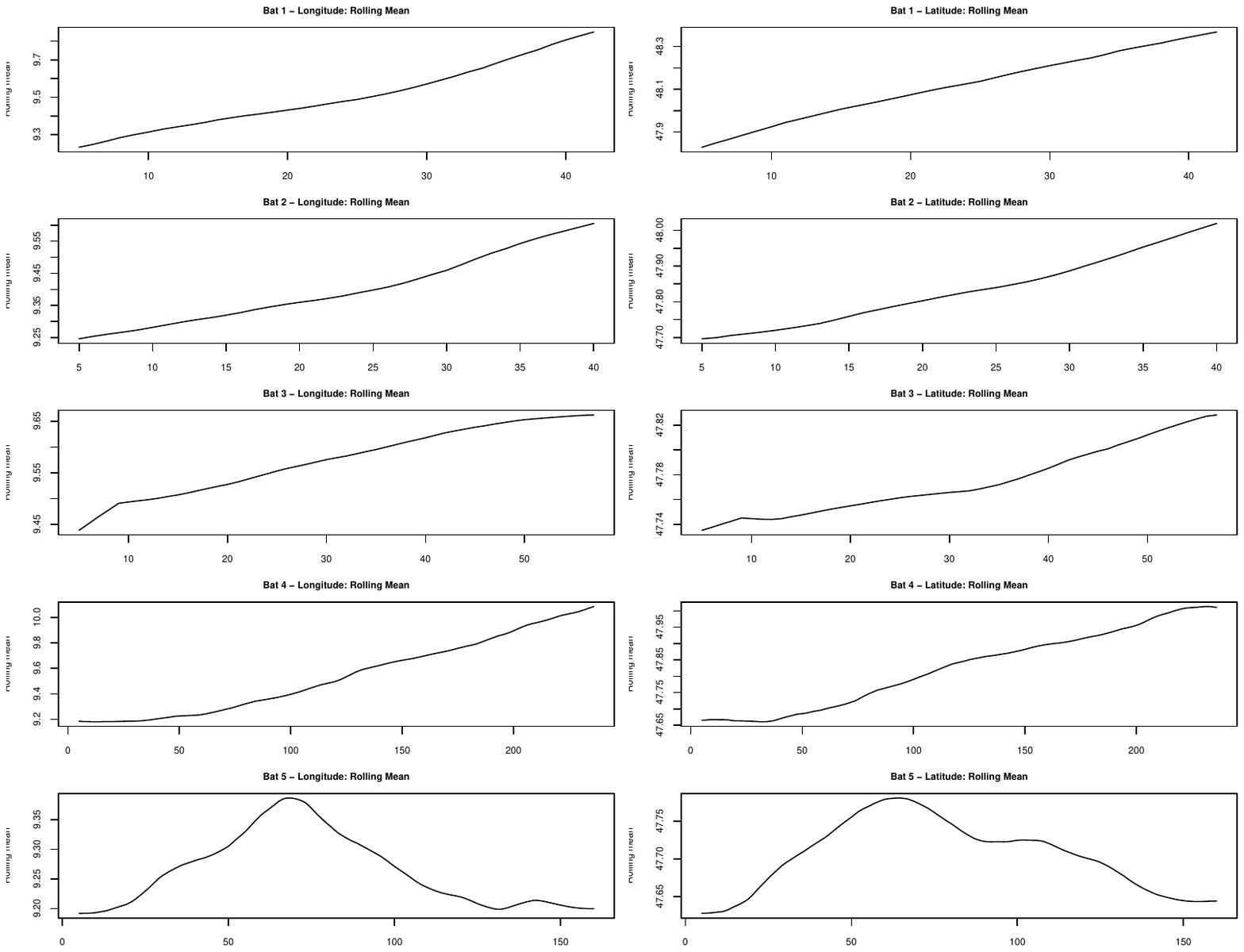}
    \caption{Rolling mean of the latitude and longitude trajectories for all bats, computed with a moving window of width $10$. The trajectories do not fluctuate around a constant level, which is consistent with non-stationarity in the mean.}
    \label{fig:rolling_mean}
\end{figure}

\begin{figure}[H]
    \centering
    \includegraphics[width=0.80\textwidth]{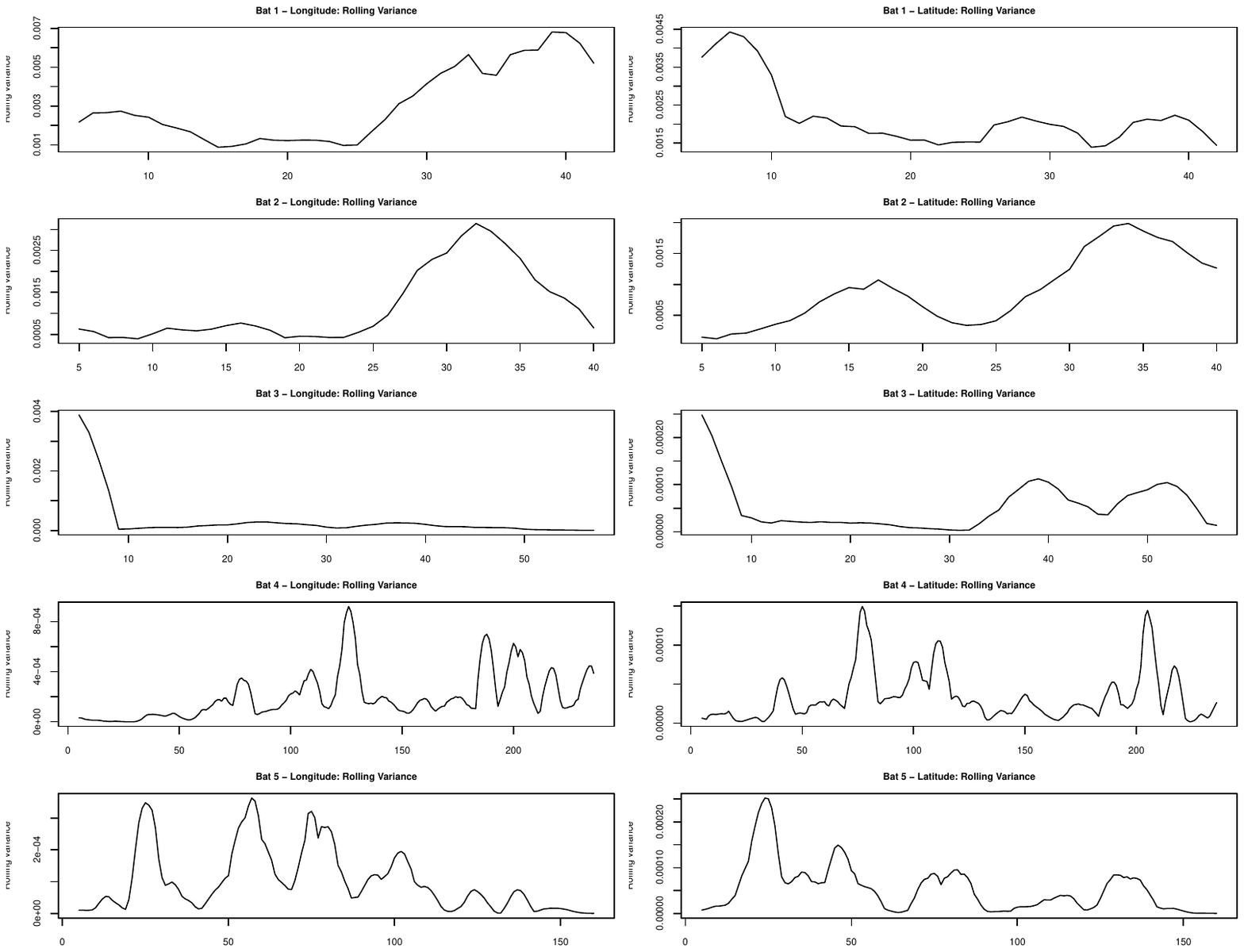}
    \caption{Rolling variance of the latitude and longitude trajectories for all bats, computed with a moving window of width $10$. The clear temporal variation in dispersion is consistent with non-stationarity in the variance.}
    \label{fig:rolling_variance}
\end{figure}

\begin{figure}[H]
    \centering
    \includegraphics[width=0.9\textwidth]{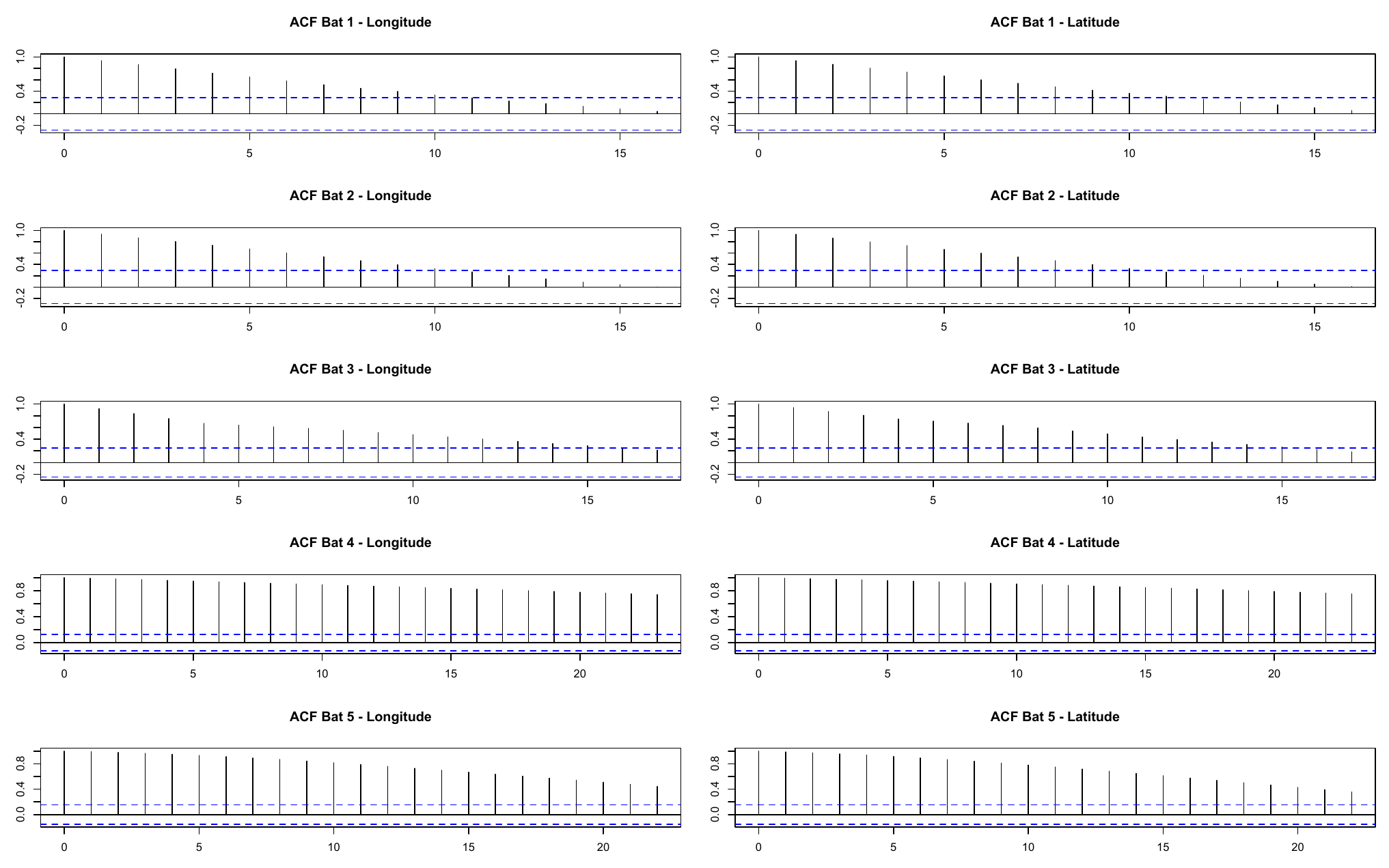}
    \caption{Sample autocorrelation functions of the raw longitude and latitude series for all bats. Black lines correspond to the empirical autocorrelations, and blue lines indicate the 95\% confidence bounds. In several cases, the autocorrelations remain significant over multiple lags.}
    \label{fig:acf_raw}
\end{figure}

\begin{figure}[H]
    \centering
    \includegraphics[width=0.85\textwidth,height=0.7\textheight,keepaspectratio]{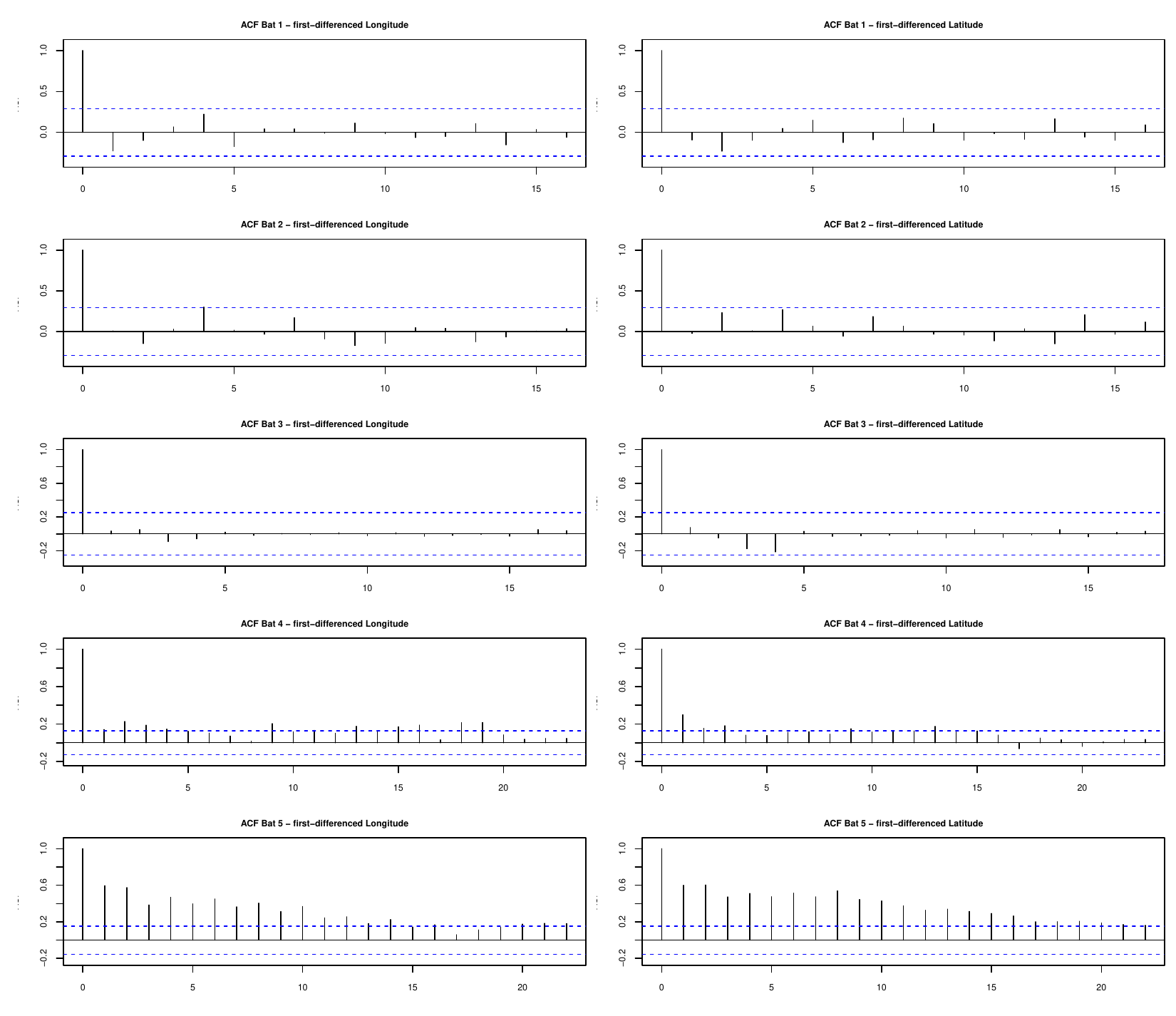}
    \caption{Sample autocorrelation functions of the first-differenced longitude and latitude series for all bats. Black lines correspond to the empirical autocorrelations, and blue lines indicate the 95\% confidence bounds. Although differencing reduces dependence, significant autocorrelations remain in several datasets.}
    \label{fig:acf_diff}
\end{figure}
\FloatBarrier

Let \((X_i)_{i=0}^{n}\) denote any one of the observed longitude or latitude series. Define the first differences by \(\Delta X_i:=X_i-X_{i-1}\), \(i=1,\ldots,n\), and the second differences by \(\Delta^2X_i:=\Delta X_i-\Delta X_{i-1}=X_i-2X_{i-1}+X_{i-2}\), \(i=2,\ldots,n\). To investigate whether an additional difference removes the serial correlation remaining after first differencing, we computed the sample autocorrelation functions of the second-differenced longitude and latitude series for all five bats.

\begin{figure}[H]
\centering
\includegraphics[width=0.65\textwidth]{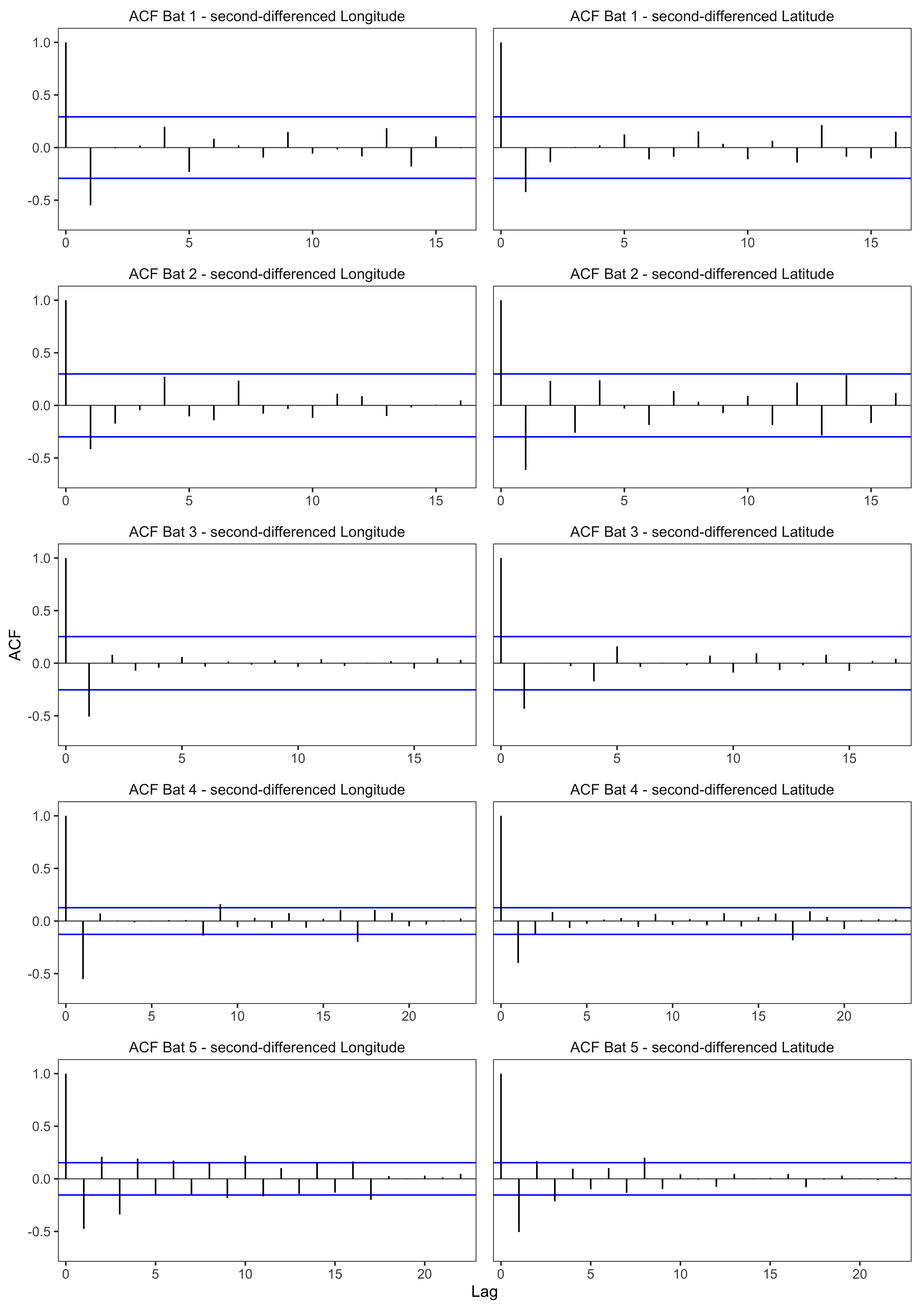}
\caption{Sample autocorrelation functions of the second-differenced longitude and latitude series for the five bats. The blue horizontal lines represent the approximate 95\% confidence bounds under a white-noise reference. In all ten series, the lag-one autocorrelation falls below the lower confidence bound, while several series retain autocorrelations outside the bounds at additional lags.}
\label{fig:acf-second-differences}
\end{figure}

As shown in Figure~\ref{fig:acf-second-differences}, second-order differencing does not uniformly whiten the telemetry series. Although it reduces some of the positive serial correlation remaining after first differencing, all ten second-differenced series exhibit a pronounced negative autocorrelation at lag one. Several series, particularly those corresponding to Bats~4 and~5, also retain autocorrelations outside the approximate confidence bounds at additional lags.

The systematic negative lag-one correlation is consistent with over-differencing. Indeed, suppose that the first differences form a zero-mean finite-variance white-noise sequence, so that \(\Delta X_i=\varepsilon_i\), where \(\mathbb{E}[\varepsilon_i]=0\), \(\operatorname{Var}(\varepsilon_i)=\tau^2\in(0,\infty)\), and \(\operatorname{Cov}(\varepsilon_i,\varepsilon_j)=0\) for \(i\neq j\). Then \(\Delta^2X_i=\varepsilon_i-\varepsilon_{i-1}\), which is an MA(1) process. Moreover,
\[
\operatorname{Var}(\Delta^2X_i)=2\tau^2,
\qquad
\operatorname{Cov}(\Delta^2X_i,\Delta^2X_{i-1})=-\tau^2,
\]
and therefore
\[
\operatorname{Corr}(\Delta^2X_i,\Delta^2X_{i-1})
=
\frac{-\tau^2}{2\tau^2}
=
-\frac12.
\]
Thus, a negative lag-one autocorrelation near \(-1/2\) may be induced by applying an unnecessary additional difference. These results do not support treating all trajectories as uniformly integrated of order two. More generally, autocorrelations remaining after differencing indicate serial dependence, but they do not by themselves establish nonstationarity, since a stationary process may also be serially correlated.

To complement the graphical analysis, we applied the Augmented
Dickey--Fuller (ADF) and Kwiatkowski--Phillips--Schmidt--Shin (KPSS)
tests. The results are summarized in Table~\ref{tab:adf_kpss}. For the
raw series, the KPSS test rejects stationarity in all ten cases, whereas
the ADF test rejects the unit-root null in only two cases. Hence, the raw
trajectories provide strong overall evidence against weak stationarity.
For the first-differenced series, the conclusions are mixed: both tests are
simultaneously consistent with stationarity in four of the ten cases,
whereas in the remaining six series at least one of the two tests indicates
non-stationarity. Thus, the evidence concerning intrinsic stationarity is
not uniform across the ten trajectories. The joint ADF--KPSS conclusions in
four cases, together with the absence of a KPSS rejection in two additional
cases, make a stationary-increment description plausible for a substantial
subset of the data. The remaining outcomes support retaining alternatives
without stationary increments. These diagnostics therefore motivate a model
comparison containing both types of non-stationary position processes rather
than imposing a common increment structure on all ten series.

\begin{table}[!h]
\centering
\normalsize
\setlength{\tabcolsep}{3pt}
\renewcommand{\arraystretch}{0.85}
\begin{tabular}{@{}cccccc@{}}
\toprule
\textbf{Bat} & \textbf{Series} & \textbf{ADF $p$} & \textbf{KPSS $p$} & \textbf{ADF} & \textbf{KPSS} \\
\midrule
\multicolumn{6}{l}{\textit{Panel A: Raw series}} \\
\midrule
1 & Longitude & 0.9724 & 0.0100 & Non-stat. & Non-stat. \\
1 & Latitude  & 0.4993 & 0.0100 & Non-stat. & Non-stat. \\
2 & Longitude & 0.8021 & 0.0100 & Non-stat. & Non-stat. \\
2 & Latitude  & 0.8010 & 0.0100 & Non-stat. & Non-stat. \\
3 & Longitude & 0.0100 & 0.0100 & Stat.     & Non-stat. \\
3 & Latitude  & 0.3952 & 0.0100 & Non-stat. & Non-stat. \\
4 & Longitude & 0.0415 & 0.0100 & Stat.     & Non-stat. \\
4 & Latitude  & 0.4264 & 0.0100 & Non-stat. & Non-stat. \\
5 & Longitude & 0.5748 & 0.0100 & Non-stat. & Non-stat. \\
5 & Latitude  & 0.4934 & 0.0100 & Non-stat. & Non-stat. \\
\midrule
\multicolumn{6}{l}{\textit{Panel B: First differences}} \\
\midrule
1 & $\Delta$ Longitude & 0.3341 & 0.1000 & Non-stat. & Stat.     \\
1 & $\Delta$ Latitude  & 0.0100 & 0.1000 & Stat.     & Stat.     \\
2 & $\Delta$ Longitude & 0.3980 & 0.1000 & Non-stat. & Stat.     \\
2 & $\Delta$ Latitude  & 0.1770 & 0.0178 & Non-stat. & Non-stat. \\
3 & $\Delta$ Longitude & 0.0100 & 0.0991 & Stat.     & Stat.     \\
3 & $\Delta$ Latitude  & 0.0100 & 0.1000 & Stat.     & Stat.     \\
4 & $\Delta$ Longitude & 0.0100 & 0.0100 & Stat.     & Non-stat. \\
4 & $\Delta$ Latitude  & 0.0262 & 0.1000 & Stat.     & Stat.     \\
5 & $\Delta$ Longitude & 0.3486 & 0.0100 & Non-stat. & Non-stat. \\
5 & $\Delta$ Latitude  & 0.2016 & 0.0100 & Non-stat. & Non-stat. \\
\bottomrule
\end{tabular}
\caption{ADF and KPSS results for the raw longitude and latitude series and for their first differences. For the ADF test, the null hypothesis is non-stationarity, whereas for the KPSS test, the null hypothesis is stationarity. The entries ``Stat.'' and ``Non-stat.'' indicate the corresponding decision at the 5\% significance level.}
\label{tab:adf_kpss}
\end{table}
\endgroup

\newpage
\subsection{Empirical DFA evidence of large-scale dependence}
\label{appC:DFA}

As a complementary empirical diagnostic, we applied detrended fluctuation
analysis (DFA) to each telemetry series. The DFA methodology was introduced
by \cite{Peng1994} and later systematically studied by
\cite{Kantelhardt2001} as a tool for analyzing scaling behavior and
correlation structure in time series, especially in the presence of trends
or other nonstationary features. Briefly, DFA is based on the cumulative
profile of the series, the estimation and removal of local trends over
windows of different sizes, and the study of the resulting fluctuation
function \(F(s)\) across scales \(s\).

In a finite sample, when the log--log plot of \(F(s)\) is approximately
linear over a selected range of scales, the empirical relation may be
summarized as
\[
F(s)\approx C s^\alpha ,
\]
where the fitted slope \(\alpha\) describes the scaling behavior over that
range. As a classical point of reference in DFA-based interpretations,
values of \(\alpha\) near \(0.5\) are commonly associated with weak or
negligible serial correlation, whereas larger values indicate stronger
dependence across scales. In the present application, these benchmarks are
used descriptively, since the telemetry series are non-stationary and the
stationarity of their first differences is not uniformly supported by the
ADF and KPSS tests. Moreover, as emphasized by \cite{Kantelhardt2001},
trends and nonstationary effects may affect the observed scaling behavior
and should therefore be taken into account when interpreting the results.

For this reason, in the present appendix DFA is used only as an additional
empirical diagnostic and is interpreted jointly with the rolling summaries,
ACF, ADF, and KPSS analyses reported above.

\begin{figure}[H]
\centering
\includegraphics[width=0.7\textwidth]{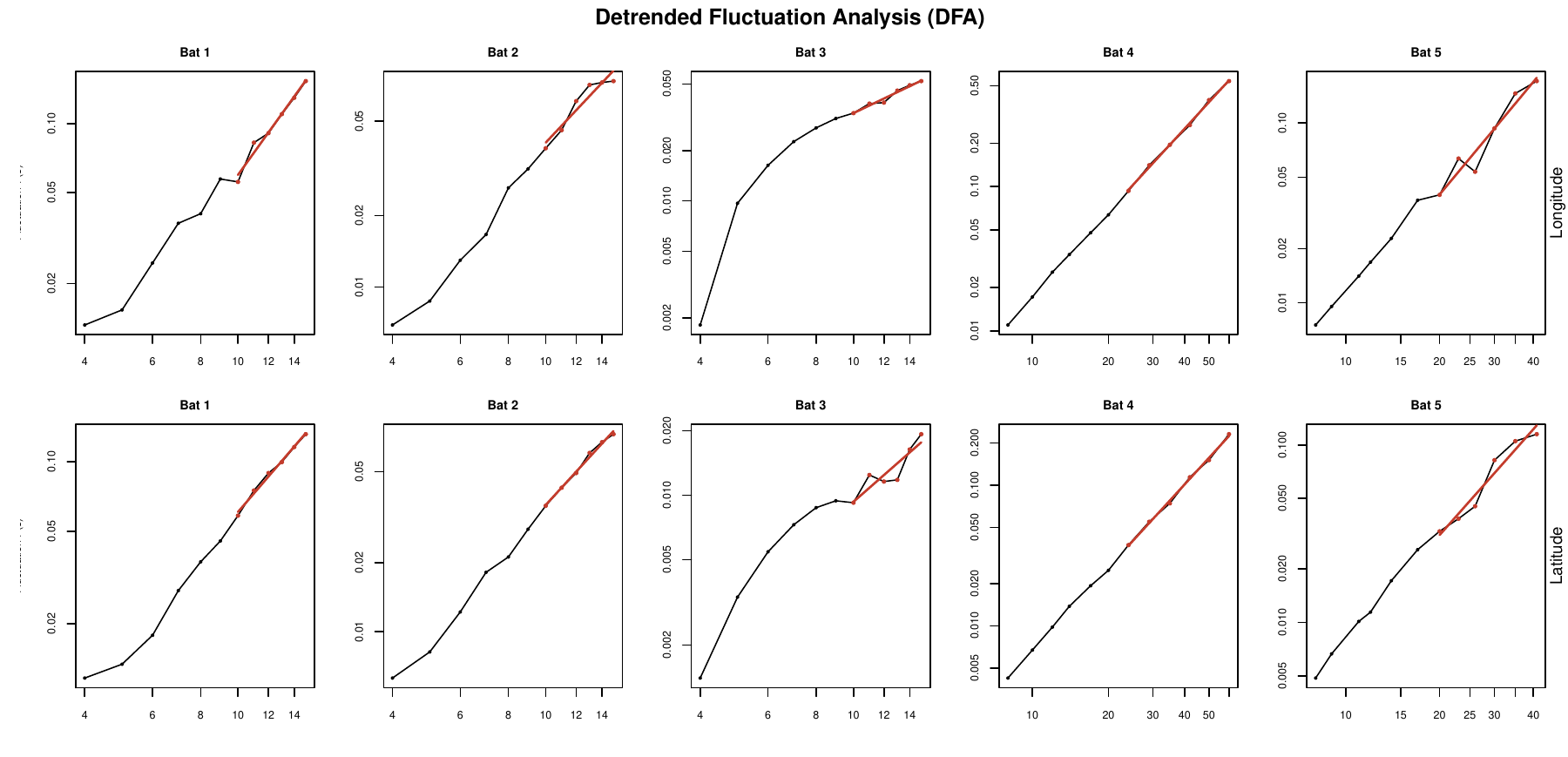}
\caption{Detrended fluctuation analysis for the telemetry series. For each
bat and coordinate, the black curve represents the empirical fluctuation
function \(F(s)\) as a function of the window size \(s\), both on logarithmic
scales. The red points and red line indicate the large-scale range used in
the regression for estimating the DFA exponent.}
\label{fig:dfa_curves_appendix}
\end{figure}

\begin{figure}[H]
\centering
\includegraphics[width=0.7\textwidth]{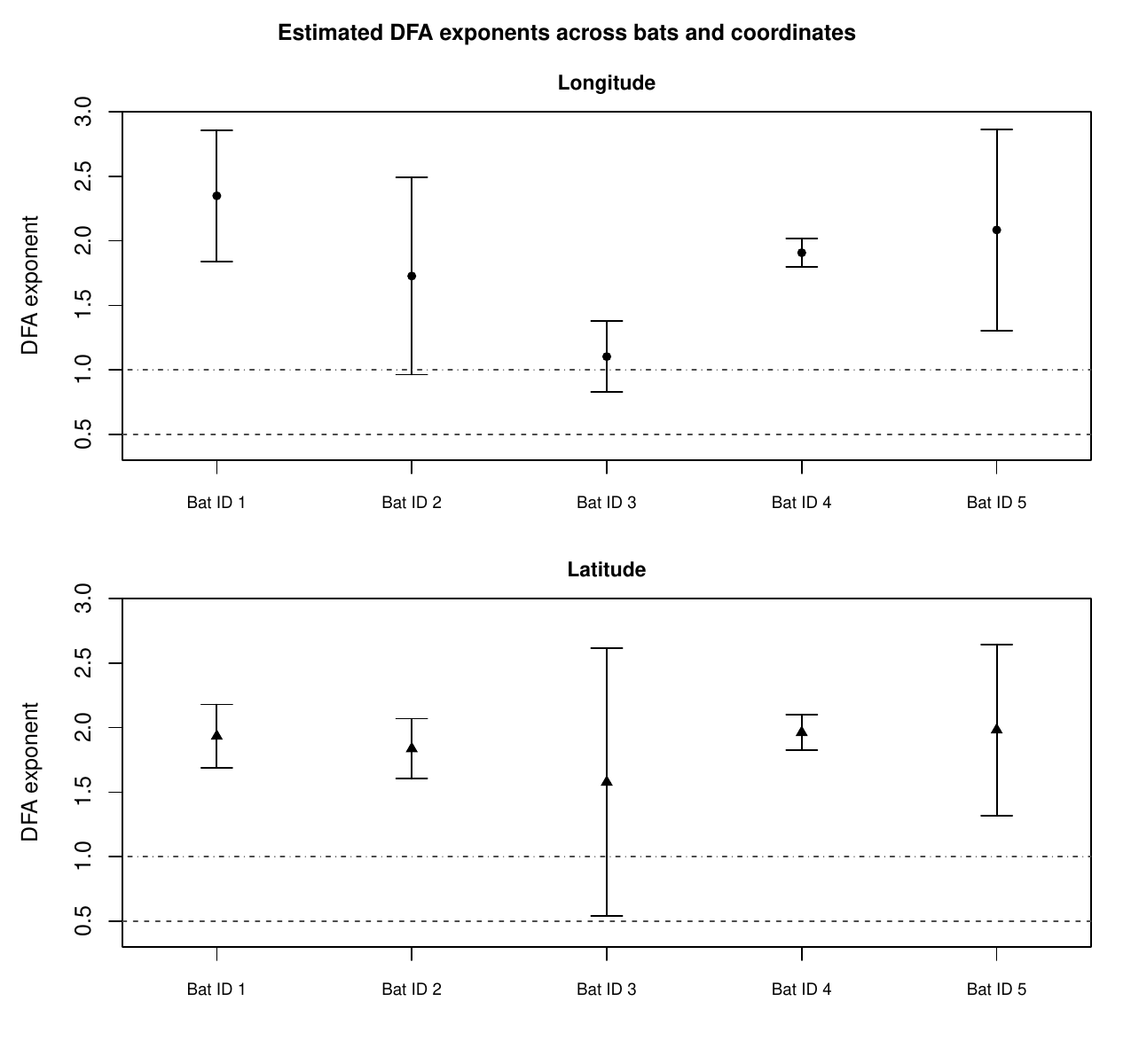}
\caption{Estimated DFA exponents across bats and coordinates. The dashed
horizontal line corresponds to \(\alpha=0.5\), a classical DFA reference
value, while the dot-dash line marks \(\alpha=1\). All point estimates lie
above \(1\); these finite-range scaling estimates are interpreted jointly
with the non-stationarity diagnostics reported above.}
\label{fig:dfa_alpha_appendix}
\end{figure}

Figure~\ref{fig:dfa_curves_appendix} displays the fluctuation functions
\(F(s)\) against the window size \(s\) for all bats and both coordinates,
on logarithmic scales. In each panel, the black curve represents the
empirical fluctuation function, while the red points and the red line
indicate the large-scale range used in the log--log regression from which
the DFA exponent was estimated. Figure~\ref{fig:dfa_alpha_appendix}
summarizes the resulting estimates \(\widehat\alpha\) for all bats and
both coordinates, together with the reference lines \(\alpha=0.5\) and
\(\alpha=1\). The estimated exponents, their regression-based confidence
intervals, and the corresponding \(R^2\) values are reported in
Table~\ref{tab:dfa_summary_appendix}.

Over the selected scaling ranges, the estimated DFA exponents exceed one
for all ten telemetry series, although some of the corresponding confidence
intervals include one. These results provide complementary empirical
evidence of dependence over large temporal scales, to be interpreted
together with the non-stationarity diagnostics above. A more detailed
inspection of Table~\ref{tab:dfa_summary_appendix} shows that Bat ID~4 has
the clearest scaling behavior in both coordinates, with estimated exponents
close to \(2\) and \(R^2=0.998\). The least stable linear fit occurs for
the latitude coordinate of Bat ID~\#3, where the confidence interval is
wider and the \(R^2\) value is \(0.816\). Even in this case, the point
estimate remains above one.

\begin{table}[H]
\centering
\scriptsize
\renewcommand{\arraystretch}{1.15}
\setlength{\tabcolsep}{4pt}
\begin{tabular}{|c|c|c|c|c|}
\hline
\textbf{Bat ID} & \textbf{Longitude: \(\widehat\alpha\) [95\% CI]} & \textbf{\(R^2\)} & \textbf{Latitude: \(\widehat\alpha\) [95\% CI]} & \textbf{\(R^2\)}\\
\hline
1 & 2.349 [1.839, 2.858] & 0.976 & 1.933 [1.688, 2.179] & 0.992 \\
2 & 1.727 [0.963, 2.490] & 0.908 & 1.836 [1.603, 2.069] & 0.992 \\
3 & 1.103 [0.830, 1.376] & 0.969 & 1.577 [0.539, 2.615] & 0.816 \\
4 & 1.907 [1.799, 2.016] & 0.998 & 1.962 [1.827, 2.098] & 0.998 \\
5 & 2.084 [1.305, 2.863] & 0.932 & 1.982 [1.320, 2.644] & 0.945 \\
\hline
\end{tabular}
\caption{Summary of the DFA exponents estimated from the telemetry series,
together with regression-based 95\% confidence intervals and coefficients of
determination \(R^2\).}
\label{tab:dfa_summary_appendix}
\end{table}

\paragraph{Scope of the empirical diagnostics}
The diagnostics in this appendix are not designed to identify the
asymptotic order of a covariance function. To state this precisely, let
\(Z=(Z_t)_{t\geq0}\) denote a candidate univariate position process and let
\(K_Z(r,t):=\operatorname{Cov}(Z_r,Z_t)\). The DFA fluctuation function
\(F(s)\), computed from an observed coordinate series, summarizes empirical
scaling across window sizes \(s\), but it is not an estimator of
\(K_Z(r,u+T)\). Therefore, without an additional model-specific relationship
between these quantities, an empirical relation of the form
\(F(s)\approx C s^\alpha\) over a selected finite range does not determine
whether \(K_Z(r,u+T)\) has logarithmic or polynomial asymptotic growth as
\(T\to\infty\).

Accordingly, the rolling summaries, sample autocorrelations, stationarity
tests, and DFA results are used as finite-sample diagnostics of
nonstationarity and dependence across temporal scales. The two asymptotic
regimes are instead represented by the complete stochastic models included
in the parametric comparison. Specifically, for observation times
\(0<t_1<\cdots<t_n\), let
\(Z^{(m,\theta)}=(Z_t^{(m,\theta)})_{t\geq0}\) denote the Gaussian process
associated with candidate model \(m\) and parameter vector \(\theta\). This
model specifies the covariance matrix
\[
\Sigma_{m,\theta}
:=
\left[
\operatorname{Cov}
\left(
Z_{t_i}^{(m,\theta)},
Z_{t_j}^{(m,\theta)}
\right)
\right]_{i,j=1}^{n}.
\]
For each bat and coordinate, \(\theta\) is estimated by Gaussian maximum
likelihood whenever a finite maximizer exists. For the flat, increasing
profile likelihoods under \(\mu_H\), the plateau representative \(\beta^*\)
defined in the main text is used instead. The candidate models are compared
using AIC evaluated at the resulting estimates. In particular,
the proposed process \(\zeta\) incorporates the logarithmic-order covariance
growth established in Proposition~1(i). The fitted fractional
Ornstein--Uhlenbeck position models \(\mu_H\), all of which satisfy
\(\widehat H>1/2\), incorporate polynomial-order covariance growth. The
integrated Ornstein--Uhlenbeck model \(\mu_{1/2}\) is treated separately.
Fractional Brownian motion represents the stationary-increment alternative;
its fitted Hurst parameters also satisfy \(\widehat H>1/2\). The stationary
Confluent and Stein covariance models complete the comparison. The likelihood
and AIC results therefore provide a relative finite-sample comparison among
complete models with different stationarity, covariance-growth, and
dependence structures.

\section{Auxiliary Results}

We first identify the temporal reductions of the space--time covariance
in Equation~(11) of \citet{Stein2005} that are compatible with the LRD
criterion adopted here, in which the normalization exponent may be any
\(b\in\mathbb R\). We consider the fixed-site temporal marginal because
longitude and latitude are modeled separately as scalar responses indexed
by time, rather than as observations of a space--time random field at
external spatial sites. Accordingly, the spatial argument in
Equation~(11) is held fixed and is not identified with the observed
differences between coordinate values.

To distinguish it from the process \(\zeta\), denote the temporal damping parameter in \citet{Stein2005} by \(\rho\), and set
\[
q(h):=\nu+\rho|h|,
\qquad
\mathcal{M}_{q}(x):=x^{q}\mathcal{K}_{q}(x),
\]
where \(\nu>0\), \(\rho\geq0\), and \(\mathcal{K}_{q}\) denotes the modified Bessel function of the second kind of order \(q\). After normalizing the variance to \(\sigma^2\), the fixed-site temporal marginal of Equation~(11) is
\begin{equation}\label{eq:stein-fixed-site-general}
R_{\varepsilon}(h)
=
\sigma^2
\frac{\Gamma\!\left(\nu+\frac d2\right)}{\Gamma(\nu)}
\frac{
\mathcal{M}_{q(h)}
\!\left(a|\varepsilon|\,|h|\right)
}{
2^{q(h)-1}
\Gamma\!\left(q(h)+\frac d2\right)
},
\end{equation}
where \(a>0\), \(d\in\mathbb N\), and \(\varepsilon\) is the advection parameter.

\begin{prop}
\label{prop:stein-lrd-classification}
The stationary Gaussian process with covariance
\eqref{eq:stein-fixed-site-general} exhibits LRD if and only if
\[
\varepsilon=0
\qquad\text{and}\qquad
\rho>0.
\]
In this case,
\begin{equation}\label{eq:stein-gamma-ratio}
R_{d,\nu,\rho}(h)
=
\sigma^2
\frac{\Gamma\!\left(\nu+\frac d2\right)}{\Gamma(\nu)}
\frac{\Gamma(\nu+\rho|h|)}
{\Gamma\!\left(\nu+\rho|h|+\frac d2\right)}.
\end{equation}
Writing \(p=d/2\), one has
\begin{equation}\label{eq:stein-gamma-tail}
R_{d,\nu,\rho}(h)
\sim
c_{d,\nu,\rho}|h|^{-p},
\qquad
|h|\to\infty,
\end{equation}
where
\[
c_{d,\nu,\rho}
:=
\sigma^2
\frac{\Gamma(\nu+p)}
{\Gamma(\nu)\rho^p}.
\]
Moreover, for every \(r<v\) and \(s<t\),
\begin{equation}\label{eq:stein-increment-lrd}
\lim_{T\to\infty}
T^{p+2}
\operatorname{Cov}
\left(
X_v-X_r,
X_{t+T}-X_{s+T}
\right)
=
-c_{d,\nu,\rho}
p(p+1)(v-r)(t-s),
\end{equation}
and the limit is finite and nonzero.
\end{prop}

\begin{proof}
Since
\[
\mathcal{M}_{q}(0)
=
2^{q-1}\Gamma(q),
\qquad q>0,
\]
setting \(\varepsilon=0\) in
\eqref{eq:stein-fixed-site-general} yields
\eqref{eq:stein-gamma-ratio}.

Assume first that \(\rho>0\), and define
\[
p:=\frac d2,
\qquad
A_{\nu,p}
:=
\sigma^2\frac{\Gamma(\nu+p)}{\Gamma(\nu)},
\qquad
g_p(x)
:=
\frac{\Gamma(x)}{\Gamma(x+p)}.
\]
Then, for \(h>0\),
\[
R_{d,\nu,\rho}(h)
=
A_{\nu,p}g_p(\nu+\rho h).
\]
The standard asymptotic expansion for a quotient of gamma functions gives
\[
g_p(x)
=
x^{-p}\left(1+O(x^{-1})\right),
\qquad x\to\infty.
\]
Consequently,
\[
R_{d,\nu,\rho}(h)
=
A_{\nu,p}
(\nu+\rho h)^{-p}
\left(1+O(h^{-1})\right),
\]
which proves \eqref{eq:stein-gamma-tail}.

Let
\[
\psi(x)
:=
\frac{\Gamma'(x)}{\Gamma(x)}
=
\frac{d}{dx}\log\Gamma(x)
\]
denote the digamma function, and let
\[
\psi_1(x)
:=
\psi'(x)
=
\frac{d^2}{dx^2}\log\Gamma(x)
\]
denote the trigamma function. Since
\[
\frac{g_p'(x)}{g_p(x)}
=
\psi(x)-\psi(x+p)
=
-\frac{p}{x}+O(x^{-2})
\]
and
\[
\psi_1(x)-\psi_1(x+p)
=
\frac{p}{x^2}+O(x^{-3}),
\]
we obtain
\[
\begin{aligned}
\frac{g_p''(x)}{g_p(x)}
&=
\left(
\psi(x)-\psi(x+p)
\right)^2
+
\psi_1(x)-\psi_1(x+p)
\\
&=
\frac{p(p+1)}{x^2}
+
O(x^{-3}).
\end{aligned}
\]
Thus,
\[
g_p''(x)
=
p(p+1)x^{-p-2}
\left(1+O(x^{-1})\right).
\]
Since
\[
R_{d,\nu,\rho}''(h)
=
A_{\nu,p}\rho^2
g_p''(\nu+\rho h),
\]
it follows that
\begin{equation}\label{eq:stein-second-derivative}
R_{d,\nu,\rho}''(h)
\sim
c_{d,\nu,\rho}
p(p+1)h^{-p-2}.
\end{equation}

By stationarity,
\[
\begin{aligned}
D_T
:={}&
\operatorname{Cov}
\left(
X_v-X_r,
X_{t+T}-X_{s+T}
\right)
\\
={}&
R_{d,\nu,\rho}(T+t-v)
-
R_{d,\nu,\rho}(T+s-v)
\\
&-
R_{d,\nu,\rho}(T+t-r)
+
R_{d,\nu,\rho}(T+s-r).
\end{aligned}
\]
Since \(R_{d,\nu,\rho}\) is twice continuously differentiable on
\((0,\infty)\), two applications of the fundamental theorem of calculus give
\begin{equation}\label{eq:increment-double-integral}
D_T
=
-\int_r^v\int_s^t
R_{d,\nu,\rho}''(T+y-x)
\,dy\,dx.
\end{equation}
The asymptotic relation
\eqref{eq:stein-second-derivative} holds uniformly when its argument differs from \(T\) by a quantity in a fixed compact interval. Hence,
\[
T^{p+2}
R_{d,\nu,\rho}''(T+y-x)
\longrightarrow
c_{d,\nu,\rho}p(p+1)
\]
uniformly for
\[
(x,y)\in[r,v]\times[s,t].
\]
Passing to the limit in
\eqref{eq:increment-double-integral} proves
\eqref{eq:stein-increment-lrd}. Thus, the LRD criterion is satisfied with
\[
b=p+2=\frac d2+2.
\]

It remains to exclude the other parameter regimes. If
\[
\varepsilon=0
\qquad\text{and}\qquad
\rho=0,
\]
then \(q(h)\equiv\nu\) and
\[
R_{\varepsilon}(h)\equiv\sigma^2.
\]
Therefore,
\[
\operatorname{Cov}
\left(
X_v-X_r,
X_{t+T}-X_{s+T}
\right)
=0
\]
for every \(T\), and consequently
\[
T^b
\operatorname{Cov}
\left(
X_v-X_r,
X_{t+T}-X_{s+T}
\right)
=0
\]
for every \(b\in\mathbb R\).

If \(\varepsilon\neq0\) and \(\rho=0\), the order of the Bessel function is fixed and
\[
\mathcal{K}_{\nu}(x)
\sim
\sqrt{\frac{\pi}{2x}}e^{-x},
\qquad x\to\infty.
\]
Consequently,
\[
R_{\varepsilon}(h)
=
O\!\left(
|h|^{\nu-\frac12}
e^{-a|\varepsilon||h|}
\right).
\]
The covariance of the distant increments is a linear combination of four terms with the same exponential decay. Hence,
\[
T^b
\operatorname{Cov}
\left(
X_v-X_r,
X_{t+T}-X_{s+T}
\right)
\longrightarrow0,
\qquad T\to\infty,
\]
for every \(b\in\mathbb R\). For \(b\geq0\), this follows because exponential decay dominates every polynomial power, whereas for \(b<0\) it follows from the boundedness of the covariance and \(T^b\to0\).

Finally, if \(\varepsilon\neq0\) and \(\rho>0\),
\citet[p.~315]{Stein2005} shows that the corresponding fixed-site temporal covariance decays exponentially. The same argument therefore gives
\[
T^b
\operatorname{Cov}
\left(
X_v-X_r,
X_{t+T}-X_{s+T}
\right)
\longrightarrow0,
\qquad T\to\infty,
\]
for every \(b\in\mathbb R\). Hence none of the remaining parameter regimes satisfies the adopted LRD criterion.
\end{proof}

\begin{cor}
\label{cor:stein-power-closure}
Let \(d\in\mathbb N\), with \(\mathbb N=\{1,2,\ldots\}\), and set
\[
p:=\frac d2.
\]
For \(\nu>0\), \(\rho>0\), and \(\sigma^2>0\), consider the
LRD-compatible Stein temporal covariance
\[
R_{d,\nu,\rho}(h)
:=
\sigma^2
\frac{\Gamma(\nu+p)}{\Gamma(\nu)}
\frac{\Gamma(\nu+\rho|h|)}
{\Gamma(\nu+\rho|h|+p)},
\qquad h\in\mathbb R.
\]
For \(\lambda>0\), define
\begin{equation}\label{eq:Sd-covariance}
R^{S,d}_{\sigma^2,\lambda}(h)
:=
\frac{\sigma^2}{(1+\lambda|h|)^p},
\qquad h\in\mathbb R.
\end{equation}
Then the following statements hold.

\begin{enumerate}
\item[(i)]
The function \(R^{S,d}_{\sigma^2,\lambda}\) is a valid stationary
covariance function on \(\mathbb R\). Accordingly, let
\(S^{(d)}_{\sigma^2,\lambda}\) denote a centered stationary Gaussian
process with covariance \eqref{eq:Sd-covariance}.

\item[(ii)]
The process \(S^{(d)}_{\sigma^2,\lambda}\) exhibits LRD. More precisely,
for every \(r<v\) and \(s<t\),
\begin{equation}\label{eq:Sd-lrd-limit}
\begin{aligned}
&\lim_{T\to\infty}T^{p+2}
\operatorname{Cov}\left(
S^{(d)}_{\sigma^2,\lambda}(v)
-
S^{(d)}_{\sigma^2,\lambda}(r),
\right.\\[-1mm]
&\hspace{4.3cm}\left.
S^{(d)}_{\sigma^2,\lambda}(t+T)
-
S^{(d)}_{\sigma^2,\lambda}(s+T)
\right)\\
&\qquad
=
-\sigma^2\lambda^{-p}p(p+1)(v-r)(t-s).
\end{aligned}
\end{equation}

\item[(iii)]
Under the parametrization
\[
\rho=\nu\lambda,
\]
the Stein family has two nondegenerate boundary closures. As
\(\nu\to\infty\),
\begin{equation}\label{eq:stein-to-Sd}
R_{d,\nu,\nu\lambda}(h)
\longrightarrow
R^{S,d}_{\sigma^2,\lambda}(h)
=
\frac{\sigma^2}{(1+\lambda|h|)^p},
\qquad h\in\mathbb R,
\end{equation}
whereas, as \(\nu\downarrow0\),
\begin{equation}\label{eq:stein-small-nu-limit}
R_{d,\nu,\nu\lambda}(h)
\longrightarrow
R^{S,2}_{\sigma^2,\lambda}(h)
=
\frac{\sigma^2}{1+\lambda|h|},
\qquad h\in\mathbb R.
\end{equation}
The convergence in \eqref{eq:stein-small-nu-limit} is uniform on every
compact subset of \(\mathbb R\). In particular, the small-\(\nu\)
closure is independent of \(d\). For convenience, write
\[
S_{\sigma^2,\lambda}
:=
S^{(2)}_{\sigma^2,\lambda}.
\]

\item[(iv)]
The small-\(\nu\) limiting process \(S_{\sigma^2,\lambda}\) exhibits
LRD and, for every \(r<v\) and \(s<t\),
\begin{equation}\label{eq:stein-small-nu-lrd}
\begin{aligned}
&\lim_{T\to\infty}T^3
\operatorname{Cov}\left(
S_{\sigma^2,\lambda}(v)
-
S_{\sigma^2,\lambda}(r),
\right.\\[-1mm]
&\hspace{4.3cm}\left.
S_{\sigma^2,\lambda}(t+T)
-
S_{\sigma^2,\lambda}(s+T)
\right)\\
&\qquad
=
-\frac{2\sigma^2}{\lambda}(v-r)(t-s).
\end{aligned}
\end{equation}

\item[(v)]
The scaling \(\rho=\nu\lambda\) is essential for the nondegenerate
small-\(\nu\) closure. If instead \(\rho>0\) is kept fixed, then
\begin{equation}\label{eq:stein-fixed-rho-small-nu}
R_{d,\nu,\rho}(h)
\longrightarrow
R_{\mathrm{nug}}(h)
:=
\begin{cases}
\sigma^2, & h=0,\\
0, & h\neq0,
\end{cases}
\qquad \nu\downarrow0.
\end{equation}
Thus, the fixed-\(\rho\) limit is the nugget covariance and does not
exhibit LRD.

\item[(vi)]
When \(d=2\), one has \(p=1\), and the reduction is exact for every
\(\nu>0\) and \(\rho>0\):
\begin{equation}\label{eq:stein-d2-exact}
R_{2,\nu,\rho}(h)
=
\frac{\sigma^2\nu}{\nu+\rho|h|}
=
\frac{\sigma^2}{1+(\rho/\nu)|h|}.
\end{equation}
Equivalently,
\[
R_{2,\nu,\rho}
=
R^{S,2}_{\sigma^2,\rho/\nu}.
\]
In particular, if \(\rho=\nu\lambda\), then
\[
R_{2,\nu,\nu\lambda}(h)
=
\frac{\sigma^2}{1+\lambda|h|}
\]
for every \(\nu>0\), so the small- and large-\(\nu\) closures coincide.

\item[(vii)]
When \(d=3\), one has \(p=3/2\), and
\begin{equation}\label{eq:stein-d3-covariance}
R_{3,\nu,\rho}(h)
=
\sigma^2
\frac{\Gamma\!\left(\nu+\frac32\right)}{\Gamma(\nu)}
\frac{\Gamma(\nu+\rho|h|)}
{\Gamma\!\left(\nu+\rho|h|+\frac32\right)}.
\end{equation}
Under \(\rho=\nu\lambda\), its two boundary closures are
\begin{equation}\label{eq:d3-small-nu}
R_{3,\nu,\nu\lambda}(h)
\longrightarrow
\frac{\sigma^2}{1+\lambda|h|},
\qquad \nu\downarrow0,
\end{equation}
and
\begin{equation}\label{eq:d3-large-nu}
R_{3,\nu,\nu\lambda}(h)
\longrightarrow
\frac{\sigma^2}{(1+\lambda|h|)^{3/2}},
\qquad \nu\to\infty.
\end{equation}
Thus, the closure of the \(d=3\) family contains both
\(S_{\sigma^2,\lambda}\) and \(S^{(3)}_{\sigma^2,\lambda}\), where
\begin{equation}\label{eq:S3-covariance}
R^{S,3}_{\sigma^2,\lambda}(h)
=
\frac{\sigma^2}{(1+\lambda|h|)^{3/2}}.
\end{equation}
Their covariance tails are proportional to \(|h|^{-1}\) and
\(|h|^{-3/2}\), respectively. Consequently, their LRD normalizations
are \(T^3\) and \(T^{7/2}\). In particular,
\begin{equation}\label{eq:S3-lrd-limit}
\begin{aligned}
&\lim_{T\to\infty}T^{7/2}
\operatorname{Cov}\left(
S^{(3)}_{\sigma^2,\lambda}(v)
-
S^{(3)}_{\sigma^2,\lambda}(r),
\right.\\[-1mm]
&\hspace{4.3cm}\left.
S^{(3)}_{\sigma^2,\lambda}(t+T)
-
S^{(3)}_{\sigma^2,\lambda}(s+T)
\right)\\
&\qquad
=
-\frac{15}{4}\sigma^2\lambda^{-3/2}
(v-r)(t-s).
\end{aligned}
\end{equation}
\end{enumerate}
\end{cor}

\begin{proof}
Since \(p=d/2>0\), the gamma-integral representation gives
\[
(1+\lambda|h|)^{-p}
=
\frac{1}{\Gamma(p)}
\int_0^\infty
u^{p-1}e^{-u}e^{-\lambda u|h|}
\,du.
\]
For every \(u>0\), the function
\[
h\longmapsto e^{-\lambda u|h|}
\]
is positive definite on \(\mathbb R\). Hence
\eqref{eq:Sd-covariance} is a positive mixture of covariance functions,
which proves~(i).

For \(h>0\),
\[
\frac{d^2}{dh^2}
R^{S,d}_{\sigma^2,\lambda}(h)
=
\sigma^2p(p+1)\lambda^2
(1+\lambda h)^{-p-2}
\sim
\sigma^2p(p+1)\lambda^{-p}h^{-p-2}.
\]
Applying the double-integral identity
\eqref{eq:increment-double-integral}, together with the uniformity of
this asymptotic relation under bounded shifts of its argument, yields
\eqref{eq:Sd-lrd-limit}. Therefore, the adopted LRD criterion holds
with
\[
b=p+2=\frac d2+2,
\]
proving~(ii).

For the large-\(\nu\) limit, substituting
\(\rho=\nu\lambda\) gives
\[
R_{d,\nu,\nu\lambda}(h)
=
\sigma^2
\frac{\Gamma(\nu+p)}{\Gamma(\nu)}
\frac{
\Gamma\!\left(\nu(1+\lambda|h|)\right)
}{
\Gamma\!\left(\nu(1+\lambda|h|)+p\right)
}.
\]
For fixed \(p>0\) and \(h\in\mathbb R\), the gamma-ratio asymptotics
yield
\[
\frac{\Gamma(\nu+p)}{\Gamma(\nu)}
\sim\nu^p
\]
and
\[
\frac{
\Gamma\!\left(\nu(1+\lambda|h|)\right)
}{
\Gamma\!\left(\nu(1+\lambda|h|)+p\right)
}
\sim
\nu^{-p}(1+\lambda|h|)^{-p}.
\]
Their product proves \eqref{eq:stein-to-Sd}.

For the small-\(\nu\) limit, let
\[
B(x,p)
:=
\frac{\Gamma(x)\Gamma(p)}{\Gamma(x+p)},
\qquad x,p>0.
\]
Then
\[
\frac{R_{d,\nu,\nu\lambda}(h)}{\sigma^2}
=
\frac{
B\!\left(\nu(1+\lambda|h|),p\right)
}{
B(\nu,p)
}.
\]
Define
\[
F_p(x):=xB(x,p).
\]
Since
\[
x\Gamma(x)\longrightarrow1
\qquad\text{and}\qquad
\Gamma(x+p)\longrightarrow\Gamma(p)
\]
as \(x\downarrow0\), one has
\[
F_p(x)\longrightarrow1.
\]
Writing \(c(h):=1+\lambda|h|\), we obtain
\[
\frac{
B\!\left(\nu c(h),p\right)
}{
B(\nu,p)
}
=
\frac{F_p\!\left(\nu c(h)\right)}{F_p(\nu)}
\frac{1}{c(h)}
\longrightarrow
\frac{1}{1+\lambda|h|},
\]
which proves \eqref{eq:stein-small-nu-limit}.

To verify local uniform convergence, fix \(H>0\). For \(|h|\leq H\),
\[
1\leq c(h)\leq1+\lambda H.
\]
Hence
\[
0<\nu c(h)\leq\nu(1+\lambda H),
\]
and \(F_p(x)\to1\) uniformly for
\(0<x\leq\nu(1+\lambda H)\) as \(\nu\downarrow0\). Therefore,
\eqref{eq:stein-small-nu-limit} is uniform on \([-H,H]\). Applying
\eqref{eq:Sd-lrd-limit} with \(d=2\), and hence \(p=1\), gives
\eqref{eq:stein-small-nu-lrd}, proving~(iii) and~(iv).

Suppose now that \(\rho>0\) remains fixed. For \(h\neq0\),
\[
\frac{\Gamma(\nu+p)}{\Gamma(\nu)}
\sim
\Gamma(p)\nu,
\qquad \nu\downarrow0,
\]
whereas
\[
\frac{\Gamma(\nu+\rho|h|)}
{\Gamma(\nu+\rho|h|+p)}
\longrightarrow
\frac{\Gamma(\rho|h|)}
{\Gamma(\rho|h|+p)}.
\]
Thus,
\[
R_{d,\nu,\rho}(h)\longrightarrow0,
\qquad h\neq0.
\]
At \(h=0\),
\[
R_{d,\nu,\rho}(0)=\sigma^2
\]
for every \(\nu>0\), proving
\eqref{eq:stein-fixed-rho-small-nu}. For the resulting nugget
covariance, the covariance between the two distant increments is zero
for all sufficiently large \(T\); hence it does not satisfy the adopted
LRD criterion. This proves~(v).

If \(d=2\), then \(p=1\), and the recurrence relation
\[
\Gamma(x+1)=x\Gamma(x)
\]
applied to \(R_{d,\nu,\rho}\) gives
\eqref{eq:stein-d2-exact}, proving~(vi).

Finally, if \(d=3\), then
\[
p=\frac32,
\qquad
p+2=\frac72,
\qquad
p(p+1)
=
\frac32\frac52
=
\frac{15}{4}.
\]
Substitution into the general Stein covariance gives
\eqref{eq:stein-d3-covariance}. Equations
\eqref{eq:d3-small-nu} and \eqref{eq:d3-large-nu} follow from
\eqref{eq:stein-small-nu-limit} and \eqref{eq:stein-to-Sd},
respectively. Finally, substituting \(p=3/2\) into
\eqref{eq:Sd-covariance} and \eqref{eq:Sd-lrd-limit} yields
\eqref{eq:S3-covariance} and \eqref{eq:S3-lrd-limit}.
\end{proof}

\subsection{Matérn-type local behavior and LRD in the confluent
hypergeometric class}

Let \(U(a,b,z)\) denote Tricomi's confluent hypergeometric function.

\begin{prop}
\label{prop:CH-local-LRD}
Let \(\sigma^2>0\), \(\eta>0\), \(\alpha>0\), and \(\beta>0\), and define
\begin{equation}\label{eq:CH-covariance}
R_A(h)
:=
\sigma^2
\frac{\Gamma(\eta+\alpha)}{\Gamma(\eta)}
U\!\left(
\alpha,1-\eta,
\eta\left(\frac{|h|}{\beta}\right)^2
\right),
\qquad h\in\mathbb R.
\end{equation}
Then \(R_A(0)=\sigma^2\), and the following statements hold.

\begin{enumerate}
\item[(i)]
The function \(R_A\) is a valid stationary covariance function on
\(\mathbb R\). Accordingly, let
\(A_{\sigma^2,\eta,\alpha,\beta}\) denote a centered stationary Gaussian
process with covariance \eqref{eq:CH-covariance}. Moreover, \(R_A\)
belongs to the confluent hypergeometric class obtained through the
Matérn scale-mixture construction of \citet{MaBhadra2023}.

\item[(ii)]
For every integer \(m\geq0\), the process
\(A_{\sigma^2,\eta,\alpha,\beta}\) is \(m\) times mean-square
differentiable if and only if
\[
\eta>m.
\]
Thus, its local mean-square regularity is determined by \(\eta\), as in
the Matérn class.

\item[(iii)]
The covariance has the polynomial tail
\begin{equation}\label{eq:CH-tail}
R_A(h)
\sim
c_A|h|^{-2\alpha},
\qquad |h|\to\infty,
\end{equation}
where
\[
c_A
:=
\sigma^2
\frac{\Gamma(\eta+\alpha)}{\Gamma(\eta)}
\eta^{-\alpha}\beta^{2\alpha}.
\]

\item[(iv)]
Under the LRD criterion adopted here, the process exhibits LRD. More
precisely, for every \(r<v\) and \(s<t\),
\begin{equation}\label{eq:CH-lrd-limit}
\begin{aligned}
&\lim_{T\to\infty}
T^{2\alpha+2}
\operatorname{Cov}
\left(
A_{\sigma^2,\eta,\alpha,\beta}(v)
-
A_{\sigma^2,\eta,\alpha,\beta}(r),
\right.
\\[-1mm]
&\hspace{3.7cm}\left.
A_{\sigma^2,\eta,\alpha,\beta}(t+T)
-
A_{\sigma^2,\eta,\alpha,\beta}(s+T)
\right)
\\
&\qquad
=
-c_A(2\alpha)(2\alpha+1)(v-r)(t-s).
\end{aligned}
\end{equation}
\end{enumerate}
\end{prop}

\begin{proof}
Set
\[
c:=\frac{\eta}{\beta^2},
\qquad
B:=\sigma^2\frac{\Gamma(\eta+\alpha)}{\Gamma(\eta)}.
\]
For \(z>0\), Tricomi's function admits the integral representation
\begin{equation}\label{eq:CH-U-integral}
U(\alpha,1-\eta,z)
=
\frac{1}{\Gamma(\alpha)}
\int_0^\infty
e^{-zu}u^{\alpha-1}(1+u)^{-\alpha-\eta}
\,du.
\end{equation}
Consequently,
\[
\frac{R_A(h)}{\sigma^2}
=
\int_0^\infty e^{-cuh^2}\pi_{\alpha,\eta}(u)\,du,
\]
where
\[
\pi_{\alpha,\eta}(u)
:=
\frac{\Gamma(\alpha+\eta)}
{\Gamma(\alpha)\Gamma(\eta)}
u^{\alpha-1}(1+u)^{-\alpha-\eta},
\qquad u>0,
\]
is a probability density. Since \(h\mapsto e^{-cuh^2}\) is positive
definite for every \(u>0\), \(R_A\) is a positive mixture of covariance
functions and is therefore a valid covariance. Setting \(h=0\) gives
\(R_A(0)=\sigma^2\), proving~(i).

For an integer \(m\geq0\), the \(m\)-th moment of the mixing variable is
\[
\int_0^\infty u^m\pi_{\alpha,\eta}(u)\,du
=
\frac{\Gamma(\alpha+m)\Gamma(\eta-m)}
{\Gamma(\alpha)\Gamma(\eta)}
\]
when \(\eta>m\), and is infinite when \(\eta\leq m\). Since
\[
\left.
\frac{d^{2m}}{dh^{2m}}e^{-cuh^2}
\right|_{h=0}
=
(-1)^m\frac{(2m)!}{m!}(cu)^m,
\]
the standard covariance criterion for mean-square differentiability
shows that the process is \(m\) times mean-square differentiable
exactly when \(\eta>m\). This proves~(ii).

For \(k=0,1,2\), differentiation under the integral sign in
\eqref{eq:CH-U-integral} gives
\[
\frac{d^k}{dz^k}U(\alpha,1-\eta,z)
=
\frac{(-1)^k}{\Gamma(\alpha)}
\int_0^\infty
e^{-zu}u^{\alpha+k-1}(1+u)^{-\alpha-\eta}
\,du.
\]
After the change of variables \(w=zu\),
\[
\frac{d^k}{dz^k}U(\alpha,1-\eta,z)
=
\frac{(-1)^kz^{-\alpha-k}}{\Gamma(\alpha)}
\int_0^\infty
e^{-w}w^{\alpha+k-1}
\left(1+\frac{w}{z}\right)^{-\alpha-\eta}
\,dw.
\]
Dominated convergence therefore yields
\begin{equation}\label{eq:CH-U-derivatives}
\frac{d^k}{dz^k}U(\alpha,1-\eta,z)
\sim
(-1)^k(\alpha)_kz^{-\alpha-k},
\qquad z\to\infty,
\end{equation}
where
\[
(\alpha)_k
:=
\frac{\Gamma(\alpha+k)}{\Gamma(\alpha)}.
\]
Taking \(k=0\) and \(z=ch^2\) gives
\[
R_A(h)
=
B\,U(\alpha,1-\eta,ch^2)
\sim
Bc^{-\alpha}|h|^{-2\alpha}.
\]
Since
\[
Bc^{-\alpha}
=
\sigma^2
\frac{\Gamma(\eta+\alpha)}{\Gamma(\eta)}
\eta^{-\alpha}\beta^{2\alpha}
=
c_A,
\]
this proves \eqref{eq:CH-tail}.

Let
\[
F(z):=U(\alpha,1-\eta,z).
\]
For \(h>0\),
\[
R_A''(h)
=
B\left[
2cF'(ch^2)+4c^2h^2F''(ch^2)
\right].
\]
By \eqref{eq:CH-U-derivatives},
\[
F'(z)\sim-\alpha z^{-\alpha-1},
\qquad
F''(z)\sim\alpha(\alpha+1)z^{-\alpha-2},
\]
and hence
\[
R_A''(h)
\sim
c_A(2\alpha)(2\alpha+1)h^{-2\alpha-2}.
\]
By stationarity,
\[
\begin{aligned}
&\operatorname{Cov}
\left(
A(v)-A(r),
A(t+T)-A(s+T)
\right)
\\
&\qquad
=
-\int_r^v\int_s^t
R_A''(T+y-x)\,dy\,dx,
\end{aligned}
\]
where \(A=A_{\sigma^2,\eta,\alpha,\beta}\). The preceding asymptotic
relation is uniform when its argument differs from \(T\) by a bounded
quantity. Multiplication by \(T^{2\alpha+2}\) and passage to the limit
therefore give \eqref{eq:CH-lrd-limit}, proving~(iv).
\end{proof}

\begin{rem}
\label{rem:CH-Matern-LRD}
The ordinary Matérn covariance with smoothness parameter \(\eta>0\) and
inverse range \(\kappa>0\) is
\[
R_{\mathrm{Mat}}(h)
=
\sigma^2
\frac{2^{1-\eta}}{\Gamma(\eta)}
(\kappa|h|)^\eta
\mathcal{K}_{\eta}(\kappa|h|),
\]
where \(\mathcal{K}_{\eta}\) is the modified Bessel function of the
second kind. Since
\[
\mathcal{K}_{\eta}(x)
\sim
\sqrt{\frac{\pi}{2x}}e^{-x},
\qquad x\to\infty,
\]
one has
\[
R_{\mathrm{Mat}}(h)
=
O\!\left(
|h|^{\eta-\frac12}e^{-\kappa|h|}
\right),
\qquad |h|\to\infty .
\]
Thus, if \(M\) is a stationary process with Matérn covariance, then for
every \(b\in\mathbb R\), \(r<v\), and \(s<t\),
\[
T^b
\operatorname{Cov}
\left(
M_v-M_r,
M_{t+T}-M_{s+T}
\right)
\longrightarrow0.
\]
Indeed, exponential decay dominates every polynomial when \(b\geq0\),
whereas for \(b<0\) the conclusion follows from boundedness of the
covariance and \(T^b\to0\). Hence the ordinary Matérn class with
\(\kappa>0\) does not exhibit LRD under the adopted criterion.

The confluent hypergeometric class separates the local and long-lag
features that the ordinary Matérn class cannot control independently.
The parameter \(\sigma^2\) controls the marginal variance, the parameter
\(\eta\) determines the local mean-square differentiability, and the
parameter \(\alpha\) determines the polynomial tail exponent \(2\alpha\)
and the LRD normalization \(T^{2\alpha+2}\). The parameter \(\beta\)
controls the temporal range and the multiplicative constant in the
long-lag asymptotics, but not either exponent. Therefore,
``Matérn-type local behavior with LRD'' means precisely that the CH class
combines Matérn local regularity, governed by \(\eta\), with polynomial
long-lag dependence, governed by \(\alpha\).
\end{rem}

\end{document}